\documentclass[final,10pt]{elsarticle}
\usepackage{amsbsy,amsfonts,amssymb,multirow}
\usepackage{color}
\usepackage{amsmath,fontenc,epsfig,subfigure,lettrine,stmaryrd,enumerate}
\usepackage[latin1]{inputenc}
\usepackage{graphics}
\usepackage{psfrag}
\usepackage[crop=pdfcrop]{pstool}
\usepackage[margin=8pt,font=small,labelfont=bf,format=plain,indention=.5cm,labelsep=quad,skip=6pt]{caption}
\usepackage{url}
\usepackage{float}
\usepackage{algorithm}
\usepackage{algorithmic}
\usepackage{fullpage}

\usepackage{soul,color,ulem}

 \usepackage{amsmath,amssymb,epsfig}
  \usepackage{subfigure}
 \usepackage{subfigmat}

  \def\beq{\begin{equation}}
\def\eeq{\end{equation}}
\def\esplit{\end{split}}
\def\beqalign{\begin{array}{rl}}
\def\eeqalign{\end{array}}
\DeclareMathOperator*{\argmin}{\arg\!\min}

\def\Abold{\mathbf{A}}
\def\Bbold{\mathbf{B}}
\def\Cbold{\mathbf{C}}

\def\Ebold{\mathbf{E}}
\def\Fbold{\mathbf{F}}
\def\Gbold{\mathbf{G}}
\def\Hbold{\mathbf{H}}
\def\Ibold{\mathbf{I}}

\def\Kbold{\mathbf{K}}

\def\Mbold{\mathbf{M}}
\def\Nbold{\mathbf{N}}

\def\Qbold{\mathbf{Q}}
\def\Rbold{\mathbf{R}}
\def\Sbold{\mathbf{S}}
\def\Tbold{\mathbf{T}}
\def\Ubold{\mathbf{U}}
\def\Vbold{\mathbf{V}}
\def\Wbold{\mathbf{W}}
\def\Xbold{\mathbf{X}}
\def\Ybold{\mathbf{Y}}
\def\Zbold{\mathbf{Z}}

\def\dbold{\mathbf{d}}

\def\fbold{\mathbf{f}}

\def\qbold{\mathbf{q}}

\def\sbold{\mathbf{s}}
\def\ubold{\mathbf{u}}

\def\wbold{\mathbf{w}}
\def\xbold{\mathbf{x}}
\def\ybold{\mathbf{y}}

\def\mubold{\boldsymbol{\mu}}
\def\nubold{\boldsymbol{\nu}}

\def\Sigmabold{\mathbf{\Sigma}}

\def\Lambdabold{\mathbf{\Lambda}}

\def\0bold{\boldsymbol{0}}
\def\0{\mathbf{\0}}

\begin{document}
\begin{frontmatter}


\title{Real-time solution of computational problems using databases of parametric linear reduced-order models with arbitrary underlying meshes}

\author[aa]{David Amsallem\corref{cor1}}
\ead{amsallem@stanford.edu}
\ead[url]{stanford.edu/${\sim}$amsallem}

\author[aa]{Radek Tezaur}
\ead{rtezaur@stanford.edu}

\author[cf]{Charbel Farhat}
\ead{cfarhat@stanford.edu}


\address[aa]{Department of Aeronautics and Astronautics, Durand Building, 496 Lomita Mall. Stanford University, Stanford, 94305-4035, USA}

\address[cf]{Department of Aeronautics and Astronautics, Institute for Computational and Mathematical Engineering, Department of Mechanical Engineering. Stanford University, Mail Code 4035, Stanford, CA 94305, U.S.A.}


\begin{abstract}
A comprehensive approach for real-time computations using a database of parameterized linear reduced-order models (ROMs) is proposed. The method proceeds by sampling offline ROMs for specific values of the parameters and interpolating online the associated reduced operators. In the offline phase, a pre-processing step transforms the reduced operators into consistent sets of generalized coordinates prior to their interpolation. The present paper also introduces a consistency enforcement approach for models defined on arbitrary underlying meshes. In the online phase, the operators are interpolated on matrix manifolds. The proposed framework is illustrated on two realistic multi-physics problems: an inverse acoustic scattering problem around a submarine and flutter predictions for a wing-tank system. The second application is implemented on a mobile device, illustrating the capability of the proposed framework to operate in real-time.
\end{abstract}

\end{frontmatter}


\section{Introduction}
Many engineering applications require the ability to generate predictions of the behavior of physical systems in real-time. Among those applications, one can mention design optimization, optimal control, the solution of inverse problems as well as uncertainty quantification. All of these applications require a large number of predictions for varying values of operating conditions. The operating conditions, usually described by a set of parameters, may define boundary conditions, initial conditions, physical or shape parameters that define the problem of interest and its underlying differential equations. However, each of these predictions usually demands computationally intensive calculations as accurate discretization of the underlying differential equations often leads to large scale systems of equations.

Projection-based model reduction~\cite{moore81,sirovich87} reduces the large computational cost associated with each solution of the underlying high-dimensional model (HDM) by reducing the number of degrees of freedom in the computation. For this purpose, a reduced-order basis (ROB) is defined and the solution is restricted to the subspace described by the ROB.
The current most challenging model reduction problems are those associated with nonlinear systems and parameter variations. Nonlinear systems require additional levels of approximation to enable large computational speedups~\cite{ryckelynck05,chaturantabut10,carlberg11,amsallem12:localROB,carlberg13,farhat14}. The model reduction of parameterized systems is also challenging due to the non-robustness of reduced-order model with respect to parameter variations that requires an appropriate offline training phase~\cite{veroy05,amsallem08,amsallem10,negri15}. Approaches addressing the model reduction of nonlinear parameterized systems are proposed in~\cite{pdt14,amsallem14:smo,amsallem14:morepas,wu15}. The focus of the present paper is on the efficient model reduction of linear parameterized systems.

More specifically, for linear parameterized systems, database approaches can be developed by pre-computing in an offline phase the reduced linear operators of a common reduced dimension for specific values of the parameters~\cite{amsallem08,amsallem09,amsallemthesis,amsallem10,degroote10,panzer10,amsallem11}.These linear operators are subsequently interpolated in the online phase for values of the parameters not present in the database. The small dimensionality of the reduced operators leads to real-time interpolation and predictions on-the-fly.

The interpolation of local reduced operators is however a challenging because each reduced operator is written in terms of a distinct set of generalized coordinates  corresponding to the local ROBs associated with each reduced-order model (ROM). To address this issue, approaches based on congruent transformations are proposed in~\cite{amsallemthesis,panzer10,amsallem11} when the underlying HDMs are defined on a common mesh. These approaches cannot, however, be applied when each  HDM is defined on a different mesh. The present paper introduces a novel approach, also based on congruent transformations, that addresses the challenge associated with arbitrary underlying meshes. 

Special care is also given in this paper to the interpolation step of the proposed procedure. The preservation of properties associated with the linear reduced operators can indeed be enforced by interpolating these operators on appropriate matrix manifolds~\cite{amsallem08,amsallem09,amsallemthesis,degroote10,amsallem11}. In that case, after appropriately mapping the reduced operators, interpolation can be carried out in the tangent space to the matrix manifold. As such, as long as the interpolation procedure preserves the tangent space, the interpolated quantity will also belong to the tangent space and can be mapped back to the manifold, leading to an interpolated reduced quantity that preserves its properties.

This paper is organized as follows. The problem of interest and the ROM database approach are formulated in Section~\ref{sec:pb}. The issue of consistency of reduced-order operators and its enforcement in the case of common and arbitrary underlying meshes are then investigated in Section~\ref{sec:consistency}. The interpolation of the resulting consistent database of ROMs is  developed in Section~\ref{sec:manifold}. Special care is given to the sampling, storage and exploitation of the database. The proposed approach is  applied in Section~\ref{sec:appli} to the model reduction of two parameterized systems. The first one is the analysis of a parameterized acoustic scattering system defined on arbitrary underlying meshes. In that case, the database approach is applied to the online solution of inverse problems. The second problem is the real-time flutter analysis of an aeroelastic system for flight conditions ranging from the subsonic to supersonic regimes. It is shown that the proposed approach successfully enables real-time predictions on a mobile device. Finally, conclusions are given in Section~\ref{sec:conclu}.

\section{Problem formulation and solution approach}\label{sec:pb}
In this paper, linear-time invariant parametric (LTIP) systems of one of the following two forms are considered
\begin{enumerate}
\item First-order LTIP systems of the form
\begin{equation}\label{eq:FOLTIP}
\Ebold(\mubold) \frac{d\wbold}{dt}(t) = \Abold(\mubold) \wbold(t) + \Bbold(\mubold)\ubold(t)\\
\end{equation}
and their formulation in the frequency domain
\begin{equation}
(j\omega \Ebold(\mubold)  - \Abold(\mubold)) \wbold(\omega) = \Bbold(\mubold)\ubold(\omega).
\end{equation}
The high-dimensional state vector is $\wbold\in\mathbb{R}^N$, $j^2=-1$, $t\geq 0$ denotes time and $\omega\geq 0$ frequency. $\Ebold$ and $\Abold$ are square high-dimensional matrices of dimension $N$. $\ubold\in\mathbb{R}^{N_i}$ denotes the input variable of dimension $N_i\ll N$ and $\Bbold\in\mathbb{R}^{N\times N_i}$. All operators depend on a vector of $N_{\mubold}$ parameters $\mubold\in\mathcal{D}\subset \mathbb{R}^{N_{\mubold}}$.

For both formulations, an output quantity of interest $\ybold\in\mathbb{R}^{N_o}$ is defined as
\begin{equation}\label{eq:output}
\ybold = \Gbold(\mubold) \wbold + \Hbold(\mubold)\ubold,
\end{equation}
with $N_o\ll N$ and $\Gbold\in\mathbb{R}^{N_o\times N}$ and $\Hbold\in\mathbb{R}^{N_o\times N_i}$.
\item Second-order LTIP systems of the form
\begin{equation}\label{eq:SOLTIP}
\Mbold(\mubold) \frac{d^2\wbold}{dt^2}(t) + \Cbold(\mubold) \frac{d\wbold}{dt}(t)  + \Kbold(\mubold) \wbold(t)= \Bbold(\mubold)\ubold(t)\\
\end{equation}
and their equivalent formulation in the frequency domain are considered
\begin{equation}\label{eq:SOLTIPF}
(-\omega^2 \Mbold(\mubold) + j\omega\Cbold(\mubold) + \Kbold(\mubold)  ) \wbold(\omega) = \Bbold(\mubold)\ubold(\omega).
\end{equation}
$\Mbold$, $\Cbold$ and $\Kbold$ are also parameter-dependent square linear operators of dimension $N$. For both formulations, an output quantity of interest $\ybold$ can be  defined as in~(\ref{eq:output}).
\end{enumerate}
The problem of interest is the fast computation of the output $\ybold$ for a given value of the parameters $\mubold\in\mathcal{D}$. More specifically, the computational cost associated with these computations should not scale with $N$ anymore.

To address this problem, an approach based on a database of linear projection-based ROMs is considered. This approach proceeds in two steps.
\begin{enumerate}
\item In the first offline step, $N_{\text{DB}}$ sample parameter values $\left\{\mubold_i\right\}_{i=1}^{N_{\text{DB}}} \in \mathcal{D}\subset\mathbb{R}^{N_{\mubold}}$ are selected and ROMs are constructed for each parameter value by defining right and left reduced-order bases (ROBs) $\Vbold(\mubold)\in\mathbb{R}^{N\times k}$ and $\Wbold(\mubold)\in\mathbb{R}^{N\times k}$, $k\ll N$ and approximating the state $\wbold$ as $\wbold\approx\Vbold(\mubold)\qbold$ where $\qbold\in\mathbb{R}^k$ is either solution of a reduced LTIP system in the time domain or frequency domain, as detailed below.  The output equation in terms of the reduced variable $\qbold$ is then
\begin{equation}
\ybold_r = \Gbold_r(\mubold)\qbold + \Hbold(\mubold)\ubold
\end{equation}
with $\Gbold_r(\mubold) = \Gbold(\mubold)\Vbold(\mubold)\in\mathbb{R}^{N_o\times k}$.
\begin{enumerate}
\item For first-order LTIP systems, the reduced equations in the time-domain are
\begin{equation}
\Ebold_r(\mubold) \frac{d\qbold}{dt}(t) = \Abold_r(\mubold) \qbold(t) + \Bbold_r(\mubold)\ubold(t)\\
\end{equation}
and in the frequency domain 
\begin{equation}
(j\omega \Ebold_r(\mubold)  - \Abold_r(\mubold)) \qbold(\omega) = \Bbold_r(\mubold)\ubold(\omega).
\end{equation}
where $\Ebold_r(\mubold) = \Wbold(\mubold)^T\Ebold(\mubold)\Vbold(\mubold)\in\mathbb{R}^{k\times k}$, $\Abold_r(\mubold) = \Wbold(\mubold)^T\Abold(\mubold)\Vbold(\mubold)\in\mathbb{R}^{k\times k}$ and $\Bbold_r(\mubold) = \Wbold(\mubold)^T\Bbold(\mubold)\in\mathbb{R}^{k\times N_i}$. 

\item For second-order LTIP systems, the time-domain reduced equations are
\begin{equation}
\Mbold_r(\mubold) \frac{d^2\qbold}{dt^2}(t) + \Cbold_r(\mubold) \frac{d\qbold}{dt}(t)  + \Kbold_r(\mubold) \qbold(t)= \Bbold_r(\mubold)\ubold(t)\\
\end{equation}
where  $\Mbold_r(\mubold) = \Wbold(\mubold)^T\Mbold(\mubold)\Vbold(\mubold)$, $\Cbold_r(\mubold) = \Wbold(\mubold)^T\Cbold(\mubold)\Vbold(\mubold)$ and $\Kbold_r(\mubold) = \Wbold(\mubold)^T\Kbold(\mubold)\Vbold(\mubold)$ are square reduced operators of dimension $k$. 
The reduced equations in the frequency domain are
\begin{equation}
(-\omega^2 \Mbold_r(\mubold) + j\omega\Cbold_r(\mubold) + \Kbold_r(\mubold)  ) \qbold(\omega) = \Bbold_r(\mubold)\ubold(\omega).
\end{equation}
\end{enumerate}
There are several model reduction techniques that can be applied to construct the ROBs $\Wbold(\mubold)$ and $\Vbold(\mubold)$ for LTI systems. Among those, the most popular are proper orthogonal decomposition (POD)~\cite{sirovich87,berkooz93}, balanced truncation~\cite{moore81} and moment matching~\cite{grimme97,hetmaniuk12}. 

Once the reduced operators are computed for the sampled values of the parameters, these reduced matrices are stored in a database of the form
\begin{equation}
\mathcal{DB} = \left\{ \mubold_i, \left(\Ebold_r(\mubold_i),\Abold_r(\mubold_i),\Bbold_r(\mubold_i),
\Gbold_r(\mubold_i),\Hbold(\mubold_i) \right) \right\}_{i=1}^{N_{\text{DB}}}
\end{equation}
for first-order systems and
\begin{equation}
\mathcal{DB} = \left\{ \mubold_i, \left(\Mbold_r(\mubold_i),\Cbold_r(\mubold_i),\Kbold_r(\mubold_i),\Bbold_r(\mubold_i),
\Gbold_r(\mubold_i),\Hbold(\mubold_i) \right) \right\}_{i=1}^{N_{\text{DB}}}
\end{equation}
for second-order systems.

In this work, the left and right ROBs are assumed to have orthonormal columns with respect to a common symmetric positive definite matrix $\boldsymbol{\mathcal{M}}$, that is $\Wbold(\mubold)^T\boldsymbol{\mathcal{M}}\Wbold(\mubold) = \Ibold_k$ and $\Vbold(\mubold)^T\boldsymbol{\mathcal{M}}\Vbold(\mubold) = \Ibold_k$. This property can be easily enforced~\emph{a posteriori} by applying a Gram-Schmidt orthogonalization procedure to the columns of non-orthonormal ROBs or directly in the ROB construction procedure~\cite{amsallem14:book}.

\item In the online phase, for a given value $\mubold^\star\in\mathcal{D}$ of the parameters, reduced operators are constructed  by interpolation of the elements of the database $\mathcal{DB}$.  

\end{enumerate}

Two technical issues however arise in the online interpolation step associated with the proposed approach:
\begin{enumerate}
\item The reduced quantities are not defined in the same system of reduced coordinates. Indeed, for each parameter $\mubold_i$, the system of reduced coordinates is defined by the local ROBs $\Vbold(\mubold_i)$ and $\Wbold(\mubold_i)$. A naive interpolation of the reduced operators may result in interpolation quantities that are not consistent with each other. A comprehensive approach to address this issue is proposed in Section~\ref{sec:consistency}
\item The linear operators stored in the database $\mathcal{DB}$ may have properties that should be preserved by interpolation. An approach relying on interpolation on a matrix manifold has been proposed in~\cite{amsallem09,amsallem11} to preserve these properties. The approach is briefly recalled and extended in Section~\ref{sec:manifold}.
\end{enumerate}
%

\section{Consistency between reduced-order models}\label{sec:consistency}
\subsection{Concept}
As underlined in the previous section and in~\cite{amsallem11}, the fact that local reduced operators are defined in different sets of generalized coordinates prevents their direct interpolation.  In this paper, two approaches are presented to address this issue. Both rely on a congruence transformation of the reduced operators to enforce consistency. They recognize the fact that the choice of local ROBs is not unique. Indeed, for a given right ROB $\Vbold(\mubold)$, any ROB of the form $\Vbold(\mubold)\Qbold$ with $\Qbold^T\Qbold=\Ibold_k$ defines an equally valid coordinate representation for the same ROM~\cite{amsallem11} with $\boldsymbol{\mathcal{M}}$-orthogonal columns. Similarly, for the left ROB $\Wbold(\mubold)$, any left ROB of the form $\Wbold(\mubold)\Zbold$ with $\Zbold^T\Zbold=\Ibold_k$ defines an equally valid basis. 

In turn, for a given first-order LTI ROM $\mathcal{R} = \left(\Ebold_r,\Abold_r,\Bbold_r,\Gbold_r,\Hbold\right)$, an equivalence class of ROMs under left and right multiplications by  orthogonal matrices $\Zbold$ and $\Qbold$ can be defined as
\begin{equation}
\mathcal{C}(\mathcal{R}) = \left\{  \left(\Zbold^T\Ebold_r\Qbold,\Zbold^T\Abold_r\Qbold,\Zbold^T\Bbold_r,\Gbold_r\Qbold,\Hbold\right)~\text{such that } \Zbold^T\Zbold=\Ibold,~\Qbold^T\Qbold=\Ibold\right\}.
\end{equation} 
Similarly, for a second-order LTI ROM $\mathcal{R}=\left(\Mbold_r,\Cbold_r,\Kbold_r,\Bbold_r,\Gbold_r,\Hbold\right)$, the equivalence class of ROMs under left and right multiplications by  orthogonal matrices $\Zbold$ and $\Qbold$  is
\begin{equation}
\mathcal{C}(\mathcal{R}) = \left\{  \left(\Zbold^T\Mbold_r\Qbold,\Zbold^T\Cbold_r\Qbold,\Zbold^T\Kbold_r\Qbold,\Zbold^T\Bbold_r,\Gbold_r\Qbold,\Hbold\right)~\text{such that } \Zbold^T\Zbold=\Ibold,~\Qbold^T\Qbold=\Ibold\right\}.
\end{equation}

Both approaches proposed in this paper will rely on a preprocessing step for which optimal transformations $\Qbold^\star(\mubold_i)$ and $\Zbold^\star(\mubold_i),~i=1,\cdots,N_{\text{DB}}$ are applied to the $N_{\text{DB}}$ ROMs stored in the database $\mathcal{DB}$.

The first approach, originally introduced in~\cite{amsallem11} and described in Section~\ref{sec:CommonMesh} is applicable whenever the underlying HDMs are defined on the same reference mesh. The location of the mesh nodes is potentially parameter dependent but the topology of the mesh is common to all parameter values. This requirement is relaxed in Section~\ref{sec:ArbitraryMesh} where a novel approach is introduced to enforce consistency in the case of arbitrary meshes. In that case, each HDM can have a different number of dofs.

A word of caution should however be formulated regarding interpolation of reduced linear operators. There are cases for which consistency cannot be enforced. Consider for instance the case of two configurations $\mubold_1$ and $\mubold_2$ associated with a common mesh and for which the right ROBs $\Vbold(\mubold_1)$ and $\Vbold(\mubold_2)$ are orthogonal to each other. In that case, the subspaces respectively defined by the ROBs are orthogonal as well and no transformation of the form $\{\Vbold(\mubold_i)\Qbold(\mubold_i)\}_{i=1}^2$ can define a consistent set of reduced coordinates. The degree of consistency between two ROBs will be quantitatively defined in the case of common underlying meshes in Section~\ref{sec:CommonMesh} and a truncation procedure introduced to further enforce consistency.

\subsection{Enforcement in the case of a common underlying mesh}\label{sec:CommonMesh} 

Consistency can be enforced in the case of a common underlying mesh by solving a series of Procrustes problems~\cite{amsallem11}. More specifically, given two local right reduced bases $\Vbold_i=\Vbold(\mubold_i)$ and   $\Vbold_j=\Vbold(\mubold_j)$, $\Vbold_j$ can be written in terms of $\Vbold_i$ as 
\begin{equation}
\Vbold_j = \Vbold_i \Rbold_{ij} + \Tbold_{ij}
\end{equation}
where $\Tbold_{ij}$ is the  component of $\Vbold_j$ that is $\boldsymbol{\mathcal{M}}$-orthogonal to $\Vbold_i$ i.e. $\Tbold_{ij}^T\boldsymbol{\mathcal{M}}\Vbold_i = \boldsymbol{0}$.

The subspace angles between the ROBs $\Vbold_i$ and $\Vbold_j$ define a measure of the maximum achievable consistency.   The subspace angles can be computed by the following three-step procedure:
\begin{enumerate}
\item Form $\Vbold_i^T\boldsymbol{\mathcal{M}}\Vbold_j = \Rbold_{ij}$.
\item Compute a singular value decomposition $\Rbold_{ij} = \Xbold \Sigmabold \Ybold^T$ where $\Xbold = [\xbold_1,\dots,\xbold_k]$, $\Ybold = [\ybold_1,\dots,\ybold_k]$ and $\Sigmabold= \text{diag}(\sigma_1,\cdots,\sigma_k)$.
\item Compute the subspace angles as $\theta_\ell = \arccos(\sigma_\ell),~\ell=1,\cdots,k$. The  canonical vectors associated with each angle $\theta_\ell$ are $(\Vbold_i\xbold_\ell,\Vbold_j\ybold_\ell)$. Note that the angles are ordered increasingly as $0 \leq \theta_1\leq \cdots \leq \theta_k \leq \frac{\pi}{2}$.
\end{enumerate}

A subspace angle $\theta_\ell$ that is equal to zero reflects perfect consistency between the associated vectors $\Vbold_i\xbold_\ell$ and $\Vbold_j\ybold_\ell$.
In general, angles that are greater than a threshold $\theta_{\max} = \frac{\pi}{4}$ may define cases for which consistency cannot be achieved. One option to address this issue that is discussed below is to truncate the directions associated with those large angles. An alternate option is to refine the database until smaller subspace angles are achieved.

For a given database $\mathcal{DB}$, optimal transformations $\{\Qbold(\mubold_i)\}_{i=1}^{N_{\text{DB}}}$ can be computed by fixing one of the ROBs (say $\Qbold(\mubold_{i_0}) = \Ibold_k$) and computing the transformations as the minimizers of the following series of Procrustes problems:
 \begin{gather}
 \begin{split}
 \Qbold(\mubold_i) &= \argmin_{\Sbold\in\mathbb{R}^{k\times k}} \| \Vbold_i\Sbold - \Vbold_{i_0}\|_{\boldsymbol{\mathcal{M}}},~~i=1,\cdots,N_{\text{DB}},\\
& ~~~~~\text{s.t.}~ \Sbold^T\Sbold = \Ibold_k,
 \end{split}
 \end{gather}
where $\|\Nbold\|_{\boldsymbol{\mathcal{M}}} = \|\boldsymbol{\mathcal{M}}^{\frac{1}{2}}\Nbold\|_F$. The optimal transformation $\Qbold(\mubold_i)$ can be determined analytically from the SVD of $\Rbold_{ii_0}$ as
\begin{equation}
\Qbold(\mubold_i) = \Xbold \Ybold^T.
\end{equation}
As stated above, truncation of the ROBs can be used to enforce consistency. For each ROB $\Vbold_i$, $i=1,\cdots,N_{\text{DB}}$, the maximum index $\ell_i$ for which the subspace angles with $\Vbold_{i_0}$ are smaller than $\theta_{\max}$  can be determined and all ROM truncated to the index 
\begin{equation}
L = \min_{i=1,\cdots,N_{\text{DB}}} \ell_i.
\end{equation}
Truncating the ROMs will result in a loss of accuracy of each ROM when compared to the underlying HDM. However, this truncation step will  improve accuracy after interpolation of the elements of the database as these will be more consistent.

\subsection{Enforcement in the case of arbitrary underlying meshes}\label{sec:ArbitraryMesh}

\begin{figure}[htbp] \centering
\includegraphics[width=0.5\textwidth]{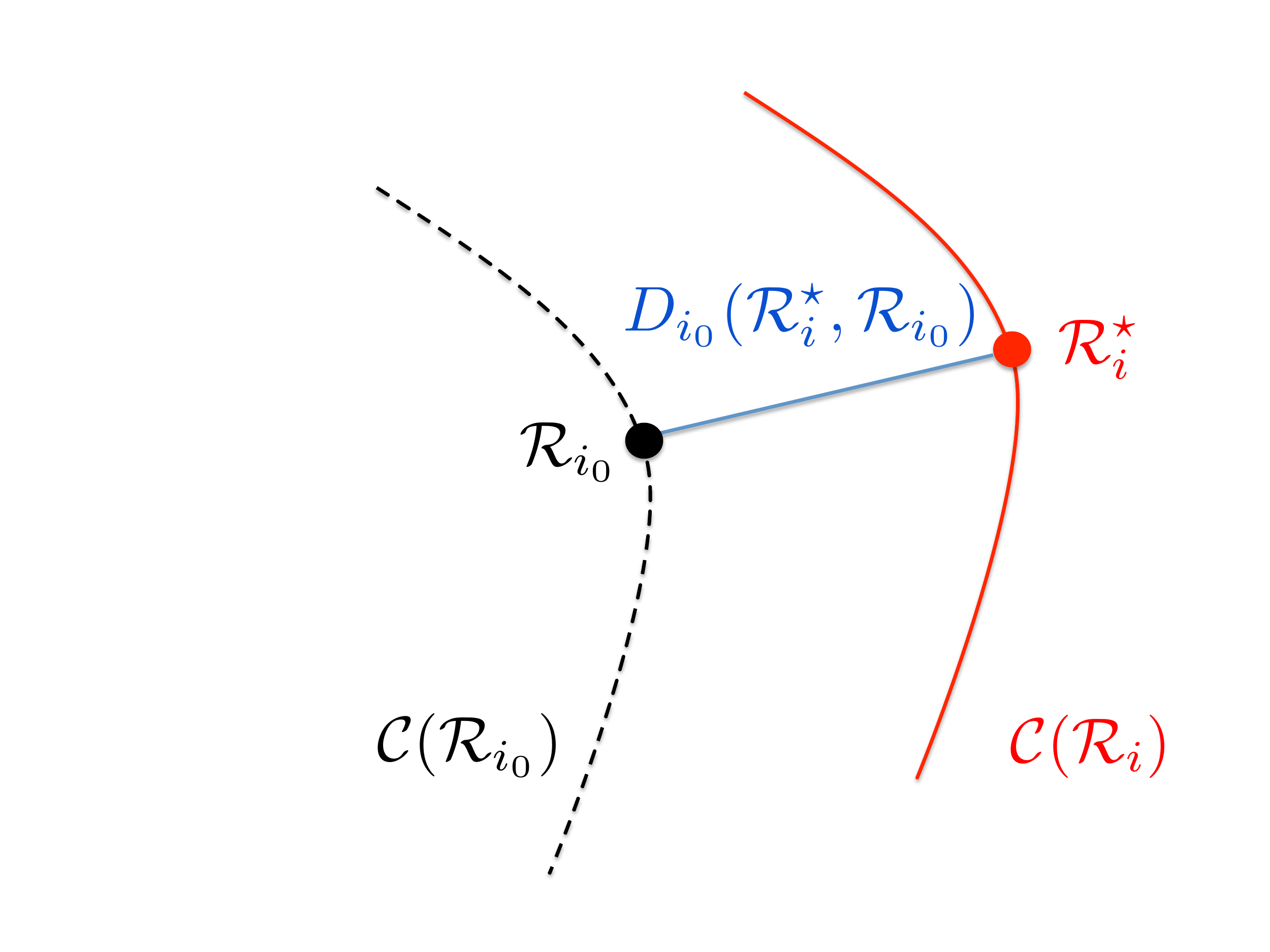}   
\caption{Schematic description of the minimization problem for consistency enforcement in the case of arbitrary meshes}
\label{fig:minPb} \end{figure}

In the case of arbitrary meshes, subspace angles cannot be defined as the underlying HDM spaces may be of different, parameter-dependent dimensions $N(\mubold_i)$. To address this issue, a heuristic procedure is developed in this section to enforce consistency for that specific scenario. As in the case of the Procrustes problem, one of the ROMs $\mathcal{R}_{i_0}$ is selected to define a reference configuration. Then, for each ROM $\mathcal{R}_i,~i=1,\cdots,N_{\text{DB}}$ in the database, a transformed ROM $\mathcal{R}_i^\star$ is determined as the minimizer of a measure of distance of the reference ROM $\mathcal{R}_{i_0}$ to the equivalence class $\mathcal{C}(\mathcal{R}_{i})$ as
\begin{equation}\label{eq:optPb}
\mathcal{R}_i^\star = \argmin_{\mathcal{R}\in\mathcal{C}(\mathcal{R}_i)} D_{i_0}(\mathcal{R},\mathcal{R}_{i_0}).
\end{equation}
The minimization problem is schematically depicted in Fig.~\ref{fig:minPb}.

The measure of distance $D_{i_0}(\mathcal{R},\mathcal{R}_{i_0})$ is defined as follows:
\begin{enumerate}
\item For a first-order system, the distance between $\mathcal{R} = (\Ebold_r,\Abold_r,\Bbold_r,\Gbold_r,\Hbold)$ and $\mathcal{R}' = (\Ebold'_r,\Abold'_r,\Bbold'_r,\Gbold'_r,\Hbold')$ is defined as the normalized expression
\begin{equation}
D_{i_0}(\mathcal{R},\mathcal{R}') = \epsilon \| \Ebold_r - \Ebold'_r\|_F^2 + \alpha  \| \Abold_r - \Abold'_r\|_F^2 +  \beta \| \Bbold_r - \Bbold'_r\|_F^2 + \gamma  \| \Gbold_r - \Gbold'_r\|_F^2 + \eta \| \Hbold - \Hbold'\|_F^2
\end{equation}
where
\begin{equation}
\epsilon = \frac{1}{\|\Ebold^0_r\|_F^2},~\alpha = \frac{1}{\|\Abold^0_r\|_F^2},~\beta = \frac{1}{\|\Bbold^0_r\|_F^2},~\gamma = \frac{1}{\|\Gbold^0_r\|_F^2},~\eta = \frac{1}{\|\Hbold^0\|_F^2}
\end{equation}
are normalization constants based on the reduced operators in $\mathcal{R}_{i_0} = (\Ebold^0_r,\Abold^0_r,\Bbold^0_r,\Gbold^0_r,\Hbold^0)$.
\item For a second-order system, the distance between $\mathcal{R} = (\Mbold_r,\Cbold_r,\Kbold_r,\Bbold_r,\Gbold_r,\Hbold)$ and $\mathcal{R}' = (\Mbold'_r,\Cbold'_r,\Kbold'_r,\Bbold'_r,\Gbold'_r,\Hbold')$ is defined as
\begin{equation}
D_{i_0}(\mathcal{R},\mathcal{R}') = \mu \| \Mbold_r - \Mbold'_r\|_F^2 + \xi  \| \Cbold_r - \Cbold'_r\|_F^2 + \kappa  \| \Kbold_r - \Kbold'_r\|_F^2 +  \beta \| \Bbold_r - \Bbold'_r\|_F^2 + \gamma  \| \Gbold_r - \Gbold'_r\|_F^2 + \eta \| \Hbold - \Hbold'\|_F^2
\end{equation}
where
\begin{equation}
\mu = \frac{1}{\|\Mbold^0_r\|_F^2},~\xi = \frac{1}{\|\Cbold^0_r\|_F^2},~\kappa = \frac{1}{\|\Kbold^0_r\|_F^2}
\end{equation}
are normalization constants based on the reduced operators in $\mathcal{R}_{i_0} = (\Mbold^0_r,\Cbold^0_r,\Kbold^0_r,\Bbold^0_r,\Gbold^0_r,\Hbold^0)$.
\end{enumerate}

In practice, since the class $\mathcal{C}(\mathcal{R}_{i})$ is parameterized by  two transformation matrices $\Qbold$ and $\Zbold$, the optimization problem~(\ref{eq:optPb}) can be explicitly written for a first-order system in terms of $\mathcal{R}_i = (\Ebold_{ri},\Abold_{ri},\Bbold_{ri},\Gbold_{ri},\Hbold_i)$, $\mathcal{R}_{i_0} $ and $(\Qbold,\Zbold)$  as 
\begin{gather}\label{eq:QZpb1}
\begin{split}
&\min_{\Qbold,\Zbold}~~ \epsilon \| \Zbold^T\Ebold_{ri}\Qbold - \Ebold^0_r\|_F^2 + \alpha  \| \Zbold^T\Abold_r\Qbold - \Abold^0_r\|_F^2 +  \beta \| \Zbold^T\Bbold_r - \Bbold^0_r\|_F^2 + \gamma  \| \Gbold_r\Qbold - \Gbold^0_r\|_F^2 + \eta \| \Hbold - \Hbold^0\|_F^2\\
&~~\text{s.t.}~\Qbold^T\Qbold=\Ibold_k,~~\Zbold^T\Zbold=\Ibold_k.
\end{split}
\end{gather}
A similar expression holds for second-order systems. The rest of this section will focus on first-order systems but the analysis directly carries over to second-order systems as well. 

\emph{Remark.} In the case of Galerkin projection, $\Wbold(\mubold_i)=\Vbold(\mubold_i)$ and $\Qbold=\Zbold$. Then~(\ref{eq:optPb}) simplifies to
\begin{gather}\label{eq:QZpb1_1}
\begin{split}
&\min_{\Qbold} ~~\epsilon \| \Qbold^T\Ebold_{ri}\Qbold - \Ebold^0_r\|_F^2 + \alpha  \| \Qbold^T\Abold_r\Qbold - \Abold^0_r\|_F^2 +  \beta \| \Qbold^T\Bbold_r - \Bbold^0_r\|_F^2 + \gamma  \| \Gbold_r\Qbold - \Gbold^0_r\|_F^2 + \eta \| \Hbold - \Hbold^0\|_F^2\\
&~~\text{s.t.}~\Qbold^T\Qbold=\Ibold_k.
\end{split}
\end{gather}

Problem~(\ref{eq:QZpb1}) is equivalent to the maximization problem~\cite{amsallemthesis}
\begin{gather}\label{eq:QZpb2}
\begin{split}
&\max_{\Qbold,\Zbold}~~  \langle\epsilon \Zbold^T\Ebold_{ri}\Qbold ,\Ebold^0_r\rangle  +  \langle\alpha  \Zbold^T\Abold_{ri}\Qbold ,\Abold^0_r\rangle +  \langle\beta  \Bbold_r\left(\Bbold^0_r\right)^T ,\Zbold\rangle  + \langle \gamma  (\Gbold^0_r)^T\Gbold_r,\Qbold \rangle \\
&~~\text{s.t.}~\Qbold^T\Qbold=\Ibold_k,~~\Zbold^T\Zbold=\Ibold_k,
\end{split}
\end{gather}
where 
\begin{equation}
\langle \Mbold,\Nbold \rangle = \text{tr}(\Mbold^T\Nbold), ~~\Mbold,\Nbold \in\mathbb{R}^{m\times n}.
\end{equation}

The problem of maximizing the first term in Eq.~(\ref{eq:QZpb2}) has been studied in the literature~\cite{fraikin08,helmkebook} in the case $\Qbold=\Zbold$ (Galerkin projection). This first term defines a correlation criterion between the matrices $\Ebold_{ri}$ and $\Ebold_r^0$. A solution to that problem developed in Ref.~\cite{fraikin08}, consists of defining an iterative algorithm whose fixed points are the critical points of the maximization problem. This approach is here extended to the optimization problem of interest for enforcing consistency between ROM operators. Both Galerkin and Petrov-Galerkin projections are considered as follows.
\begin{itemize}
\item Galerkin projection.

In this case, the functional in~(\ref{eq:QZpb2}) is
\begin{equation}
\mathcal{J}_G(\Qbold) =  \langle\epsilon \Qbold^T\Ebold_{ri}\Qbold ,\Ebold^0_r\rangle  +  \langle\alpha  \Qbold^T\Abold_{ri}\Qbold ,\Abold^0_r\rangle +  \langle\beta  \Bbold_r\left(\Bbold^0_r\right)^T ,\Qbold\rangle  + \langle \gamma  (\Gbold^0_r)^T\Gbold_r,\Qbold \rangle.
\end{equation}
Adapting the algorithm developed in~\cite{fraikin08} to the present case, the iterative algorithm is  based on an affine map defined as
\begin{equation}
\Mbold_{s,\text{G}}(\Qbold) =  \epsilon\left(\Ebold_{ri}\Qbold \left(\Ebold_{r}^0\right)^T + \Ebold_{ri}^T \Qbold \Ebold_{r}^0 \right)+ \alpha\left(\Abold_{ri}\Qbold \left(\Abold_{r}^0\right)^T + \Abold_{ri}^T \Qbold \Abold_{r}^0 \right)+ s\Qbold + \beta \Bbold_{ri}\left( \Bbold_{r}^0\right)^T + \gamma\Gbold_{ri}^T \Gbold_{r}^0,
\end{equation}
where $s$ is a fixed real parameter chosen such that $s>s_{\text{min,G}}$ with
\begin{equation}\label{eq:smindef}
s_{\text{min,G}} = 2\epsilon \left\|\Ebold_{ri}\right\|_2 \left\|\Ebold_{r}^0\right\|_2 + 2\alpha \left\|\Abold_{ri}\right\|_2 \left\|\Abold_{r}^0\right\|_2 + \left\|\beta \Bbold_{ri} \left(\Bbold_{r}^0\right)^T + \gamma\Gbold_{ri}^T \Gbold_{r}^0\right\|_2.
\end{equation}
Defining the parameter $s$ is necessary to ensure that the fixed points of the proposed iterative algorithm are exactly the critical points of the maximization problem (see Theorem 1). 

The proposed procedure then proceeds by iteratively solving the maximization problem
\begin{equation}
\Qbold_{j+1} = \text{arg} \max_{\Sbold^T\Sbold=\Ibold_k}\langle \Sbold, \Mbold_{s,\text{G}}(\Sbold_j)\rangle,~~j=0,\cdots.
\end{equation}
The solution to this problem is established in the following lemma, leading to the proposed iterative procedure presented in Algorithm~\ref{alg:iterAlgG}.

{\bf{Lemma.}}  Let the singular value decomposition of $\Mbold_{s,\text{G}}(\Qbold_j)$ be defined as
\begin{equation}
\Mbold_{s,\text{G}}(\Qbold_j) = \Ubold \Sigmabold \Vbold^T.
\end{equation}
Then
\begin{equation}
\max_{\Qbold^T\Qbold = \Ibold_k} \langle \Qbold, \Mbold_s(\Qbold_j)\rangle = \sum_{\ell=1}^k \sigma_\ell(\Mbold_{s,\text{G}}(\Qbold_j)),
\end{equation}
where $\{ \sigma_\ell(\Mbold_{s,\text{G}}(\Qbold_j))\}_{\ell=1}^k$ is the set of singular values of $\Mbold_{s,\text{G}}(\Qbold_j)$. The solution to the maximization problem is unique and equal to $\Ubold \Vbold^T$.

A proof is offered in~\cite{fraikin08} in a more general setting.
\begin{algorithm}[htbp]
\caption{Fixed-point procedure in the case of Galerkin projection}
\begin{algorithmic}[1] \label{alg:iterAlgG}
\STATE Compute $s>2\epsilon \left\|\Ebold_{ri}\right\|_2 \left\|\Ebold_{r}^0\right\|_2 + 2\alpha \left\|\Abold_{ri}\right\|_2 \left\|\Abold_{r}^0\right\|_2 + \left\|\beta \Bbold_{ri} \left(\Bbold_{r}^0\right)^T + \gamma\Gbold_{ri}^T \Gbold_{r}^0\right\|_2$.
\STATE Compute $\Fbold = \beta\Bbold_{ri} \left(\Bbold_{r}^0\right)^T +\gamma \Gbold_{ri}^T \Gbold_{r}^0$.
\STATE Choose an orthogonal initial matrix $\Qbold_0\in\mathbb{R}^{k\times k}$.
\FOR{$j=0,\cdots$}
	\STATE Compute the map
	\[ \Mbold_{s,\text{G}}(\Qbold_j) = \epsilon\left( \Ebold_{ri} \Qbold_j \left(\Ebold_{r}^0\right)^T + \Ebold_{ri}^T \Qbold_j \Ebold_{r}^0 \right) + \alpha\left( \Abold_{ri} \Qbold_j \left(\Abold_{r}^0\right)^T + \Abold_{ri}^T \Qbold_j \Abold_{r}^0 \right)+ s\Qbold_j + \Fbold\]
         \STATE Compute its SVD
	\[\Ubold_{j+1}\Sigmabold_{j+1}\Vbold_{j+1}^T =  \Mbold_{s,\text{G}}(\Qbold_j)\]
	\STATE $\Qbold_{j+1} = \Ubold_{j+1} \Vbold^T_{j+1}$
\ENDFOR
\end{algorithmic}
\end{algorithm}

In order to show that the fixed points of the recursive algorithm are the critical points of $\mathcal{J}_G$, one needs to characterize these critical points. This is done in the following theorem. 

{\bf{Theorem 1.}} The critical points $\Qbold^\star$ of $\mathcal{J}_G$ are orthogonal matrices satisfying the identity
\begin{equation}\label{eq:criticalPts}
\Qbold^\star \Sbold =  \epsilon\left(\Ebold_{ri}\Qbold^\star\left(\Ebold_{r}^0\right)^T + \Ebold_{ri}^T\Qbold^\star\Ebold_{r}^0\right) + \alpha\left(\Abold_{ri}\Qbold^\star\left(\Abold_{r}^0\right)^T + \Abold_{ri}^T\Qbold^\star\Abold_{r}^0\right) +\beta \Bbold_{ri} \left(\Bbold_{r}^0\right)^T + \gamma\Gbold_{ri}^T \Gbold_{r}^0,
\end{equation}
where $\Sbold$ is a symmetric matrix.

A proof is offered in Appendix 1. 

{\bf{Theorem 2.}} The set of the fixed points of the recursive Algorithm~\ref{alg:iterAlgG} is exactly the set of the critical points of $\mathcal{J}_G$.

 A proof of the theorem is presented in Appendix 2.

\item Petrov-Galerkin projection.

Defining the functional
\begin{equation}
\mathcal{J}_{PG}(\Qbold,\Zbold) =  \langle\epsilon \Zbold^T\Ebold_{ri}\Qbold ,\Ebold^0_r\rangle  +  \langle\alpha  \Zbold^T\Abold_{ri}\Qbold ,\Abold^0_r\rangle +  \langle\beta  \Bbold_r\left(\Bbold^0_r\right)^T ,\Zbold\rangle  + \langle \gamma  (\Gbold^0_r)^T\Gbold_r,\Qbold \rangle,
\end{equation}
the iterative algorithm is now  based on the block-affine map defined as
\begin{gather}
\begin{split}
\Mbold_{s,\text{PG}}(\Qbold,\Zbold) &=  \left[\begin{array}{cc}\Mbold^Q_{s,\text{PG}}(\Qbold,\Zbold)  & \boldsymbol{0} \\\boldsymbol{0} &\Mbold^Z_{s,\text{PG}}(\Qbold,\Zbold) \end{array}\right]  \\
&=\epsilon\left[\begin{array}{cc}\boldsymbol{0} & \Ebold_{ri} \\\Ebold_{ri}^T & \boldsymbol{0}\end{array}\right]  \left[\begin{array}{cc}\Qbold  & \boldsymbol{0} \\\boldsymbol{0} & \Zbold\end{array}\right]    \left[\begin{array}{cc}\boldsymbol{0} & \left(\Ebold_{r}^0\right)^T \\\left(\Ebold_{r}^0\right)^T & \boldsymbol{0}\end{array}\right]  \\
& ~~~~+ \alpha\left[\begin{array}{cc}\boldsymbol{0} & \Abold_{ri} \\\Abold_{ri}^T & \boldsymbol{0}\end{array}\right]  \left[\begin{array}{cc}\Qbold  & \boldsymbol{0} \\\boldsymbol{0} & \Zbold\end{array}\right]    
  \left[\begin{array}{cc}\boldsymbol{0} & \left(\Abold_{r}^0\right)^T \\\left(\Abold_{r}^0\right)^T & \boldsymbol{0}\end{array}\right] \\
  &~~~~+ s \left[\begin{array}{cc}\Qbold  & \boldsymbol{0} \\\boldsymbol{0} & \Zbold\end{array}\right]   + \left[\begin{array}{cc} \beta \Bbold_{ri}\left( \Bbold_{r}^0\right)^T  & \boldsymbol{0} \\\boldsymbol{0} & \gamma\Gbold_{ri}^T \Gbold_{r}^0\end{array}\right],
\end{split}
\end{gather}
where $s$ is  chosen such that $s>s_{\text{min,PG}}$ with
\begin{equation}\label{eq:smindef2}
s_{\text{min,PG}} = \epsilon \left\|\Ebold_{ri}\right\|_2 \left\|\Ebold_{r}^0\right\|_2 + \alpha \left\|\Abold_{ri}\right\|_2 \left\|\Abold_{r}^0\right\|_2 + \max\left(\beta\left\| \Bbold_{ri} \left(\Bbold_{r}^0\right)^T\|_2, \gamma\|\Gbold_{ri}^T \Gbold_{r}^0\right\|_2\right).
\end{equation}
Similarly as in the case of Galerkin projection, a fixed point procedure is defined in Algorithm~\ref{alg:iterAlgPG}.
\begin{algorithm}[htbp]
\caption{Fixed-point procedure in the case of Petrov-Galerkin projection}
\begin{algorithmic}[1] \label{alg:iterAlgPG}
\STATE Compute $s> \epsilon \left\|\Ebold_{ri}\right\|_2 \left\|\Ebold_{r}^0\right\|_2 + \alpha \left\|\Abold_{ri}\right\|_2 \left\|\Abold_{r}^0\right\|_2 + \max\left(\beta\left\| \Bbold_{ri} \left(\Bbold_{r}^0\right)^T\|_2,\gamma\|\Gbold_{ri}^T \Gbold_{r}^0\right\|_2\right)$.
\STATE Compute $\Fbold_B = \beta\Bbold_{ri} \left(\Bbold_{r}^0\right)^T$ and $\Fbold_G =\gamma \Gbold_{ri}^T \Gbold_{r}^0$.
\STATE Choose orthogonal initial matrices $\Qbold_0\in\mathbb{R}^{k\times k}$ and $\Zbold_0\in\mathbb{R}^{k\times k}$.
\FOR{$j=0,\cdots$}
	\STATE Compute the maps
\[\Mbold^Q_{s,\text{PG}}(\Qbold_j,\Zbold_j) = \epsilon\Ebold_{ri}\Zbold_j\left(\Ebold_r^{0}\right)^T  +\alpha\Abold_{ri}\Zbold_j\left(\Abold_r^{0}\right)^T + s\Qbold_j + \Fbold_G \]
and
\[\Mbold^Z_{s,\text{PG}}(\Qbold_j,\Zbold_j) = \epsilon\Ebold_{ri}^T\Qbold_j\left(\Ebold_r^{0}\right)^T  +\alpha\Abold_{ri}^T\Qbold_j\left(\Abold_r^{0}\right)^T +s\Zbold_j +\Fbold_B \]

	         \STATE Compute their SVDs
	\[\Ubold^Q_{j+1}\Sigmabold^Q_{j+1}\Vbold_{j+1}^{QT} =  \Mbold^Q_{s,\text{PG}}(\Qbold_j)\]
	and
	\[\Ubold^Z_{j+1}\Sigmabold^Z_{j+1}\Vbold_{j+1}^{ZT} =  \Mbold^Z_{s,\text{PG}}(\Qbold_j)\]
	\STATE $\Qbold_{j+1} = \Ubold^Q_{j+1} \Vbold^{QT}_{j+1}$
	\STATE $\Zbold_{j+1} = \Ubold^Z_{j+1} \Vbold^{ZT}_{j+1}$
\ENDFOR
\end{algorithmic}
\end{algorithm}
The following theorems, whose proofs follow closely the ones of Theorems 1 and 2, establish the fact that the fixed point procedure in Algorithm~\ref{alg:iterAlgPG}
 can be used to find critical points of $\mathcal{J}_{PG}$.
 
 {\bf{Theorem 3}} The critical points $(\Qbold^\star,\Zbold^\star)$ of $\mathcal{J}_{PG}$ are orthogonal matrices satisfying the identities
\begin{gather}\label{eq:criticalPtsPG}
\begin{split}
\Qbold^\star \Sbold_Q &=  \epsilon\Ebold_{ri}\Zbold^\star\left(\Ebold_{r}^0\right)^T + \alpha\Abold_{ri}\Zbold^\star\left(\Abold_{r}^0\right)^T + \gamma\Gbold_{ri}^T \Gbold_{r}^0,\\
\Zbold^\star \Sbold_Z &=  \epsilon\Ebold_{ri}^T\Qbold^\star\left(\Ebold_{r}^0\right)^T + \alpha\Abold_{ri}^T\Qbold^\star\left(\Abold_{r}^0\right)^T+\beta \Bbold_{ri} \left(\Bbold_{r}^0\right)^T, 
\end{split}
\end{gather}
where $\Sbold_Q$ and  $\Sbold_Z$ are symmetric matrices.

{\bf{Theorem 4.}} The set of the fixed points of the recursive Algorithm~\ref{alg:iterAlgPG} is exactly the set of the critical points of $\mathcal{J}_{PG}$.

\end{itemize}

\emph{Remark.} For the case of arbitrary meshes, assessing consistency is a more difficult task as two conflicting factors intervene in the distance measure $D_{i_0}$: (1) the inconsistency arising from a choice of two distinct sets of coordinates and (2) the inherent variation of the ROM operators due to parameter changes. 

\subsection{Consistent set of reduced-order models}
After the computation of the optimal transformation operators $\left\{\left(\Zbold^\star(\mubold_i),\Qbold^\star(\mubold_i)\right)\right\}_{i=1}^{N_{\mubold}}$, the reduced operators in the database $\mathcal{DB}$ are transformed accordingly as
\begin{gather}
\begin{split}
\mathcal{DB} &= \left\{ \mubold_i;\left( \Ebold^\star_r(\mubold_i),\Abold^\star_r(\mubold_i),\Bbold^\star_r(\mubold_i),\Gbold^\star_r(\mubold_i),\Hbold^\star_r(\mubold_i)    \right) \right\}_{i=1}^{N_{\mubold}}\\
&= \left\{ \mubold_i;\left( \Zbold^\star(\mubold_i)^T\Ebold_r(\mubold_i)\Qbold^\star(\mubold_i),\Zbold^\star(\mubold_i)^T\Abold_r(\mubold_i)\Qbold^\star(\mubold_i),\Zbold^\star(\mubold_i)^T\Bbold_r(\mubold_i),\Gbold^\star_r(\mubold_i)\Qbold^\star(\mubold_i),\Hbold_r(\mubold_i)    \right) \right\}_{i=1}^{N_{\mubold}}.
\end{split}
\end{gather}
Similar expressions hold for the case of second-order systems.

\section{Interpolation in a database of ROMs on matrix manifolds}\label{sec:manifold}

\subsection{Interpolation}\label{sec:interp}
As indicated in Section~\ref{sec:pb}, the interpolation of the linear operators stored in the database $\mathcal{DB}$ should often preserve properties of the operators such as symmetry, positivity, orthogonality or non-singularity. An approach to preserve these properties was first presented in~\cite{amsallem08}. It is based on the interpolation on the tangent space of the appropriate manifold and was applied to the case of interpolation of reduced operators in~\cite{amsallem09,degroote10,amsallem11}.

The algorithm proceeds for each of the elements of the database in four steps as follows:
\begin{enumerate}
\item An identification of the manifold the reduced matrices belong to
\item A mapping (logarithmic map) of all the database reduced matrices to the tangent space of the manifold at one of the database points
\item An interpolation of the mapped quantities in the tangent space at the target parameter $\mubold$
\item A mapping (exponential map) of the interpolated quantity back to the manifold leading to a reduced operator at the target parameter $\mubold$
\end{enumerate}
More details on the interpolation algorithm as well as the formulas for computing the mapping are provided in~\cite{amsallem11,amsallemthesis}.

In practice there may be several choices for an interpolation procedure on matrix manifolds as underlined by the following two cases.
\begin{itemize}
\item In~\cite{degroote10}, the authors develop a heuristic technique for interpolating non-singular matrices either on the manifold on non-singular matrices or square matrices. The heuristic proceeds by selecting the manifold for which a nonlinearity criterion is the smallest. This heuristic is applied in Section~\ref{sec:flutter}.
\item An alternative to interpolating symmetric positive definite matrices on the tangent space to that manifold is to use the Choleski factorization. This novel approach is described in Algorithm~\ref{alg:Choleski} and applied in Section~\ref{sec:scattering}. This approach avoids selecting one of the database points and interpolating in its associated tangent space. It preserves the SPD properties of the matrices as long as the interpolated quantity on the diagonal of the Choleski factor are all strictly non zero.
\end{itemize}

\begin{algorithm}[htbp]
\caption{Interpolation of SPD matrices by Choleski factorization}
\begin{algorithmic}[1] \label{alg:Choleski}
\FOR{$i=1,\cdots,N_{\text{DB}}$}
	\STATE Compute the Choleski factorization
	\[ \Kbold_{ri}= \Sbold_i\Sbold_i^T\]
\ENDFOR
        \STATE Interpolate the Choleski factors $\{\Sbold_i\}_{i=1}^{N_{\text{DB}}}$, leading to an interpolated factor $\Sbold^\star$
        \STATE Compute the interpolated matrix as
        \[ \Kbold_r^\star = \Sbold^\star\Sbold^{\star T}\]
\end{algorithmic}
\end{algorithm}

There is no restriction on the interpolation technique in the tangent space to the matrix manifold of interest as long as the it leads to an interpolated quantity that preserves the tangent space~\cite{amsallem11}. In~\cite{zimmermann14}, the author identifies an interpolation technique that does not preserve that property.
When the database parameters belong to a lattice of points, spline or polynomial interpolation can be used in the tangent space. When the dimension $N_{\mubold}$ of the parameter domain is large, however, interpolating from a lattice of points is not an option and instead, interpolation based on radial basis functions or Kriging can be used instead~\cite{amsallem08,amsallemthesis,choi15}.

\subsection{Sampling}
The selection of sample points $\{\mubold_i\}_{i=1}^{N_{DB}}$ is an important step that influences the accuracy of the resulting interpolation approach. A poor choice of sample points will result in large errors of the proposed procedure in some regions of the parameter domain $\mathcal{D}$. A priori sampling techniques such as factorial and latin hypercube sampling can be used to provide a uniform coverage of $\mathcal{D}$.  Alternatively, greedy techniques that iteratively sample the regions of the parameter space associated with the largest ROM error can provide a selection of the samples that is more suited for the problem of interest. Such greedy techniques have been introduced in the context of model reduction in general in~\cite{veroy05,grepl05,buithanh08,pdt14,amsallem14:smo,amsallem14:expintegrators} and for interpolation of LTIP ROM systems in particular in~\cite{choi15}. A priori sampling will be used in the application of Section~\ref{sec:flutter} and a greedy sampling approach developed for the inverse problem application of Section~\ref{sec:scattering}.

\subsection{Storage and exploitation}
In practice the reduced operators are stored after their congruence transformation in one database $\mathcal{DB}$ or several sub-databases $\{\mathcal{DB_s}\}_{s=1}^{N_s}$ of the form
\begin{equation}
\mathcal{DB} = \bigcup_{s=1}^{N_{S}}\mathcal{DB}_s = \bigcup_{s=1}^{N_{S}} \left\{\mubold_i;\left(\Ebold^\star_r(\mubold_i),\Abold^\star_r(\mubold_i),\Bbold_r^\star(\mubold_i),\Gbold^\star_r(\mubold_i),\Hbold^\star(\mubold_i)\right)\right\}_{i=1}^{N_{DB,s}}
\end{equation}
for first-order systems and
\begin{equation}
\mathcal{DB} = \bigcup_{s=1}^{N_{S}}\mathcal{DB}_s = \bigcup_{s=1}^{N_{S}} \left\{\mubold_i;\left(\Mbold_r^\star(\mubold_i),\Cbold_r^\star(\mubold_i),\Kbold_r^\star(\mubold_i),\Bbold_r^\star(\mubold_i),\Gbold_r^\star(\mubold_i),\Hbold^\star(\mubold_i)\right)\right\}_{i=1}^{N_{DB,s}}
\end{equation}
for second-order systems.

Storing the database is inexpensive as it only involves reduced operators. In practice, a database with $N_{DB} = \sum_{s=1}^S N_{DB,s}$ contains $N_{DB}(N_{\mubold}+2k^2+k(N_i+N_o) + N_iN_o)$ matrix entries for first-order systems of the form~(\ref{eq:FOLTIP}) and $N_{DB}(N_{\mubold}+3k^2+k(N_i+N_o) + N_iN_o)$ entries for second-order systems of the form~(\ref{eq:SOLTIP}).

The parameter domain $\mathcal{D}$ is in practice subdivided in $N_S$ non-overlaping subdomains  $\{\mathcal{D}_s\}_{s=1}^{N_S}$ such that
\begin{equation}
\mathcal{D} = \bigcup_{s=1}^{N_S} \mathcal{D}_s
\end{equation}
and each subdomain $\mathcal{D}_s$ is associated with the sub-database $\mathcal{DB}_s$.
Then, in the online phase, for a new value $\widehat\mubold\in\mathcal{D}$ of the parameters, the sub-database $\mathcal{DB}_{s_0}$ it belongs to is readily identified and  reduced-operators computed for $\widehat{\mubold}$ using the ROMs stored in $\mathcal{DB}_{s_0}$.

\section{Applications and performance assessment}\label{sec:appli}
The proposed approaches are here applied to two challenging physical applications: the acoustic inverse obstacle problem and the aeroelastic flutter problem.

\subsection{Acoustic scattering analysis}\label{sec:scattering}
The acoustic inverse obstacle problem  considered here consists in determining
the shape of an obstacle or a part of this shape from the knowledge of some
scattered far-field patterns, assuming certain characteristics of the surface
of the obstacle. It is well-known \cite{coltonbook} that such an inverse problem is 
non-linear and often quite ill-posed, making its numerical solution
challenging.

To illustrate the ROM database framework proposed above, 
a parameter
identification problem is considered, where the shape of the obstacle is 
assumed to be known a-priori, but the vector of parameters $\mubold\in\mathbb{R}^{N_{\mubold}}$ of the shape
needs to be identified from the measured far-field pattern. This class
of problems is a subset of a general acoustic inverse obstacle problem. 

In order to describe the considered problem more accurately,
the corresponding direct acoustic scattering problem is recalled first.
The scattering of time-harmonic acoustic waves by an impenetrable 
obstacle with the boundary $\Sigma$  embedded in an infinite homogeneous
fluid medium  $ \Omega_e \subset \mathbb{R}^d$
can be formulated as the
following exterior boundary value problem for the unknown acoustic pressure
field $w$
in the fluid
\begin{equation}
\label{eq:scattering}
\begin{array}{r c l c l}
\Delta w +\kappa^2 w & = & 0     & \mbox{in} & \ \ \Omega_e, \\
\displaystyle{\left(a+b\frac{\partial }{\partial \nubold} \right)
 \left( w + w^{inc} \right)} & = & 0 & \mbox{on} & \ \  \Sigma, \\
\displaystyle{ \lim_{r \to \infty} r^{\frac{d-1}{2}}
\left( \frac{\partial w}{\partial r} -j\kappa w \right)} &=& 0,  &&   
\end{array}
\end{equation}
where the incident wave is given by 
$$
w^{inc} = e^{j \kappa \dbold \cdot \xbold},
$$
the unit vector $\dbold$ indicates the direction of the incident plane wave,
and either $a \neq 0$ or $b \neq 0$. Sound-hard, sound-soft, or impedance
boundary conditions can all be represented by the second equation of 
(\ref{eq:scattering}). In the example below, the sound-hard scattering
problem is used, leading to the choice of the Neumann boundary condition
($a=0$ and $b=1$). The third equation in~(\ref{eq:scattering}) is the Sommerfeld radiation
condition. It ensures, in the physical sense, that all waves are outgoing and,
mathematically, that the direct scattering problem is  well-posed for any
wavenumber $\kappa = \omega / c$, where $\omega$ is the angular frequency
of the harmonic oscillations and $c$ is the speed of sound in the fluid.

In order to discretize the direct scattering problem (\ref{eq:scattering}),
the finite element method is considered. The infinite domain is first truncated,
and a perfectly matched layer~\cite{berenger94}
near the exterior boundary is used to simulate the effect of the 
Sommerfeld condition. This converts the boundary value problem 
(\ref{eq:scattering}) into that of solving the algebraic system of
linear equations
\begin{equation}
\label{eq:discrh}
\left( \mathbf{K}(\mubold) - \kappa^2 \mathbf{M}(\mubold) \right) \mathbf{w}(\kappa,\mubold)  = 
\mathbf{f}(\kappa,\mubold)
\end{equation}
for the unknown degrees of freedom $\mathbf{w}\in\mathbb{C}^{N}$.
This system is a second-order LTI system of the form~(\ref{eq:SOLTIPF}). 
Here, $\mathbf{K}$ corresponds to the finite element discretization of
the Laplace operator, and $\mathbf{M}$ is a mass-type matrix. For an
interior problem associated with the Helmholtz equation, $\mathbf{M}$ is real,
symmetric positive-definite, and
$\mathbf{K}$ is real, symmetric non-negative. When the perfectly 
matched layer is used, the matrices $\mathbf{K}$ and $\mathbf{M}$ become
complex and non-Hermitian. The source vector $\mathbf{f}$ arises from 
the discretization of the sound-hard boundary condition.

The far field pattern characterizes the asymptotic behavior of the 
acoustic scattered field far away from the obstacle. In two dimensions,
it admits the following integral representation \cite{coltonbook}
\begin{equation}
\label{eq:ffp}
w_{\infty}(\hat {\xbold} ) \ =
 \ {{e^{j{\pi \over 4}}} \over {(8 \pi \kappa)}^{1 \over 2}}
\int_{\Gamma} \left({\partial w \over \partial \nubold}(\ybold)+
 j\kappa~\hat {\xbold}.\nubold~w(\ybold) \right){e^{-j\kappa {\hat {\xbold}}.\ybold}}  d \sigma_y~;~
\qquad \hat {\xbold} \in S^1, 
\end{equation}
where $S^1$ is the unit circle. After computing the finite-element 
solution, the integral (\ref{eq:ffp}) can be evaluated by integrating
over a suitable curve $\Gamma$ (often the boundary $\Sigma$)
in the computational domain.
The integral in (\ref{eq:ffp}) evaluated at $N_o$ locations $\hat x_1,\cdots,\hat x_{N_o}$ of the circle $S^1$ can then be in practice represented
by the action of a matrix on the solution vector 
$\ybold(\kappa,\mubold) = \mathbf{G}(\kappa,\mubold) \mathbf{w}(\kappa,\mubold)  $ with $\Gbold\in\mathbb{C}^{N_o\times N}$.
Practically, the following 
logarithmic scale quantity is usually plotted
$$
 \mathcal{S}(\hat\xbold) = 10 \log_{10}(2\pi|w_{\infty}(\hat \xbold)|^2).
$$
This quantity can also be computed for each entry of the output vector $\ybold(\kappa,\mubold) $ as
\begin{equation}
\sbold(\kappa,\mubold)  = \sbold\left(\ybold(\kappa,\mubold)\right).
\end{equation}

The inverse problem considered here consists of identifying parameters $\mubold=(L^s,t^s)\in\mathcal{D}\subset\mathbb{R}^2$ of 
a two-dimensional mockup submarine, characterized by its length $L^s$,
and the position of its tower $t^s$ from given
far-field data for several frequencies (wavenumbers)
$\{\sbold_m(\kappa_i)\}_{i=1}^{N_{\kappa}}$.
The inverse problem can be written as
\begin{gather}
\begin{split}
&\min_{\mubold\in\mathcal{D}} \sum_{i=1}^{N_{\kappa}}\alpha_i \left\|  \sbold\left(\ybold(\kappa_i,\mubold)\right) -  \sbold_m(\kappa_i)  \right\|^2_2 + \frac{\beta}{2}\|\mubold\|^2_2\\
&~~\text{s.t.}~~\left( \mathbf{K}(\mubold) - \kappa_i^2 \mathbf{M}(\mubold) \right) \mathbf{w}(\kappa_i,\mubold)  = 
\mathbf{f}(\kappa_i,\mubold),\\
&~~~~~~~~~\ybold(\kappa_i,\mubold) = \mathbf{G}(\kappa_i,\mubold)  \mathbf{w}(\kappa_i,\mubold),~~i=1,\cdots,N_{\kappa},
\end{split}
\end{gather}
where a Tikhonov regularization term has been added and $\{\alpha_i\}_{i=1}^{N_{\kappa}},\beta$ are appropriate positive weights.

Since the solution of the discrete direct problem (\ref{eq:discrh}) for 
each different wavenumber requires a costly re-factorization of the matrix 
on the left-hand side, to enable the efficient computations
for many wavenumbers, a reduced-order model is built using
a derivative-based Galerkin projection (DGP) framework~\cite{hetmaniuk12} for a given value $\mubold$ of the shape parameters.
In this method, based on moment-matching, $N_{\partial}^{\text{DGP}}$ derivatives of the solution
with respect to the wavenumber $\kappa$
are first computed by solving (\ref{eq:discrh}) with recursively constructed
right-hand sides at $N_{\kappa}^{\text{DGP}}$ interpolating wavenumbers. Then, these derivatives
are orthogonalized to achieve numerical robustness, and used to form 
a subspace of dimension $k=N_{\kappa}^{\text{DGP}}N_{\partial}^{\text{DGP}}$ for Galerkin projection, leading to reduced matrices $\Kbold_r(\mubold)$, $\Mbold_r(\mubold)$, $\fbold_r(\kappa,\mubold)$ and $\Gbold_r(\kappa,\mubold)$.

\begin{figure}[ht] \centering
\includegraphics[width=0.45\textwidth]{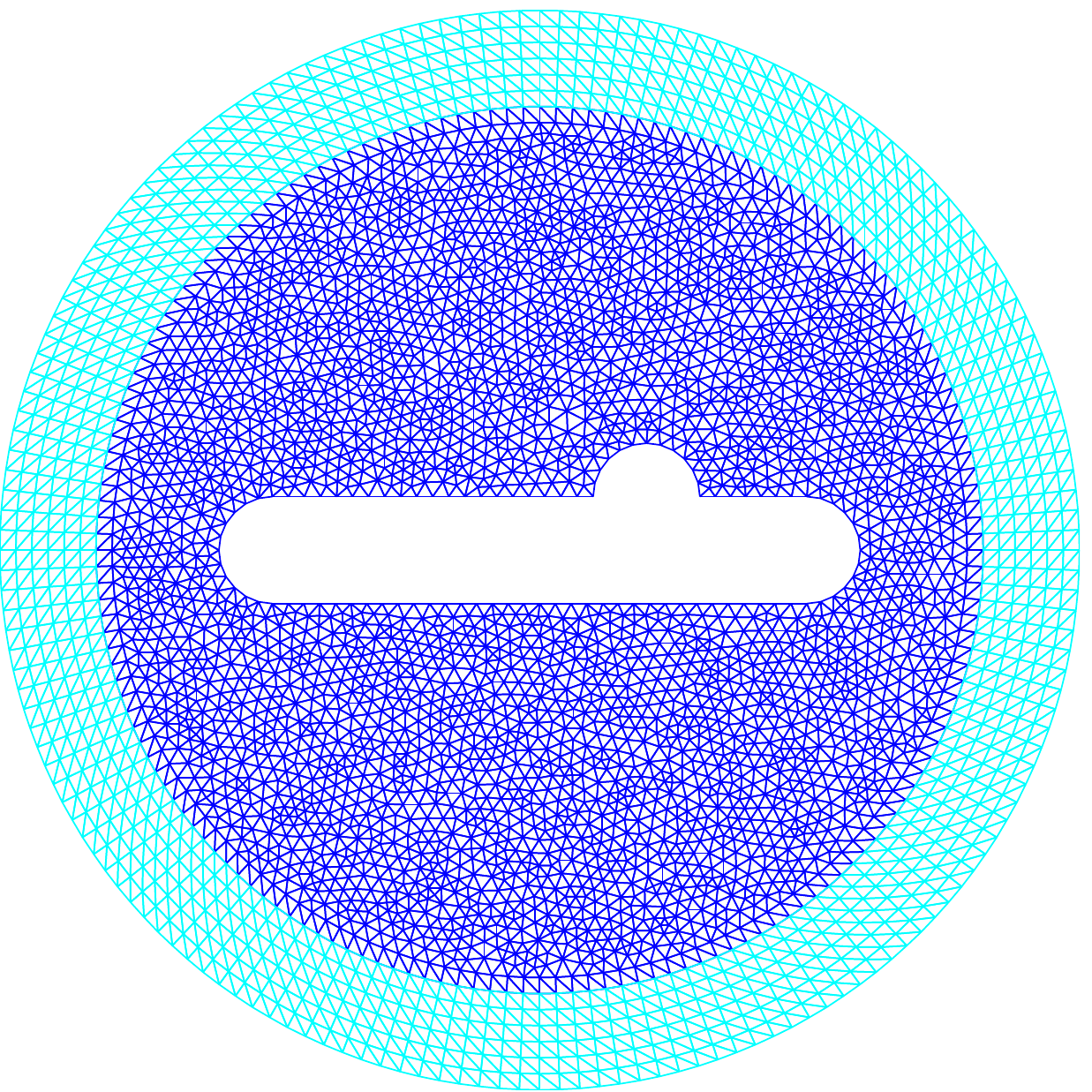}   
\includegraphics[width=0.45\textwidth]{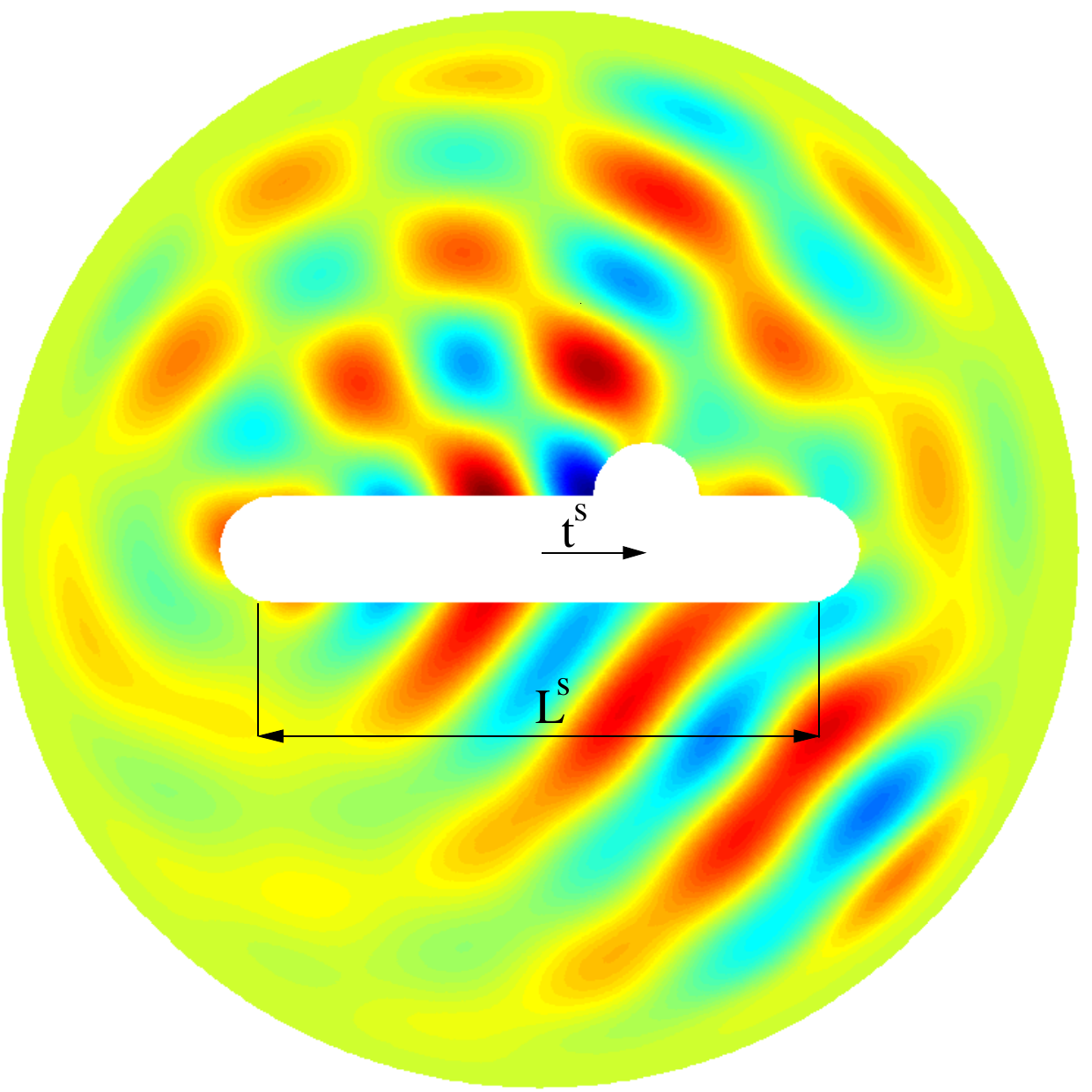}   
\caption{Computational mesh and the real part of the solution for $\kappa=20$}
\label{fig:scattermesh} \end{figure}

Figure \ref{fig:scattermesh} depicts the triangulated computational domain (left) 
with the elements in the PML layer shown in cyan; 
the real part of the solution for $\kappa=20$, $L^s = 1$, and $t^s=0.2$ 
is shown on the right. 
For different values of the shape parameters,
the computational domain is remeshed. Isoparametric cubic Finite Elements are used. All computations are done using Matlab.

The solution procedure described in detail below builds
a database of frequency-sweep ROMs offline by sampling the shape parameter
space adaptively to ensure accuracy. The database of the ROMs is then used to
efficiently solve the reduced inverse problem. For a given value of the parameters, $N_{\partial}^{\text{DGP}}=8$ (including the $0$-th derivative) are computed for the $N_{\kappa}^{\text{DGP}}=2$ frequencies $\kappa\in\{10,20\}$, leading to a ROM of dimension $k=16$. In the present case, there are $N_0=360$ outputs equidistributed on the sphere $S^1$.

In a first set of numerical experiments, the effect of consistency on predictions based on a database of ROMs is illustrated. For that purpose, a small database of $N_{\text{DB}}=4$ ROMs is first created for the values of the parameters indicated in Table~\ref{tab:smallDB} where the number of dofs for each underlying HDM is also reported. One can observe that each HDM has a different number of dofs. 

The consistency enforcement procedure for arbitrary meshes developed in Section~\ref{sec:ArbitraryMesh} is first applied to transform the reduced operators prior to their interpolation. In a second case, the operators are not transformed. In all cases, the reduced operators are interpolated on their appropriate manifold at $\mubold^\star = [0.9625,0.1125]$ by bilinear interpolation as follows:
\begin{itemize}
\item The operators $\text{Re}(\Mbold_{r}(\cdot))$, $\text{Im}(\Mbold_{r}(\cdot))$ and $\text{Re}(\Kbold_{r}(\cdot))$ are SPD matrices and as such as interpolated on their appropriate manifold using the Choleski decomposition-based approach proposed in Section~\ref{sec:interp}.
\item The operator  $\text{Im}(\Kbold_{r}(\cdot))$ is interpolated on the manifold of symmetric matrices.
\item The operators $\{\fbold_r(\kappa_i,\cdot),\Gbold_r(\kappa_i,\cdot)\}_{i=1}^{N_{\kappa}}$ are interpolated on the manifold of rectangular complex matrices.
\end{itemize}

\begin{table}[htdp]
\begin{center}
\begin{tabular}{|c|c|c|c|}
\hline $\mathcal{DB}$ point  & $\mu_1$ & $\mu_2$ & $N(\mubold)$ \\
\hline 1 & 0.95 & 0.1 & 41,235 \\
\hline 2 & 0.975 & 0.1 & 40,965 \\
\hline 3 & 0.95 & 0.125 & 41,424\\
\hline 4 & 0.975 & 0.125 & 40,929\\
\hline 
\end{tabular} 
  \caption{Database of ROMs and associated underlying HDM number of dofs}
  \label{tab:smallDB}
\end{center}
\end{table}

\begin{figure}[htbp] \centering
\includegraphics[width=0.75\textwidth]{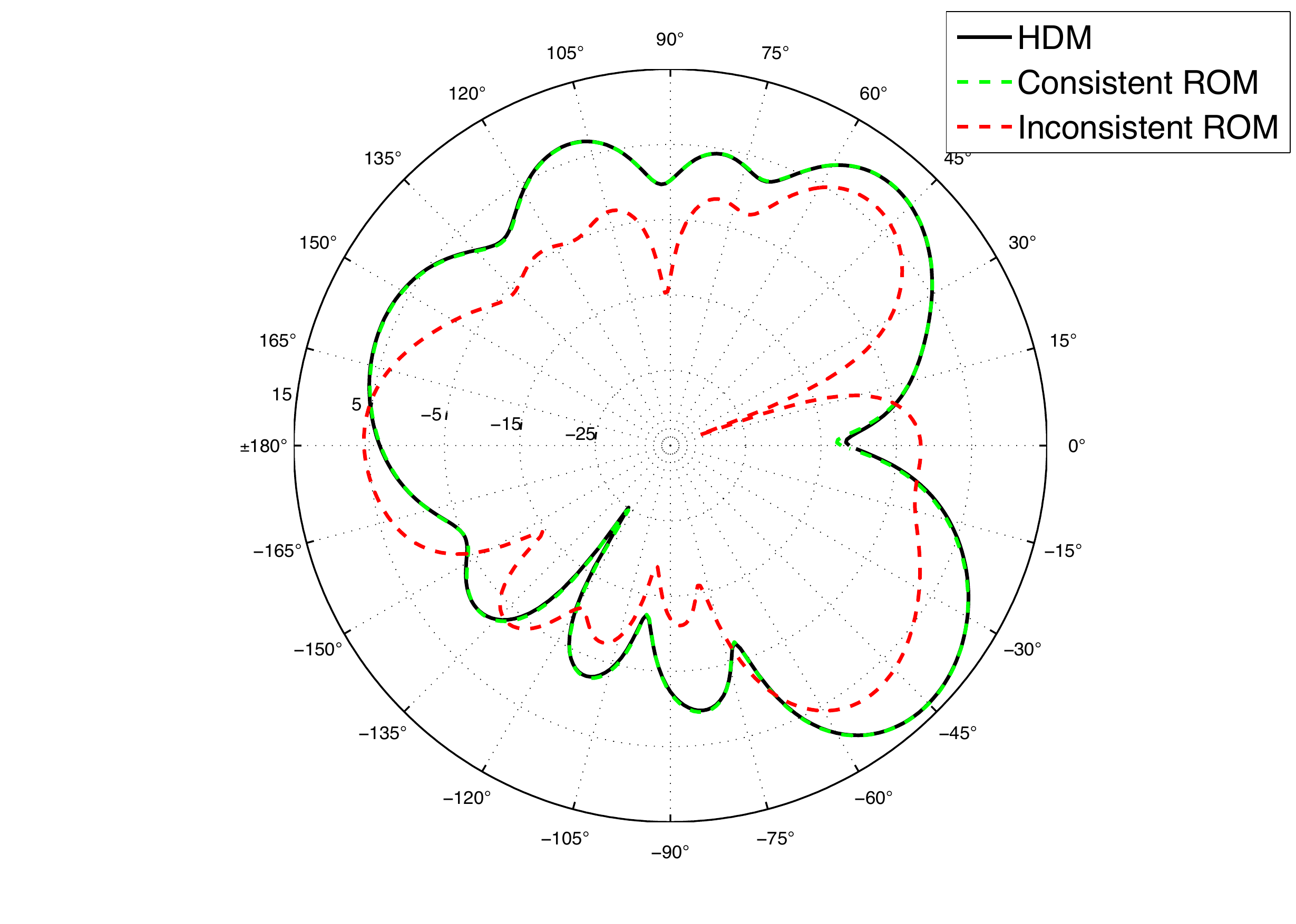}   
\caption{Comparison of the far-field pattern predictions obtained by interpolation of inconsistent and consistent databases.}
\label{fig:consistency_effect} \end{figure}

Figure~\ref{fig:consistency_effect} reports the far-field pattern for $\kappa=14$ at $\mubold^\star = [0.9625,0.1125]$ computed by the HDM and the two ROM interpolation approaches. One can observe the importance of consistency as the inconsistent ROM database leads to poor predictions whereas the consistent ROM database predictions very closely follow the HDM results.

Next an adaptive approach for constructing the ROM database is developed. The approach proceeds by comparing the predictions associated with the ROM database and HDM at the center of each hypercube of the database and refining that hypercube if the error is above a given threshold. The error can be associated with an output of interest, such as the far-field pattern, leading to a goal-oriented approach. In the following, the error will be based on the accuracy of the shape $\widehat{\mubold}$ returned by solving the reduced inverse problem by simulated annealing.
\begin{gather}
\begin{split}
&\min_{\mubold\in\mathcal{D}} \sum_{i=1}^{N_{\kappa}}\alpha_i \left\|  \sbold\left(\ybold_r(\kappa_i,\mubold)\right) -  \sbold_m(\kappa_i)  \right\|^2_2 + \frac{\beta}{2}\|\mubold\|^2_2\\
&~~\text{s.t.}~~\left( \mathbf{K}_r(\mubold) - \kappa_i^2 \mathbf{M}_{r}(\mubold) \right) \mathbf{q}(\kappa_i,\mubold)  = 
\mathbf{f}_{r}(\kappa_i,\mubold),\\
&~~~~~~~~~\ybold_r(\kappa_i,\mubold) = \mathbf{G}_r(\kappa_i,\mubold)  \mathbf{q}(\kappa_i,\mubold),~~i=1,\cdots,N_{\kappa}. 
\end{split}
\end{gather}
 The error measure is then defined as
 \begin{equation}
 \text{Error} = \left\|\frac{\widehat{\mubold}-\mubold^\star}{\mubold_{\max}-\mubold_{\min}}\right\|_\infty
 \end{equation}
where $\mubold^\star$ denotes the target shapes, $\mubold_{\min}$ and $\mubold_{\max}$ the lower and upper bounds for each shape parameters, respectively and the ratio in $\text{Error}$ is computed entry-by-entry.

\begin{figure}[htbp] \centering
\includegraphics[width=0.75\textwidth]{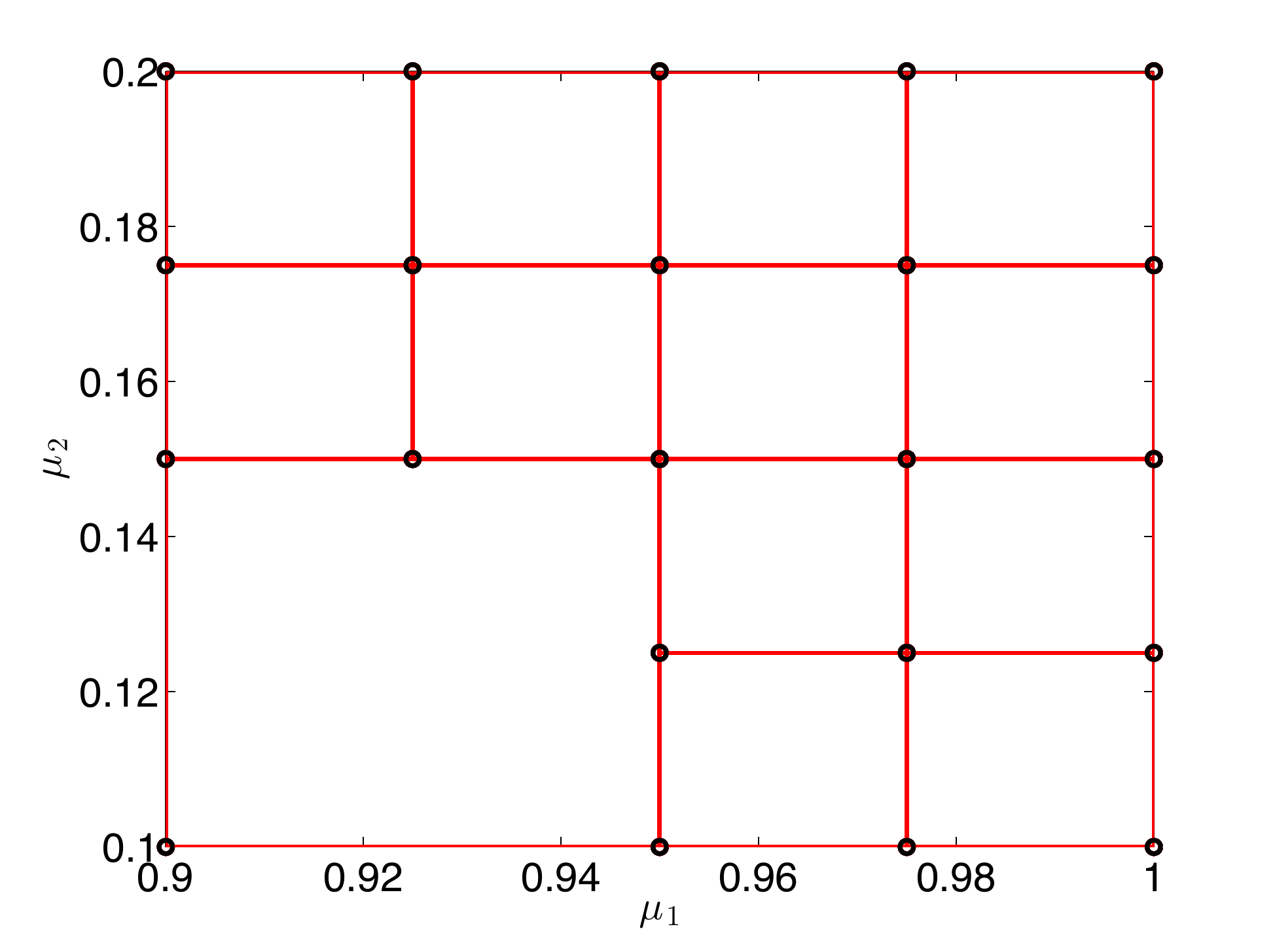}   
\caption{Pre-computed database $\mathcal{DB}$ of reduced-order models selected by the adaptive sampling procedure.}
\label{fig:scatter1} \end{figure}

\begin{figure}[htbp] \centering
\includegraphics[width=0.75\textwidth]{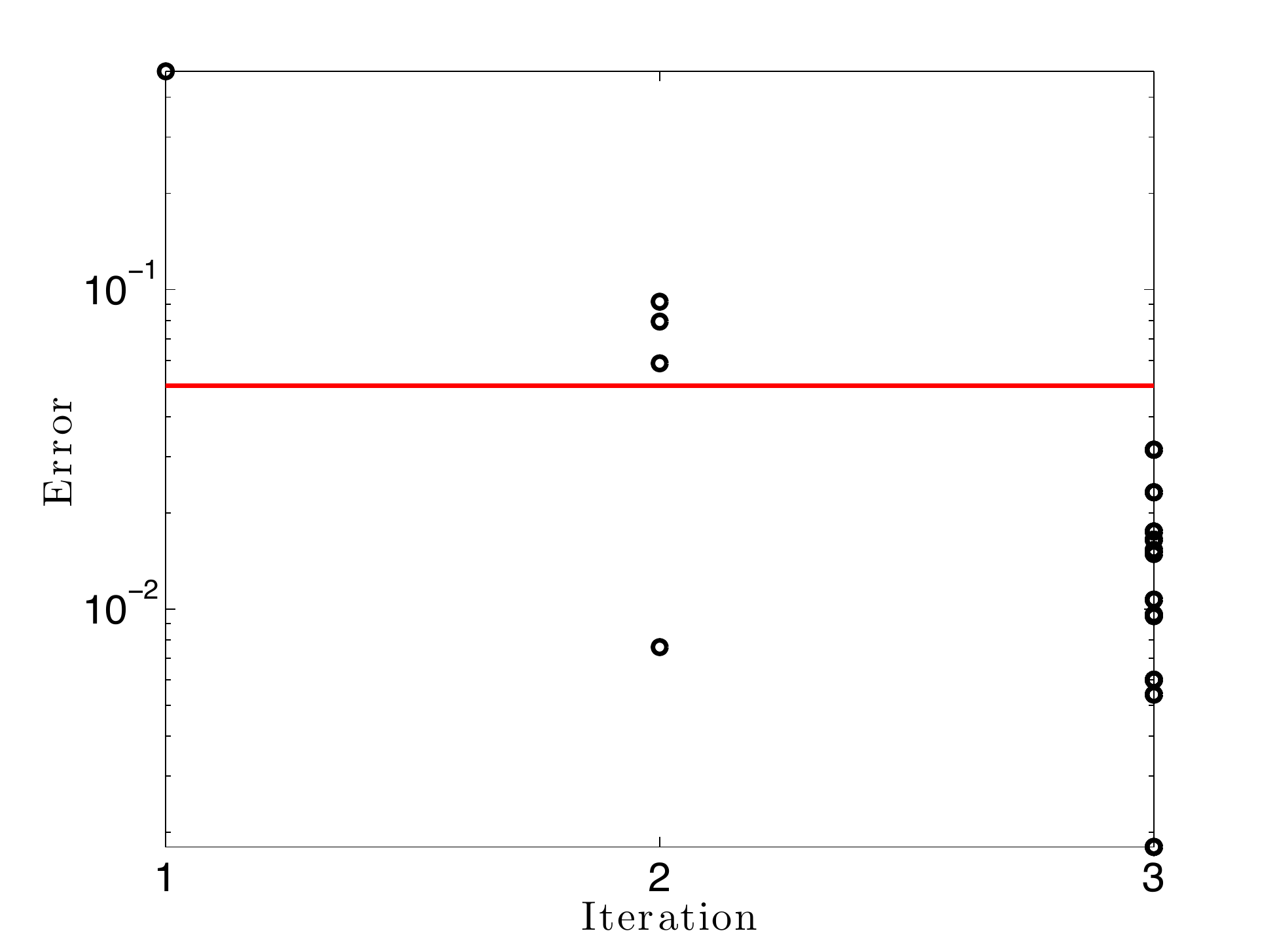}   
\caption{Training errors at each iteration of the adaptive sampling procedure.}
\label{fig:scatter2} \end{figure}

Figure~\ref{fig:scatter1} reports the refined database for an error tolerance of $5\%$ and the parameter domain $\mathcal{D} = [0.9,1]\times[0.1,0.2]$. $N_{\text{DB}}=21$ points are sampled in the domain. The training errors obtained at each iteration refinement of the procedure are reported in Figure~\ref{fig:scatter2}. One can observe that after the second refinement, all errors are below the error threshold of $5\%$.

\begin{figure}[htbp] \centering
\includegraphics[width=0.75\textwidth]{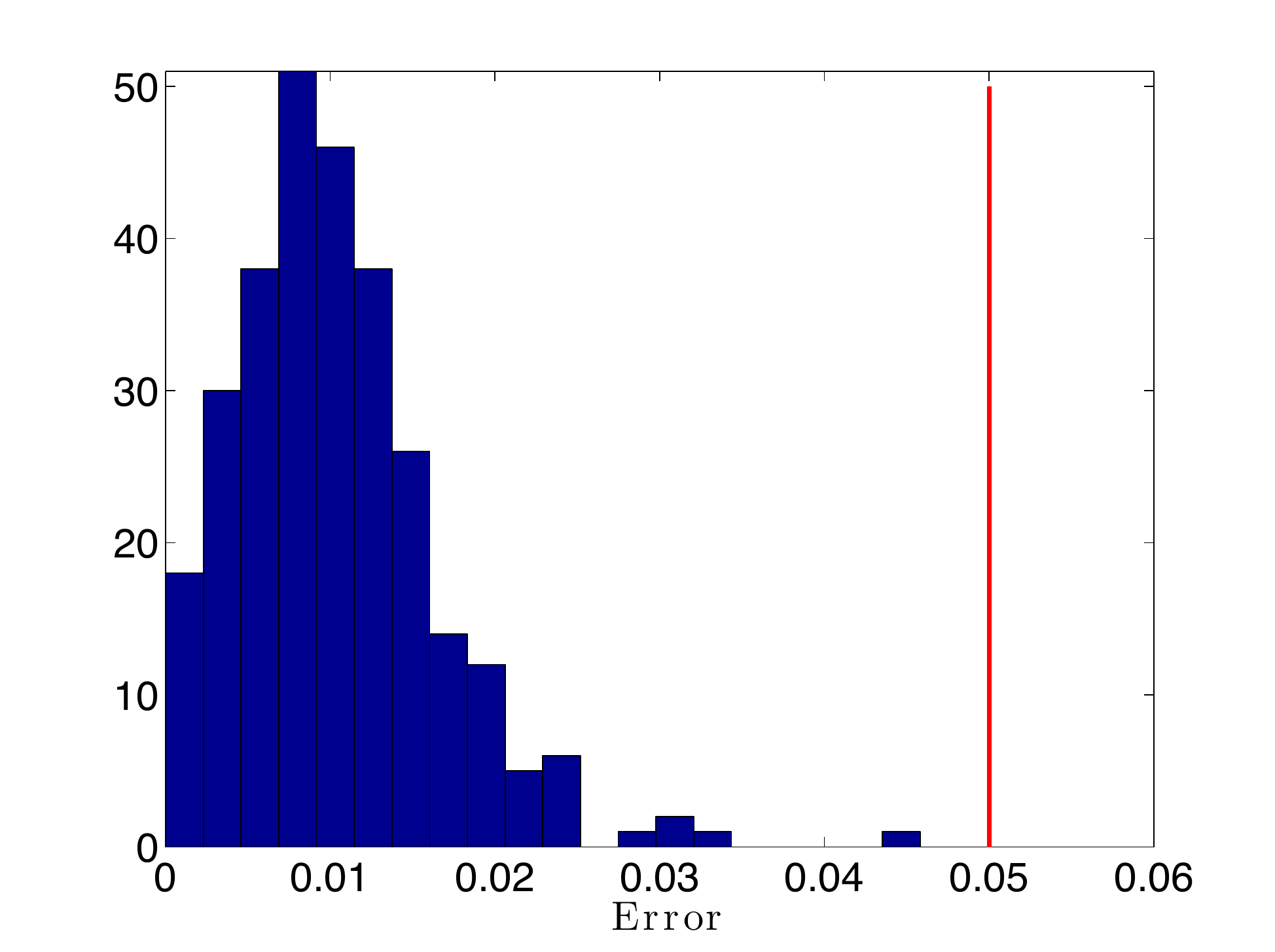}   
\caption{Distribution of prediction errors for the solution of $289$ inverse problems.}
\label{fig:scatter3} \end{figure}

\begin{figure}[htbp] \centering
\includegraphics[width=0.75\textwidth]{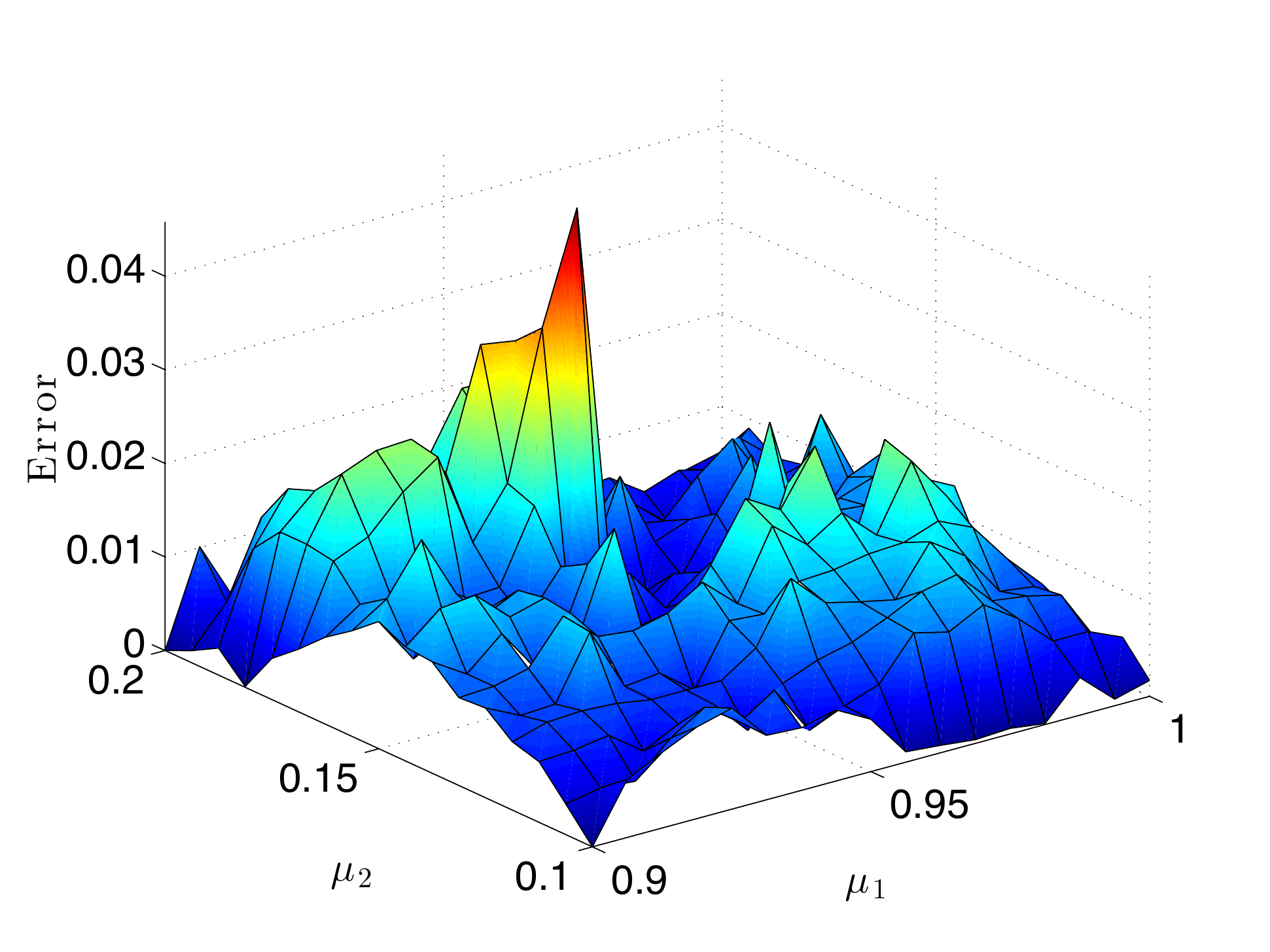}   
\caption{Prediction errors for the solution of $289$ inverse problems.}
\label{fig:scatter4} \end{figure}

To validate the accuracy of the ROM database, $289$ target shape parameters are selected in $\mathcal{D}$ and the reduced inverse problems solved for each of them. The distribution of corresponding errors are reported in Figures~\ref{fig:scatter3} and~\ref{fig:scatter4}. One can observe that all errors are below the error threshold, confirming the validity of the training procedure.

Finally, the CPU timings associated with the solution of a given inverse problem are compared for the HDM and ROM database strategy. The CPU timings are reported in Table~\ref{tab:scatteringOnline}. One can observe that an impressive speedup of $270$ is obtained thanks to the database strategy. For a given function call in the optimization problem, the speedup is equal to $207$.

\begin{table}[htdp]
\begin{center}\begin{tabular}{|c|c|c|c|c|}\hline Approach  & Error & Number of  & Online optimization  & Speedup  \\    & &function calls & Wall time & \\ \hline\hline HDM & $8\times10^{-4}$ & 1530 & 1 h 30 min & 1  \\\hline ROM database & $0.02$ & 1176 & 20.1 s & 270  \\\hline \end{tabular}
\end{center}
  \caption{Wall times associated with the solution of the inverse problem with the HDM and database of ROMs}
  \label{tab:scatteringOnline}
\end{table}


\subsection{Flutter analysis}\label{sec:flutter}
The aeroelastic analysis of a wing-store configuration flying in the subsonic, transonic and supersonic regimes is considered. Some properties of that system were originally studied in~\cite{chiu09,farhat13}. Among those, the hydroelastic effects inside the fuel tank modify the structural properties of the wing-store configuration, thereby affecting the flutter characteristic of the system. This justifies parameterizing the aeroelastic system by the fuel fill level $f$ inside the tank. Furthermore, the aeroelastic properties of the wing-store system depend heavily on the aerodynamic properties of the configuration. As a result, the system will also be parameterized by the free-stream Mach number $M_\infty$.

\begin{figure}
\begin{subfigmatrix}{2}
\subfigure[CFD surface grid]{\includegraphics{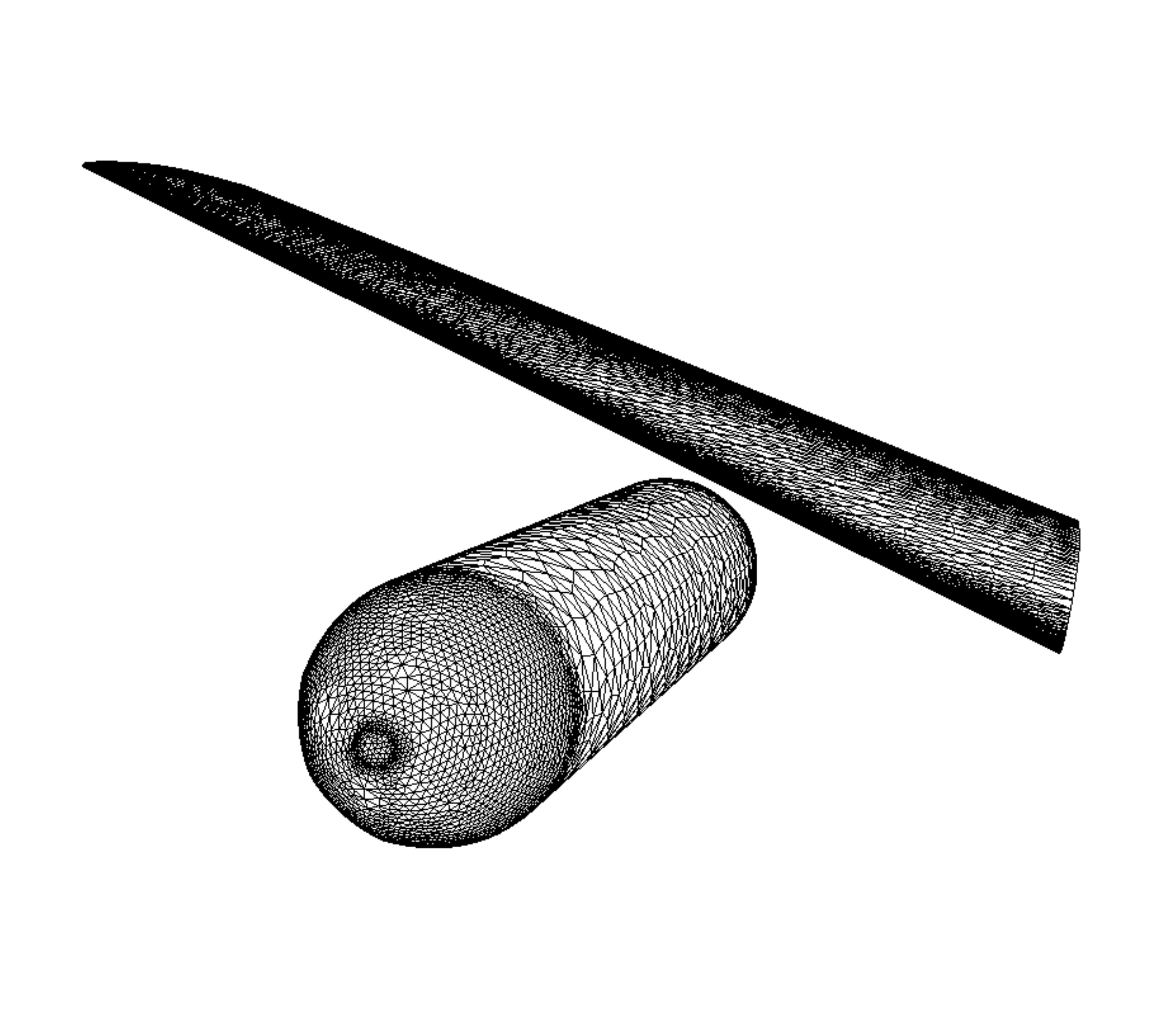}}
\subfigure[FE structural model]{\includegraphics{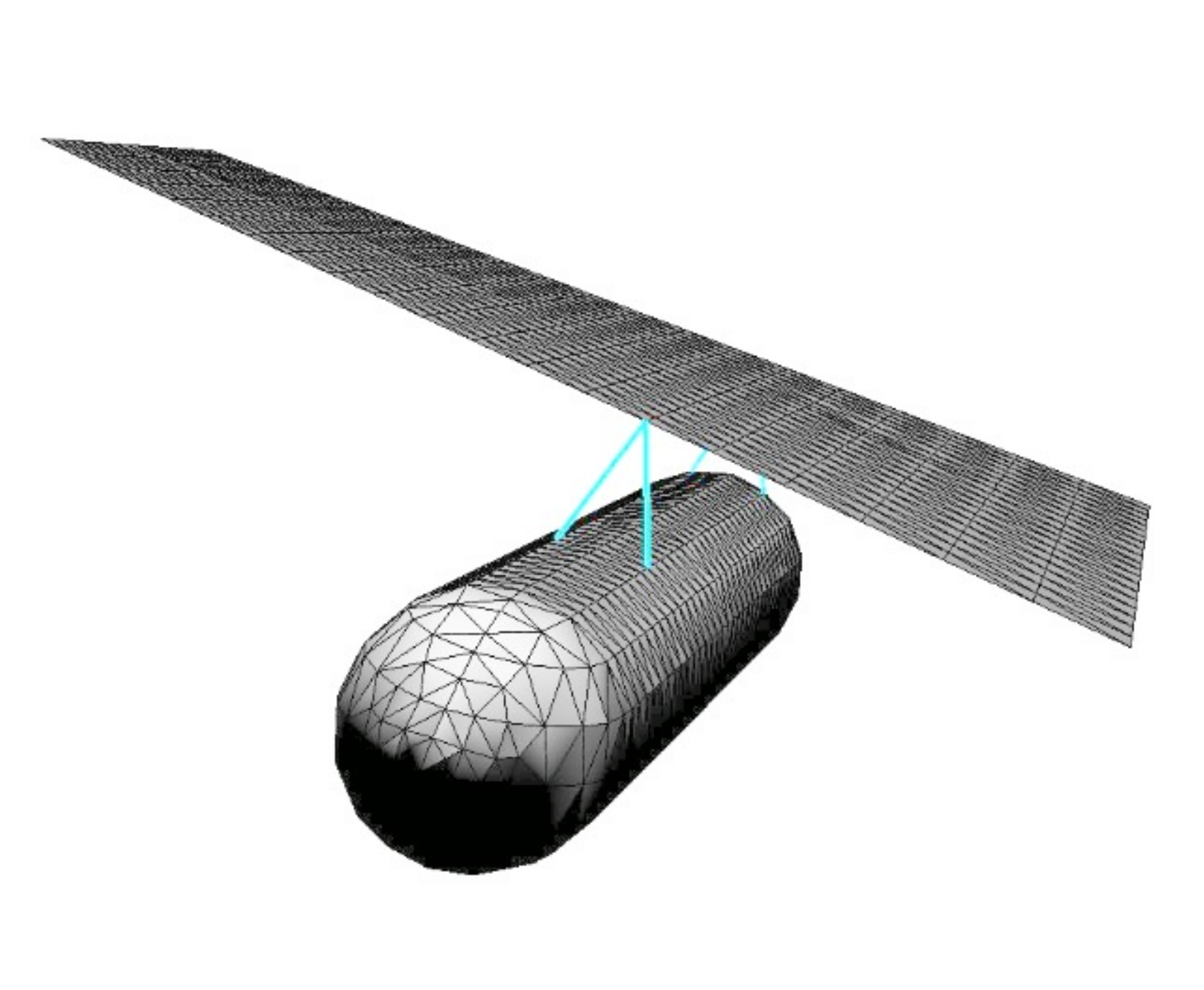}}
\end{subfigmatrix}
\caption{High-dimensional aeroelastic model of a wing-store configuration.}
\label{fig:AGARD_meshes}
\end{figure}

The structural and fluid surface models of the wing-store system are graphically depicted in Figure~\ref{fig:AGARD_meshes}.  The structural subsystem is a second-order LTIP of the form~(\ref{eq:SOLTIP}) and is modeled by the Finite Element method. For each fill level, a new structural mesh is generated inside the tank for the full physical domain. The hydroelastic effects are modeled by an added mass effect~\cite{chiu09,farhat13}, resulting in a linear HDM with $N^{(s)} = 6,834$ dofs for all values of $f$. The proposed ROM database approach will enable by-passing the re-meshing of the fuel domain everytime the fill level is varied.

The fluid subsystem is modeled by the linearized Euler equations and discretized by the Finite Volume method using a second-order accurate linear flux reconstruction and a second-order accurate implicit backward difference time integration scheme. For each operating point $\mubold = (M_\infty, f)$, the nonlinear fluid HDM is linearized around a steady-state, resulting in a first-order LTIP system of the form~(\ref{eq:FOLTIP}) with $N^{(f)}\approx 400,000$ dofs.

In this work, the operation domain of interest is $(M_\infty,f)\in \mathcal{D} = [0.6,1.1] \times [0,100]$. For each operating point, the critical values of pressure $p_\infty^{\text{cr}}$ and velocity $V_\infty^{\text{cr}}$ at the onset of flutter are sought. Once, these quantities are determined, the 
 flutter speed index (FSI) can be computed as
 \begin{equation}
 \text{FSI} = \frac{V^{\text{cr}}_\infty}{b_s \omega_\alpha \sqrt{\bar{\mu}}},
 \end{equation}
where $b_s$ is the semi-chord of the wing at its root, $\omega_\alpha$ is the first dry torsional mode of the wing-store structural system and $\bar{\mu}$ is the mass ratio as defined in~\cite{yates87,chiu09,farhat13}.

\begin{figure}[htbp] \centering
\includegraphics[width=0.45\textwidth]{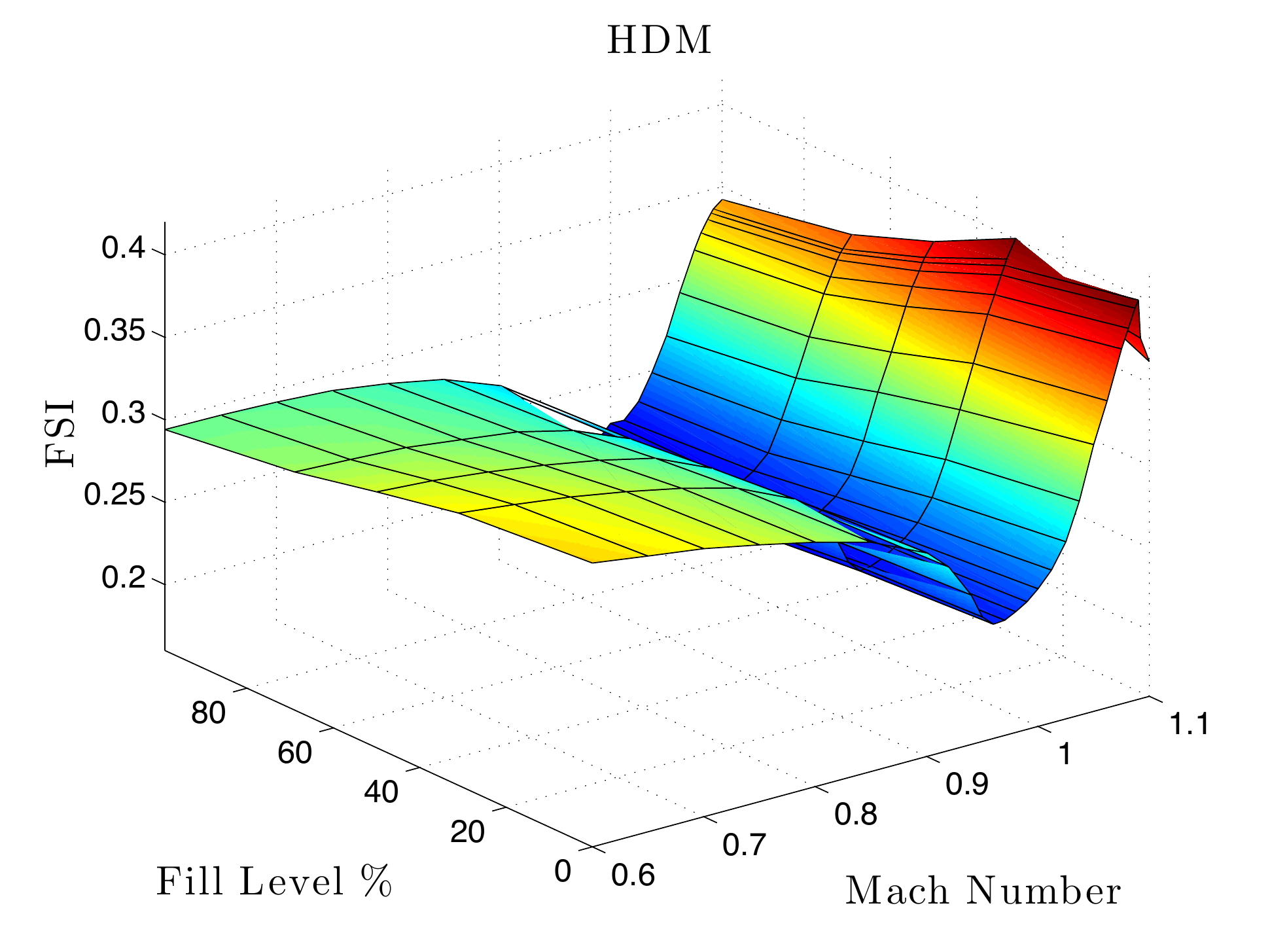}   \includegraphics[width=0.45\textwidth]{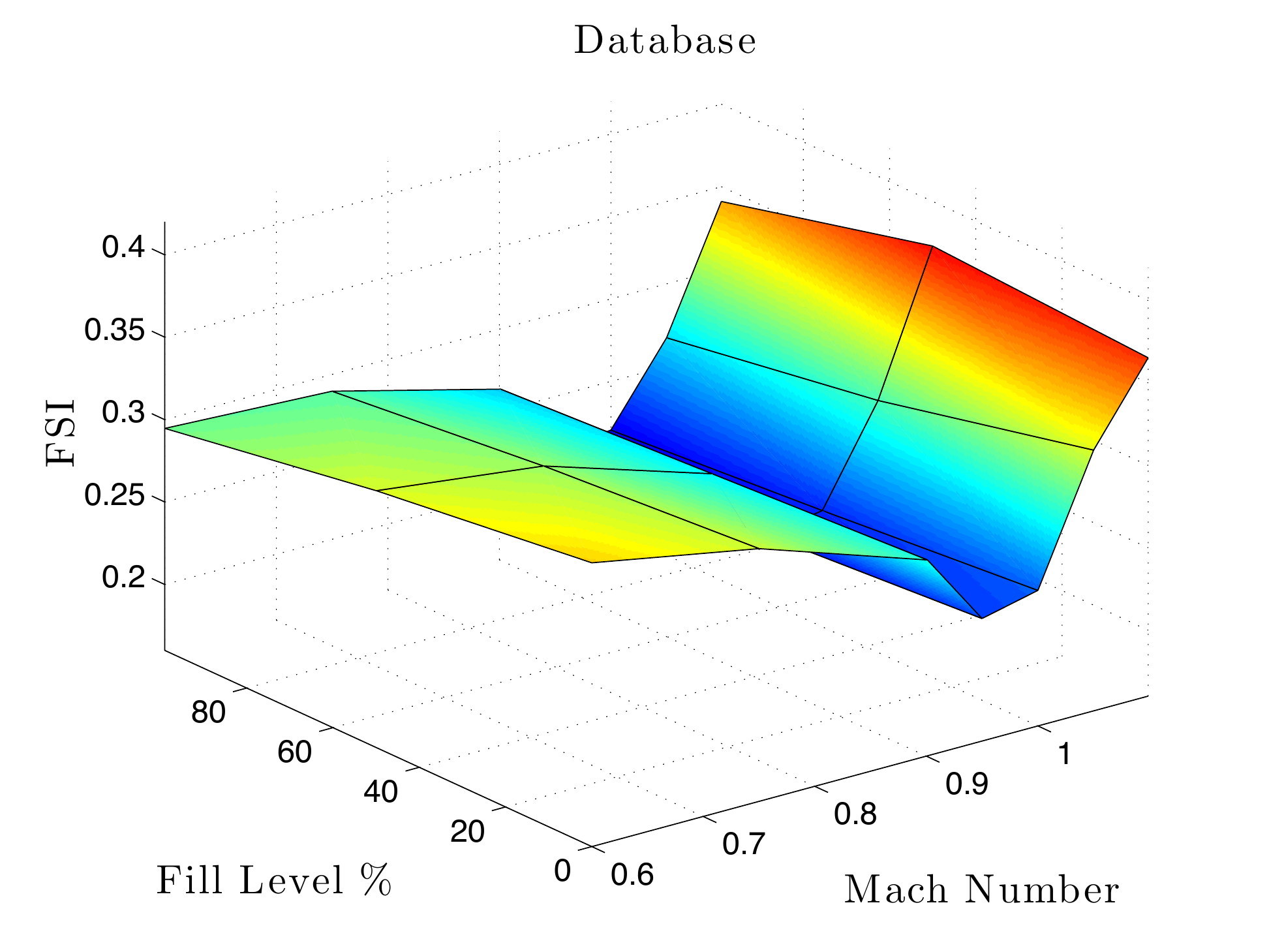}   
\caption{Comparison of the high-dimensional model and database reduced-order models flutter speed indices.}
\label{fig:FSI_DB} 
\end{figure}

The flutter speed indices of the system of interest are computed using the HDM for 26 different free-stream Mach number and 5 different fill levels in the domain $(M_\infty,f)\in \mathcal{D}$, resulting in 130 operating points,  and reported in Figure~\ref{fig:FSI_DB}. One can observe the characteristic flutter dip for $M_\infty\approx 0.96$.

The framework developed in this paper is then applied to the problem of interest to interpolate reduced aeroelastic operators. In this example, all structural and fluid HDMs are defined on the same mesh and the approach of enforcing ROM consistency developed in Section~\ref{sec:CommonMesh} is followed. For that purpose, $N_{\text{DB}}=21$ operating points are sampled and their corresponding aeroelastic ROMs constructed and stored in the offline phase in a database $\mathcal{DB}$. These points correspond to a lattice $(M_\infty,f) \in \{0.6,0.75,0.9,0.95,1.0,1.05,1.1\}\times \{0,50,100\}$. For each operating point, a structural ROM of dimension $k^{(s)}=4$ is constructed by projecting the linear structural HDM onto its first four natural modes.  Then a fluid ROM of dimension $k^{(f)}=15$ is constructed by POD using the approach described in~\cite{amsallem08,amsallem10,amsallemthesis} and a reduction of the system in descriptor form~\cite{amsallem14:book}. The FSIs predicted by those 21 aeroelastic ROMs are reported in the right portion of Figure~\ref{fig:FSI_DB}. Very good agreements can be observed at the database points when compared to their counterparts determined by the HDM that are depicted in the left portion of  that same figure.

\begin{figure}[htbp] \centering
\includegraphics[width=0.75\textwidth]{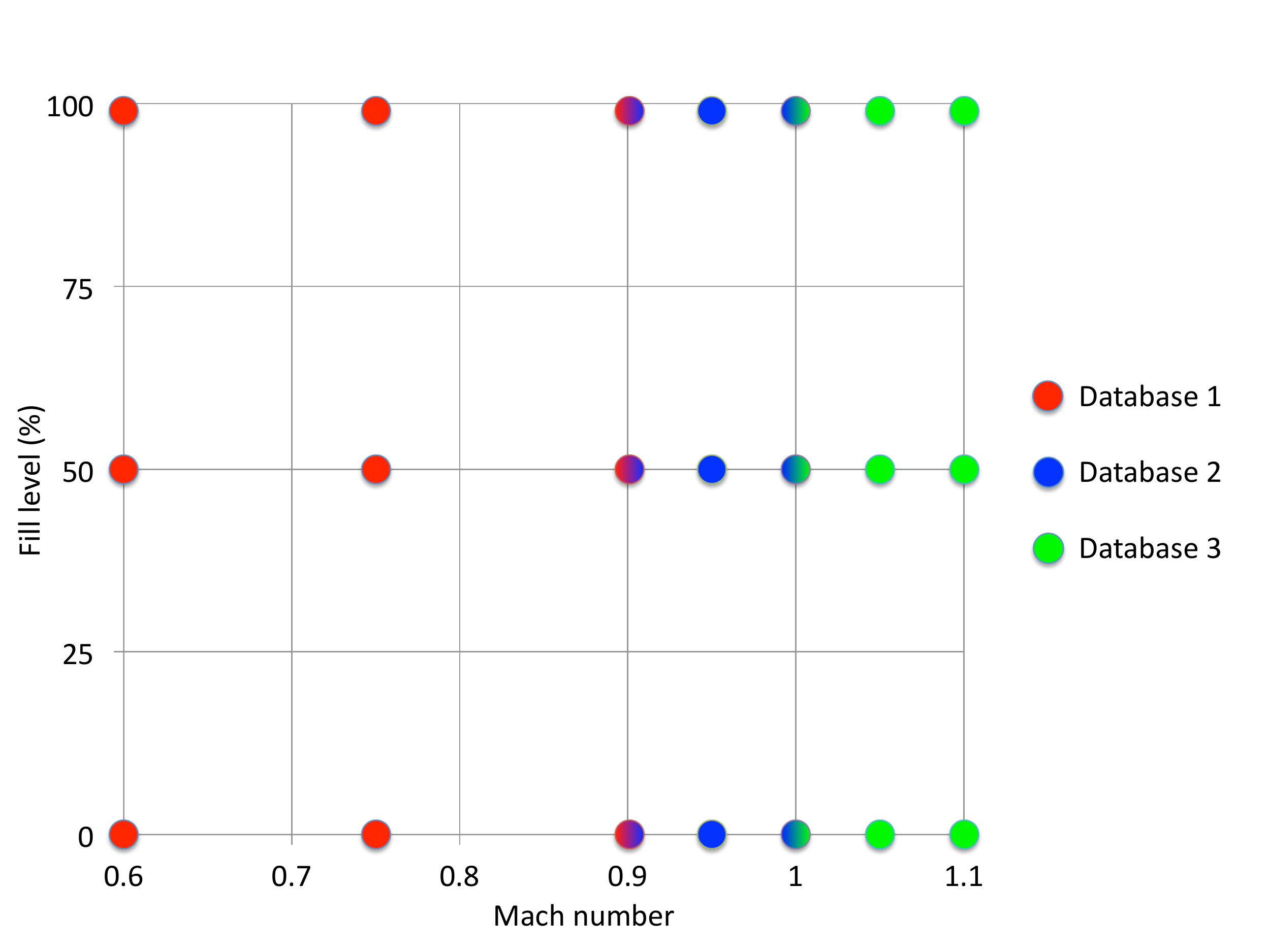}   
\caption{Pre-computed sub-databases $\{\mathcal{DB}_s\}_{s=1}^3$ of aeroelastic reduced-order models.}
\label{fig:DB} \end{figure}

The 21 pre-computed aeroelastic ROMs are then distributed in $N_s=3$ sub-databases: the first one covers the subsonic and lower transonic flow regime $M_\infty \in [0.6,0.9]$, the second one the upper transonic regime $M_\infty \in [0.9,1.0]$ and the third one the supersonic regime $M_\infty \in [1.0,1.1]$. These three databases are graphically depicted in Figure~\ref{fig:DB}.
 In the online interpolation procedure, in each database, piecewise-linear interpolation will be used in the $M_\infty$ direction and cubic spline interpolation in the $f$ direction.

Before the online phase, the heuristic developed in~\cite{degroote10} and mentioned in Section~\ref{sec:interp} is applied for choosing the manifold onto which to interpolate the reduced fluid operators. For the fluid operators, interpolation can indeed be done on the manifold $\text{GL}\left(k^{(f)}\right)$ of non-singular matrices of size $k^{(f)}$ or on $\mathbb{R}^{k^{(f)}\times k^{(f)}}$. Since the interpolation procedure of choice involves two points in the $M_\infty$ direction and three in the $f$ direction at a time, the heuristic is applied for six different regions of the parametric space, as indicated in Table~\ref{tab:heuristic}. As reported in Table~\ref{tab:heuristic}, the manifold $\text{GL}\left(k^{(f)}\right)$ is chosen in two regions while the manifold $\mathbb{R}^{k^{(f)}\times k^{(f)}}$ is chosen in four regions. 
   
   \begin{table}[htdp]
\begin{center}
\begin{tabular}{|c|c|c|c|c|c|c|}
\hline 
& \multicolumn{2}{c|}{Database 1} & \multicolumn{2}{c|}{Database 2} & \multicolumn{2}{c|}{Database 3} \\
\hline
$M_\infty\in$ &$[0.6,0.75]$  & $[0.75,0.9]$ & $[0.9,0.95]$ & $[0.95,1]$& $[1,1.05]$& $[1.05,1.1]$ \\
 \hline 
$\text{GL}\left(k^{(f)}\right)$ & 0.89 & 0.87 & 0.93 & 0.93& 0.76&0.86\\
 \hline 
$\mathbb{R}^{k^{(f)}\times k^{(f)}}$ & 0.84 & 0.86 & 0.98 &0.87 & 0.78&0.82\\
 \hline 
Manifold choice & $\mathbb{R}^{k^{(f)}\times k^{(f)}}$ &  $\mathbb{R}^{k^{(f)}\times k^{(f)}}$&  $\text{GL}\left(k^{(f)}\right)$ &  $\mathbb{R}^{k^{(f)}\times k^{(f)}}$& $\text{GL}\left(k^{(f)}\right)$ & $\mathbb{R}^{k^{(f)}\times k^{(f)}}$\\
 \hline
  \end{tabular} 
  \caption{Non-linearity indicator in each database for the manifold choice heuristic}\label{tab:heuristic}
\end{center}

\end{table}     
   \begin{figure}[htbp] \centering
\includegraphics[width=0.65\textwidth]{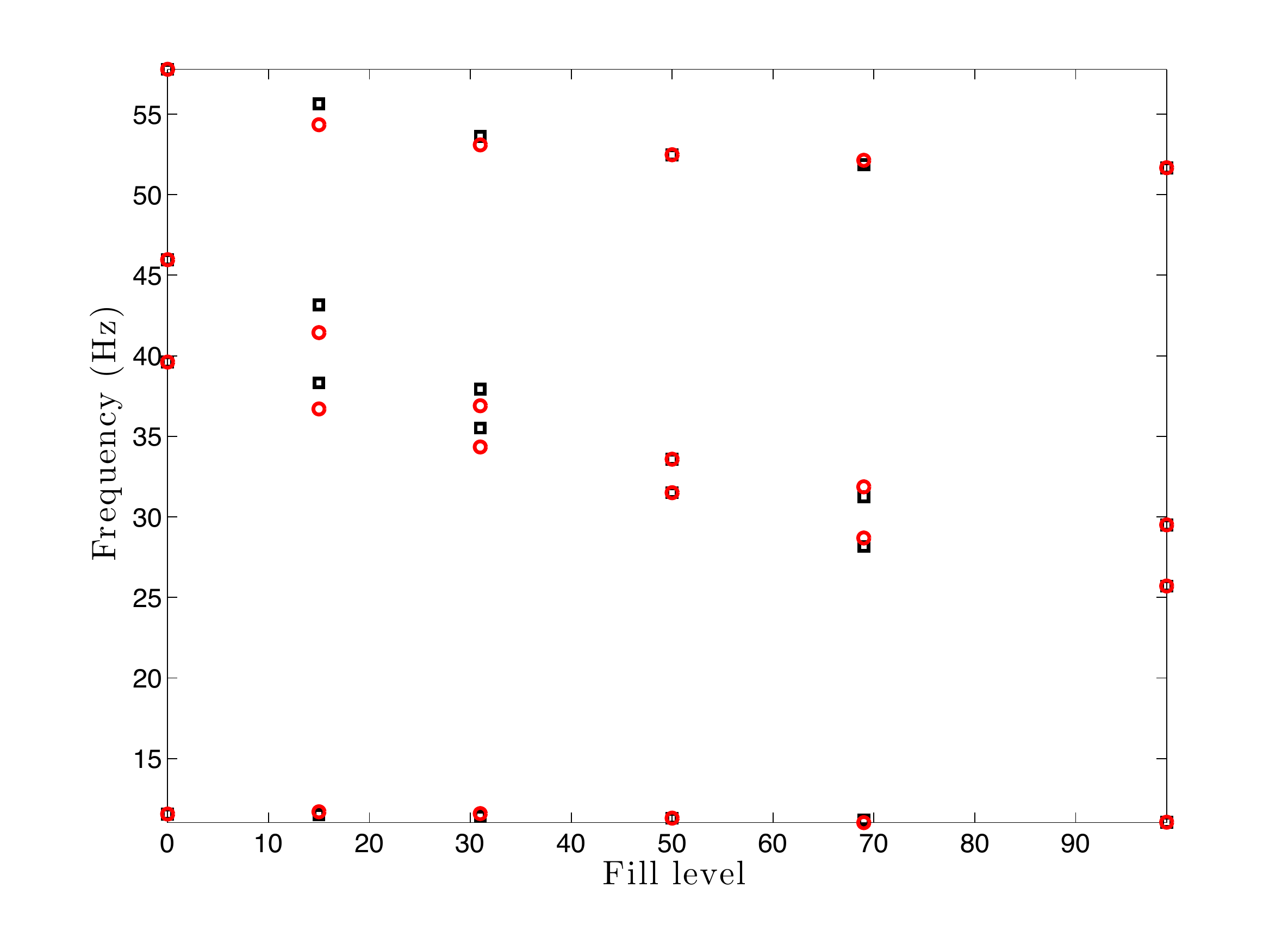}   
\caption{Comparison of eigenvalues of the structural subsystem: HDM $\mathbf{\square}$,  interpolated ROM  {\color{red} $\mathbf\circ$}.}
\label{fig:struc_ev} \end{figure}

\begin{figure}[htbp] \centering
 \includegraphics[width=0.45\textwidth]{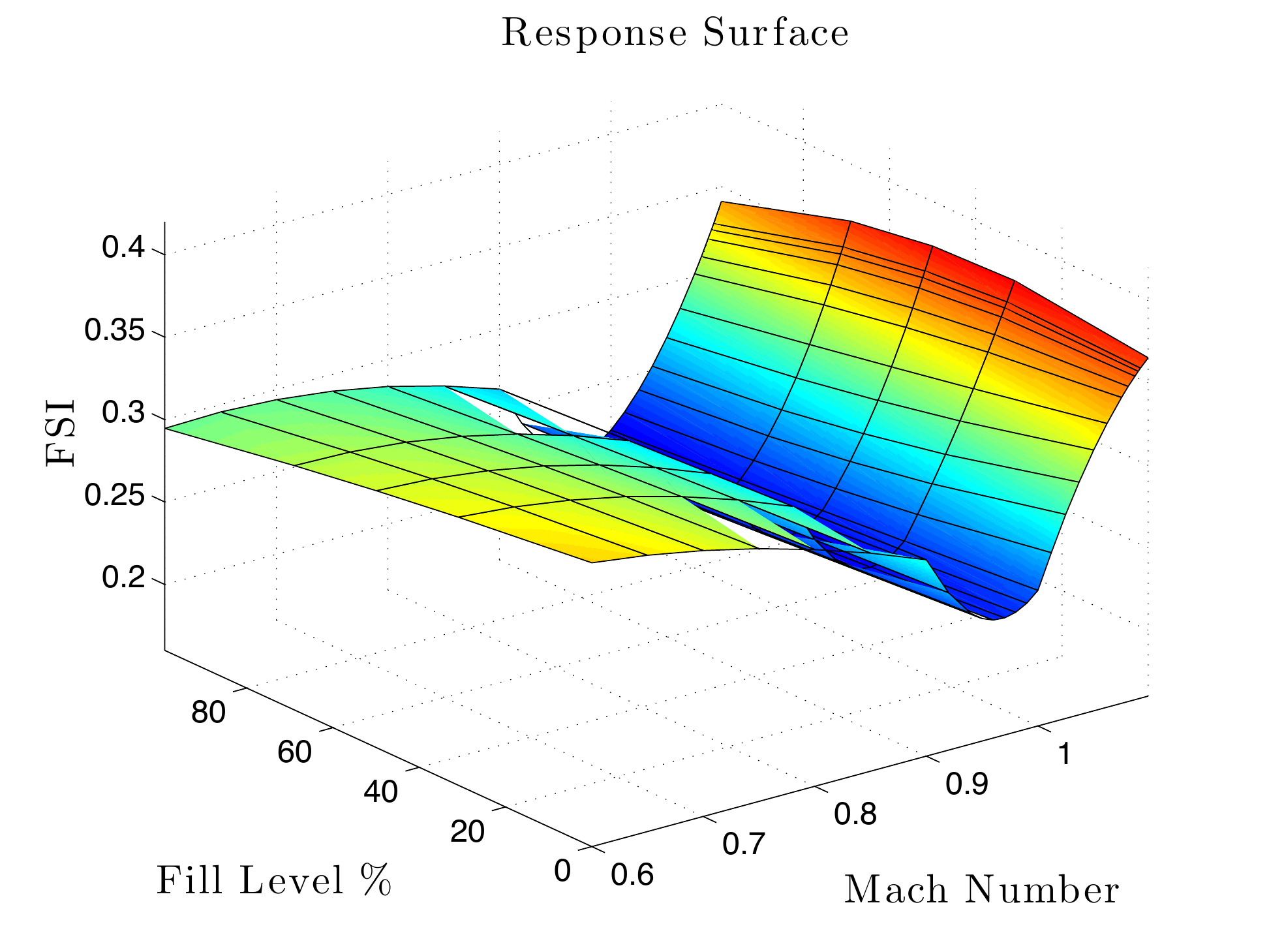}       
  \includegraphics[width=0.45\textwidth]{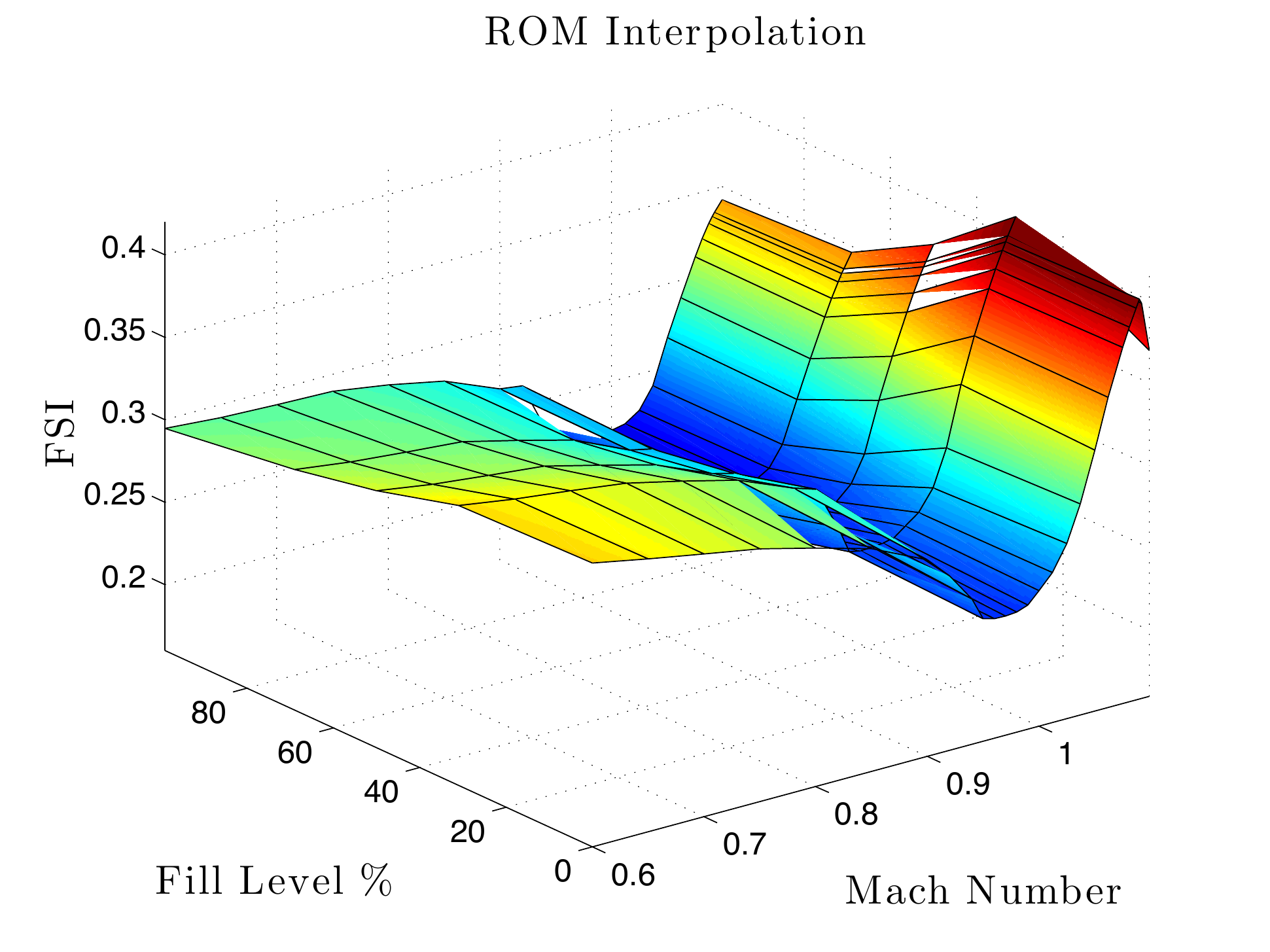}    
\caption{Comparison of  flutter speed indices predicted using (1) response surface estimation, (2) ROM interpolation.}
\label{fig:FSI_predictions} 
\end{figure}

Next, the proposed  methodology is applied to interpolate the aeroelastic ROM operators. The properties of the structural operators resulting from that interpolation are first analyzed by comparing their respective eigenfrequencies to their HDM counterparts. The corresponding results are reported in Figure~\ref{fig:struc_ev}. Good agreements can be observed, even for fill levels that are not present in the database. Next, the interpolated aeroelastic ROM operators are used to predict the onset of flutter in the entire parametric domain $(M_\infty,f) \in [0.6,1.1] \times [0,100]$. The predicted FSIs are reported in Figure~\ref{fig:FSI_predictions}. Very good qualitative and quantitative agreement can be observed for all flight conditions considered. For comparison, response surface estimation (RSE) is also applied to predict flutter using the database FSI data reported in the right portion of Figure~\ref{fig:FSI_DB}. In this case, bicubic spline interpolation is used. When compared with the predictions arising from ROM interpolation, the results from RSE are found to be much less accurate, especially near the transonic dip and in the supersonic regime. RSE cannot, in particular, predict the FSI behavior for low fill levels at supersonic speeds. It is quite remarkable that the method proposed in this paper can capture this complex phenomenon with  only the ROM database associated with the results shown in the right portion of Figure~\ref{fig:FSI_DB}. This example underlines the potential for accurate predictions of the proposed method which operates by interpolating models and not outputs, as in RSE. The offline and online CPU times associated with the prediction of the entire parametric FSI surface each of the four techniques are reported in Table~\ref{tab:CPUtimes}. These results clearly demonstrate the real-time capability of the proposed approach, as it can accurately predict the FSI for $130$ configurations in only $31$ seconds.

  \begin{table}[htdp]
\begin{center}
\begin{tabular}{|c|c|c|c|c|}
\hline 
Approach & \multicolumn{2}{c|}{Offline phase} & \multicolumn{2}{c|}{Online phase}  \\
\hline
 &Wall time  (s) & Number of processors &Wall time (s)  & Number of processors \\
 \hline 
 HDM & - & - & 9,152,000 & 32\\
 \hline 
 RSE &  28,000& 32 & 2  & 1\\
  \hline 
ROM database  & 28,000 & 32 & 31& 1\\ 
 \hline 

  \end{tabular} 
  \caption{Wall times associated with computing the complete predicted FSI surface reported in Figure~\ref{fig:FSI_DB} and~\ref{fig:FSI_predictions}}\label{tab:CPUtimes}
\end{center}
\end{table}

\begin{figure}[htbp] \centering
 \includegraphics[width=0.45\textwidth]{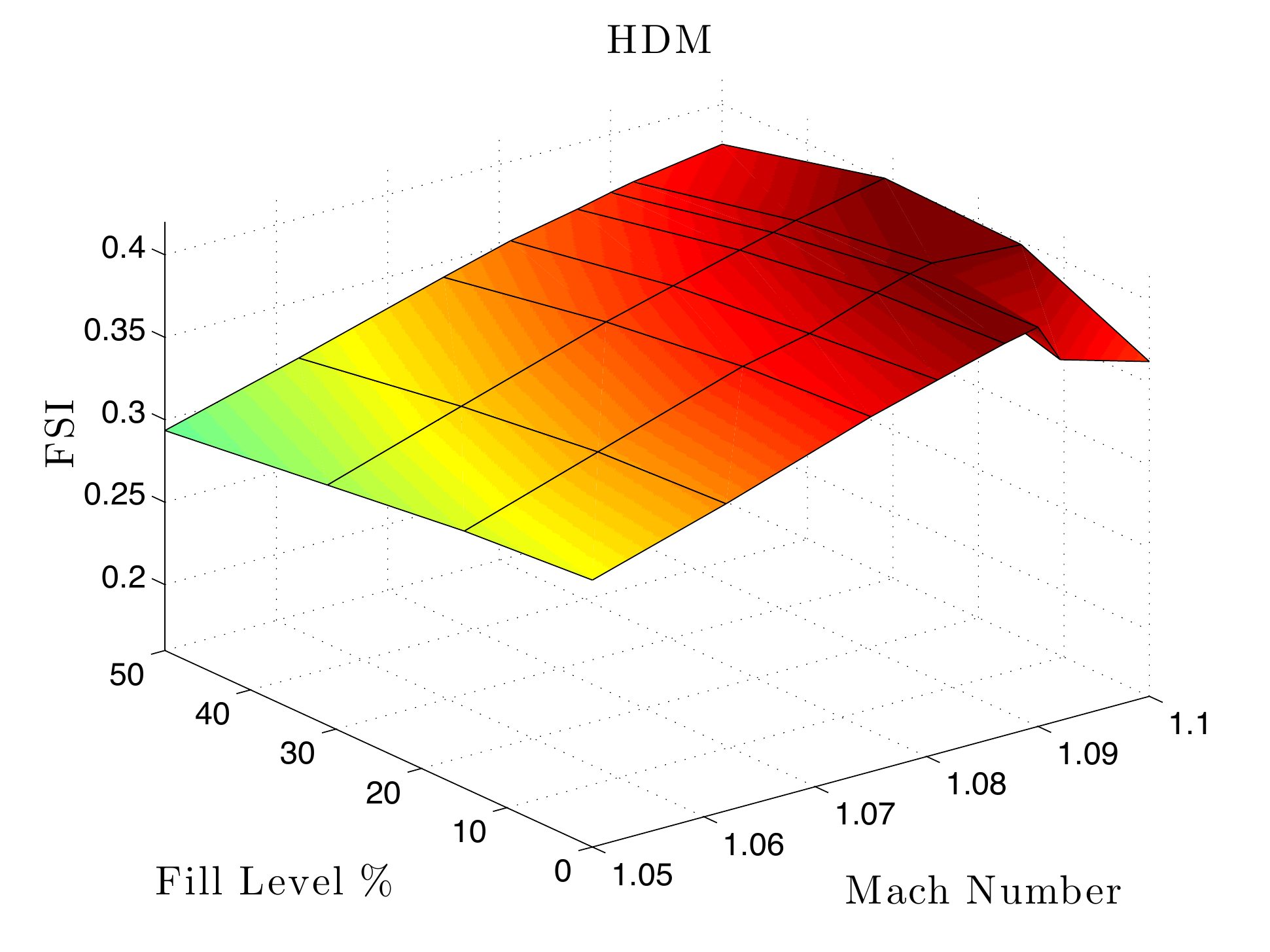}    \includegraphics[width=0.45\textwidth]{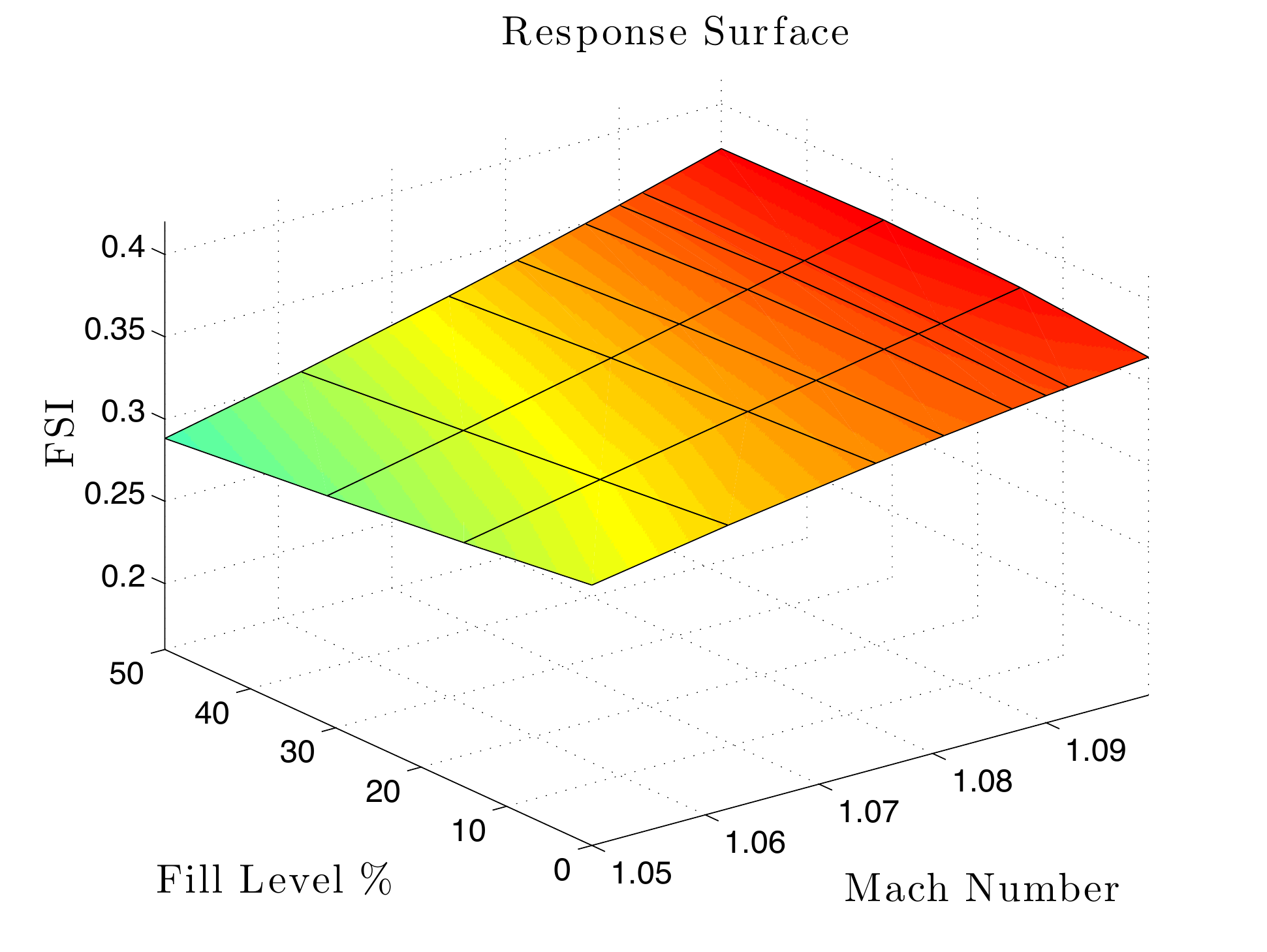}       
   \includegraphics[width=0.45\textwidth]{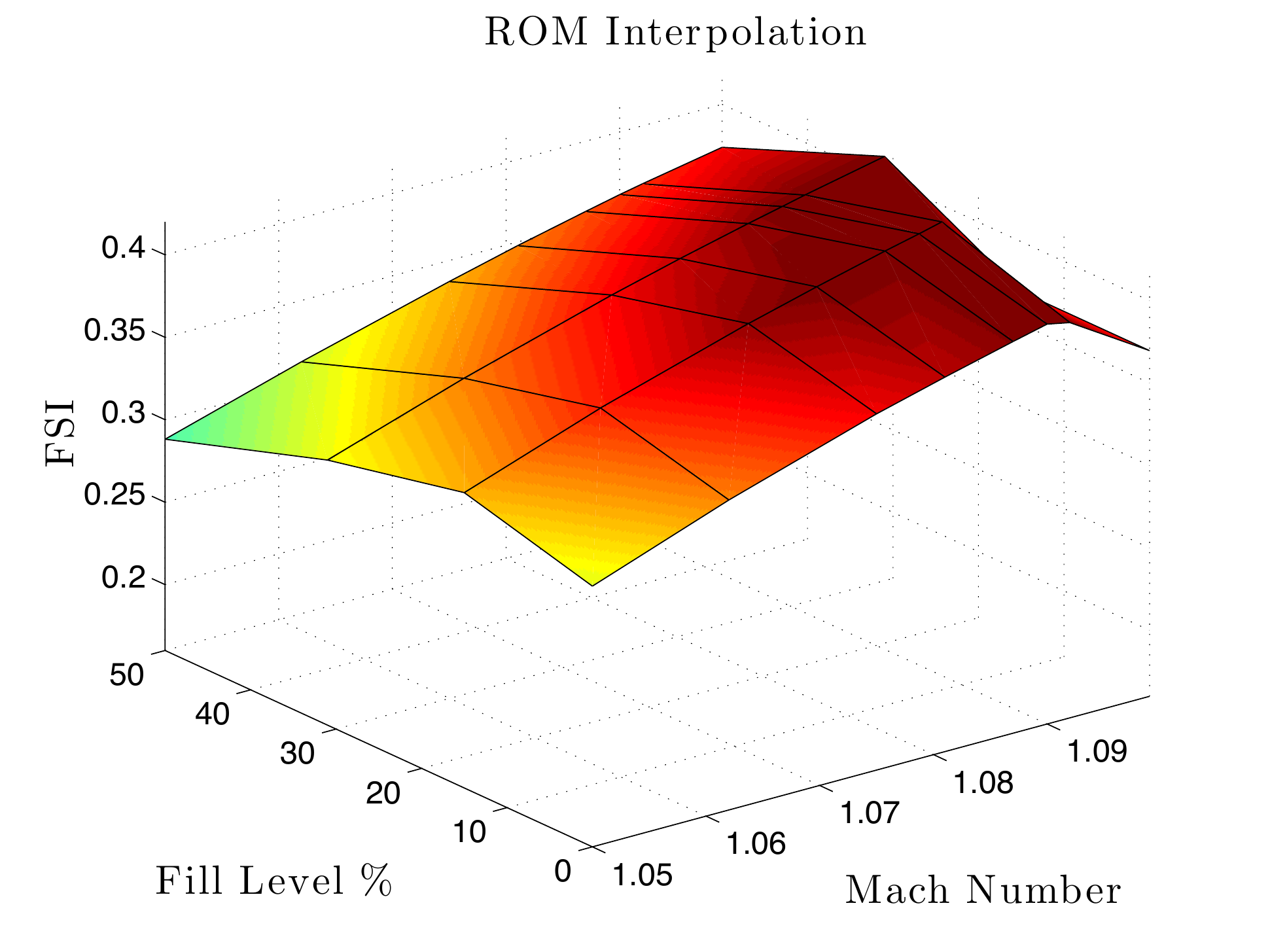}    
\caption{Comparison of high-dimensional model and predicted flutter speed indices between $M_\infty = 1.05$ and $M_\infty = 1.1$ for low fuel fill levels using (1) response surface estimation, (2) ROM interpolation}
\label{fig:FSI_predictions_zoom} 
\end{figure}

%
   
 \begin{figure}[htbp] \centering
 \includegraphics[width=0.4\textwidth]{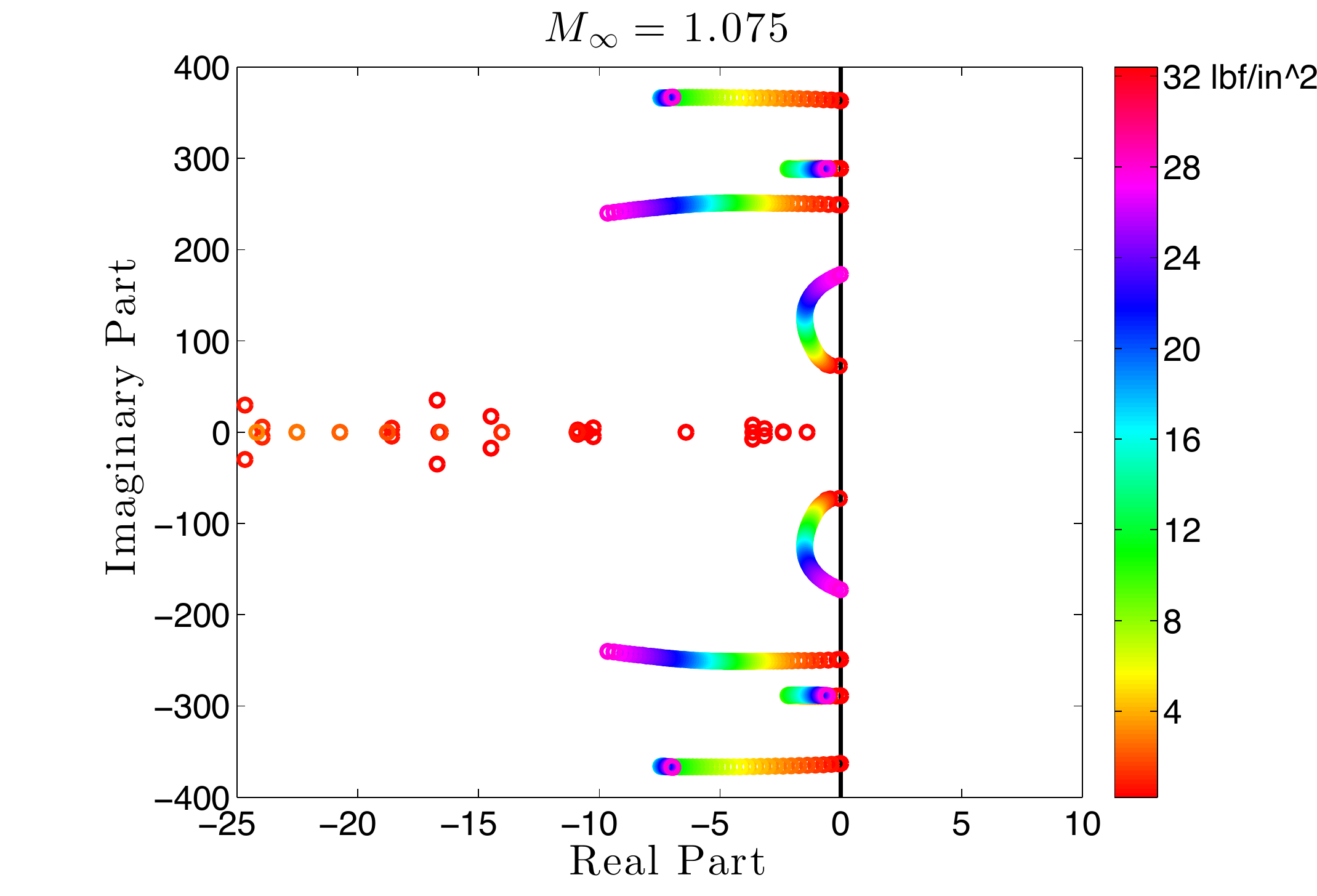}    \includegraphics[width=0.4\textwidth]{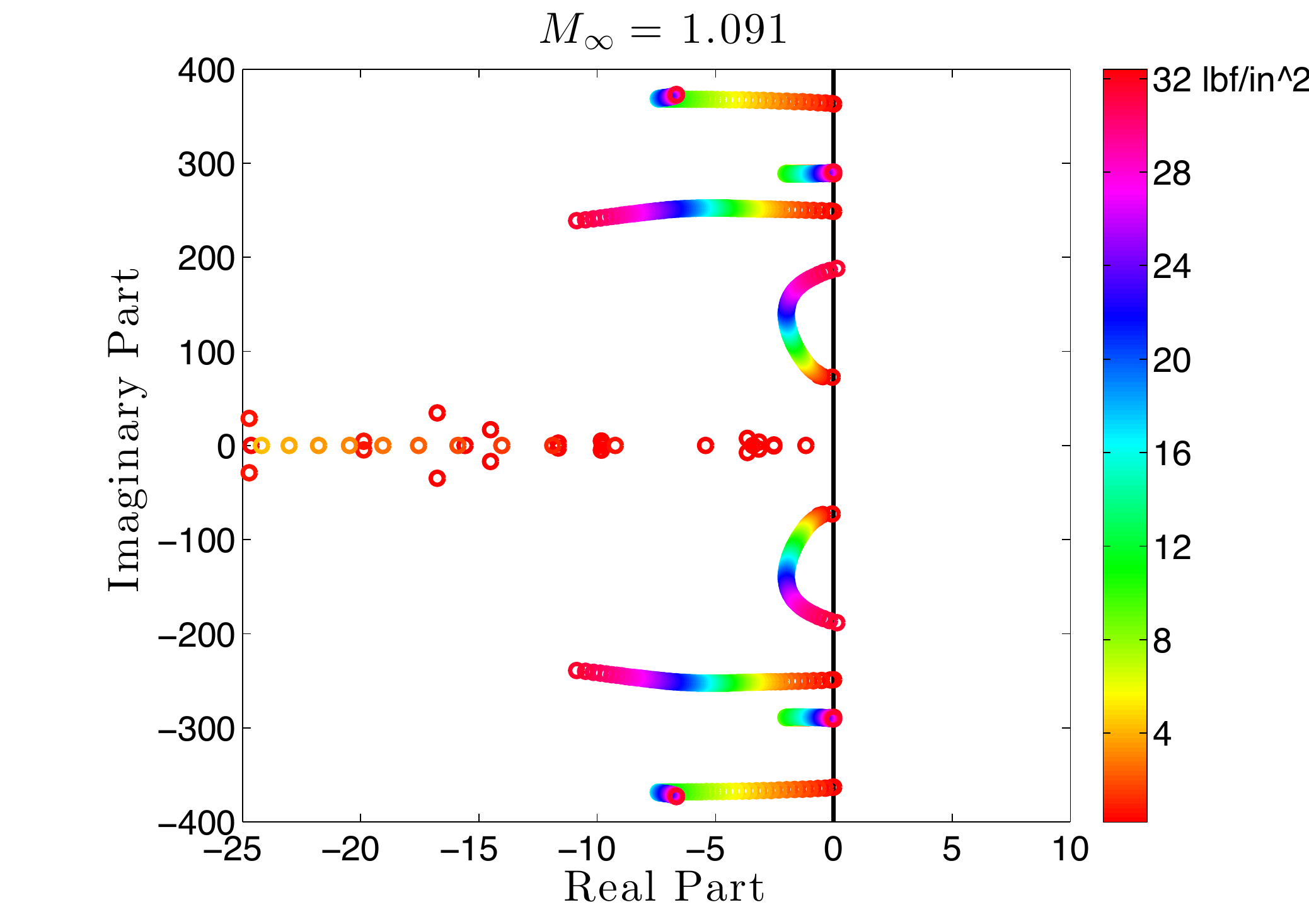}      
 \includegraphics[width=0.4\textwidth]{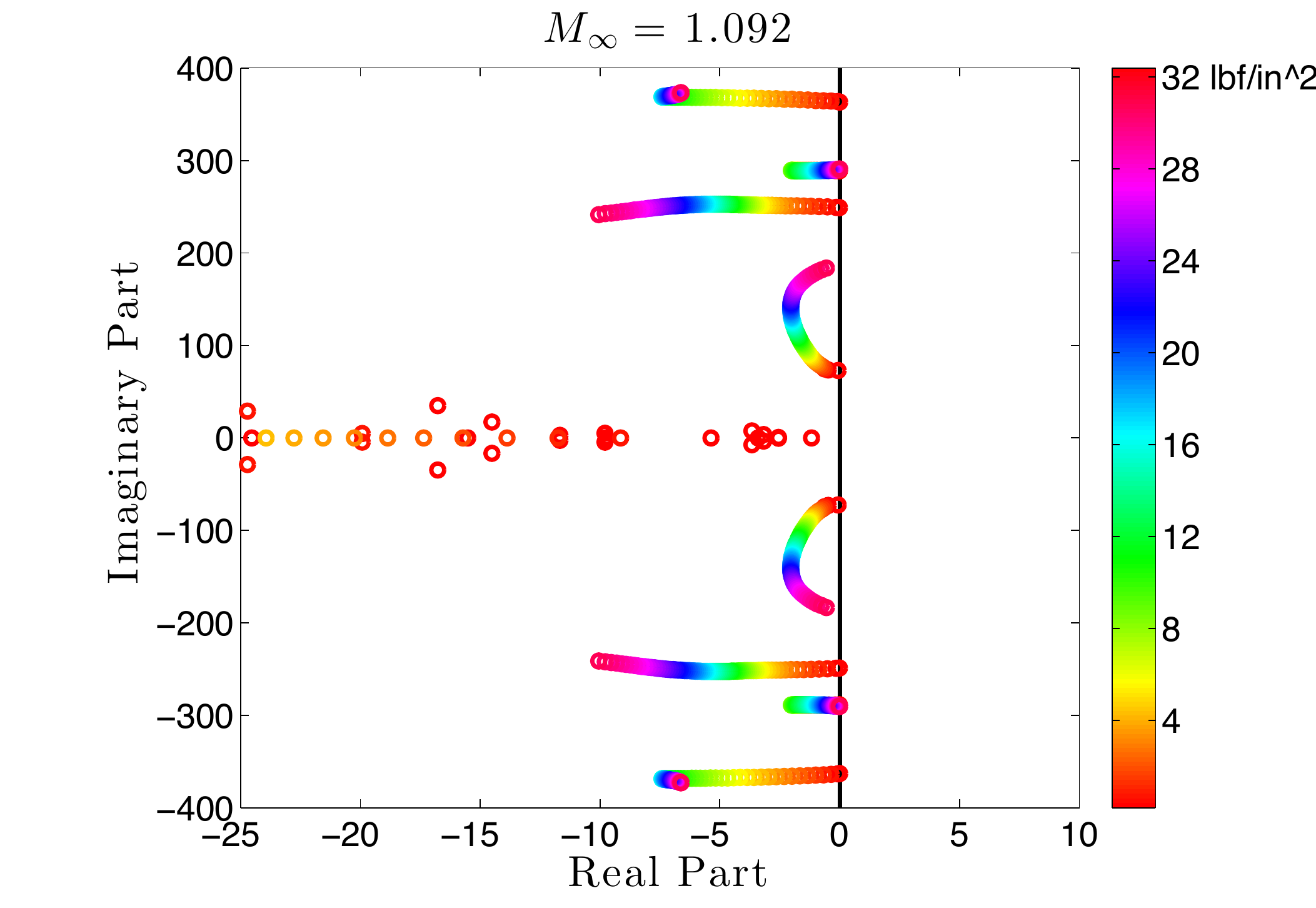}    \includegraphics[width=0.4\textwidth]{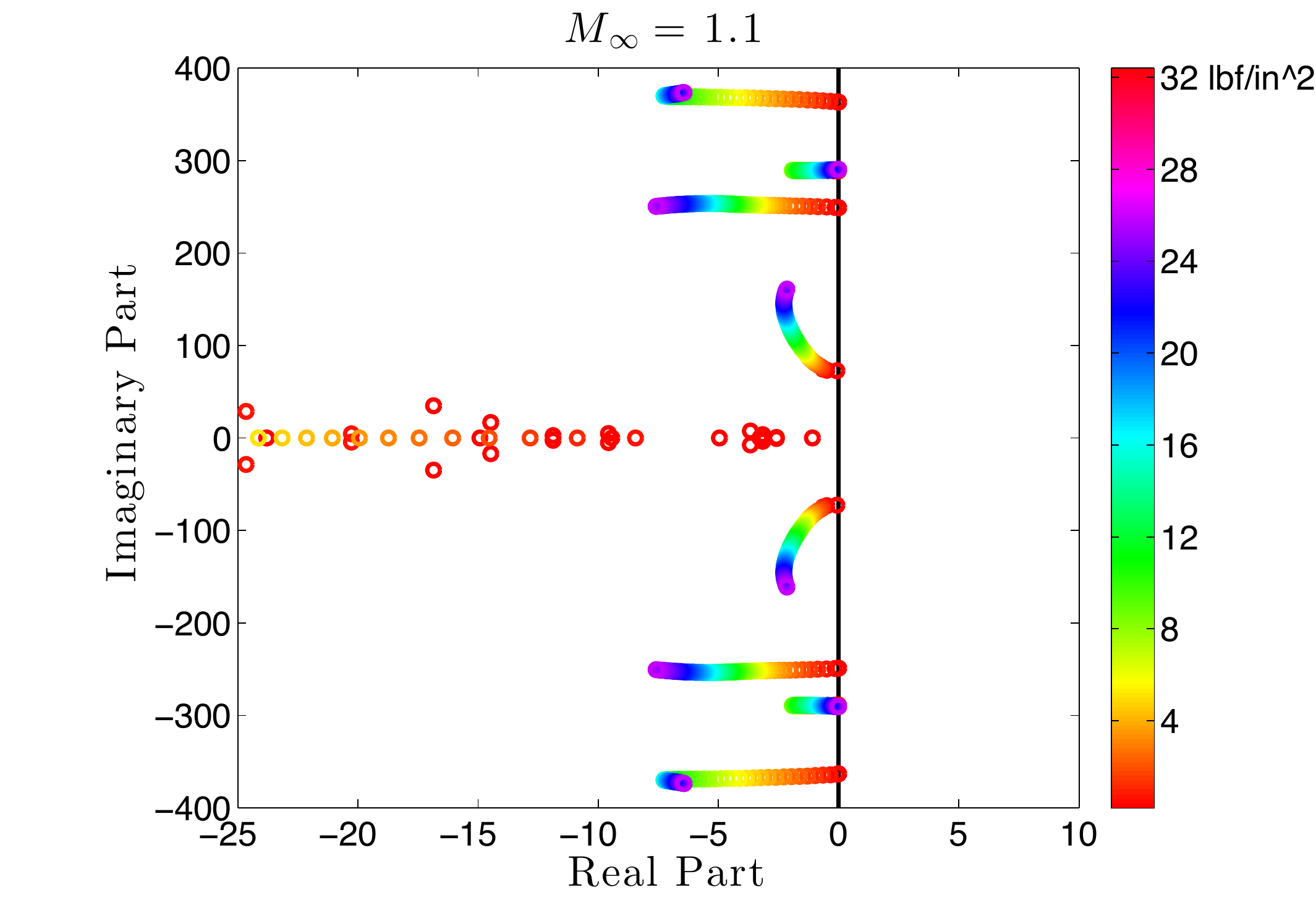}      
\caption{Aeroelastic matrix eigenvalues loci at $0\%$ fill level for various free-stream Mach numbers. }
\label{fig:bifurcation} 
\end{figure}

 \begin{figure}[htbp] \centering
 \includegraphics[width=0.4\textwidth]{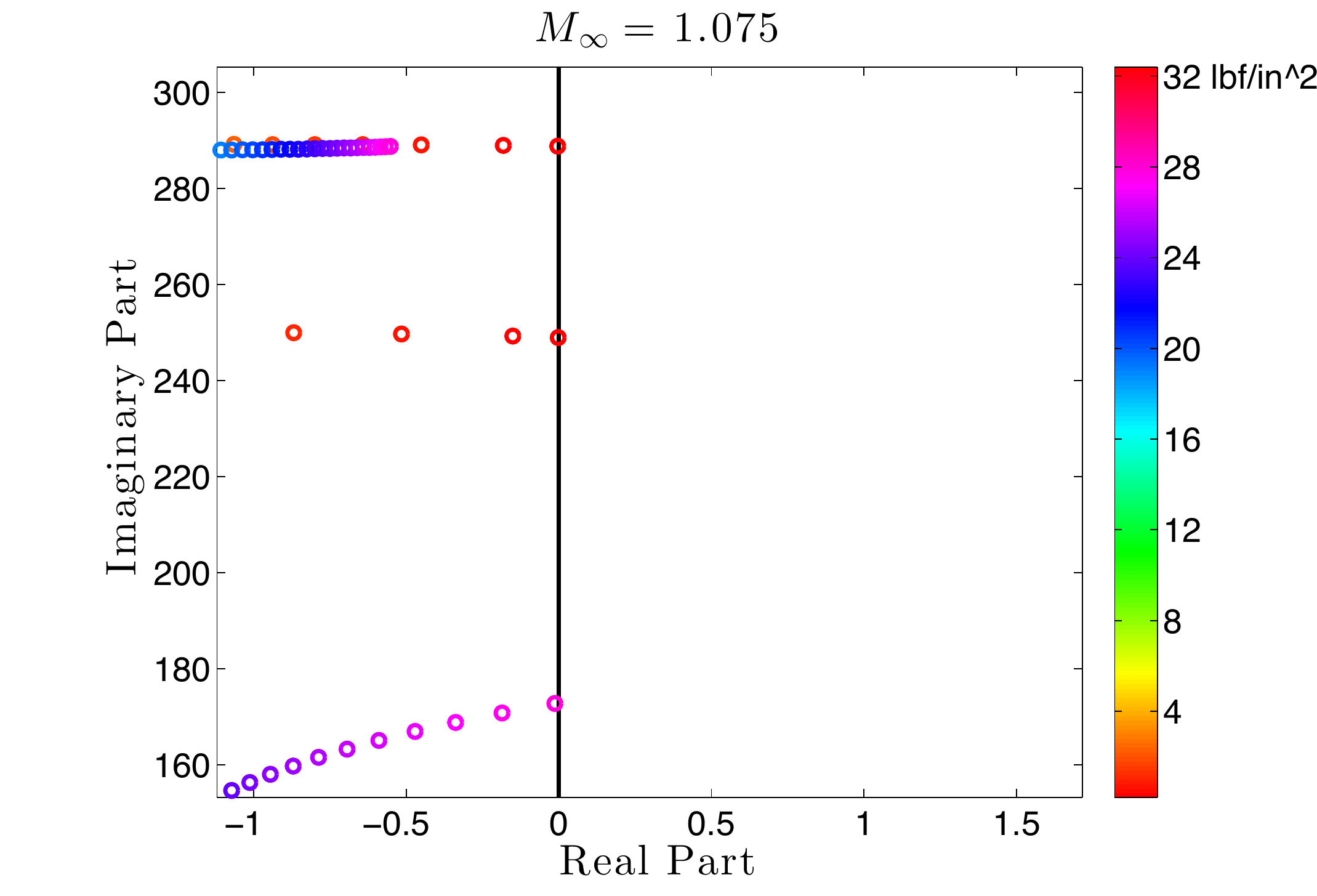}    \includegraphics[width=0.4\textwidth]{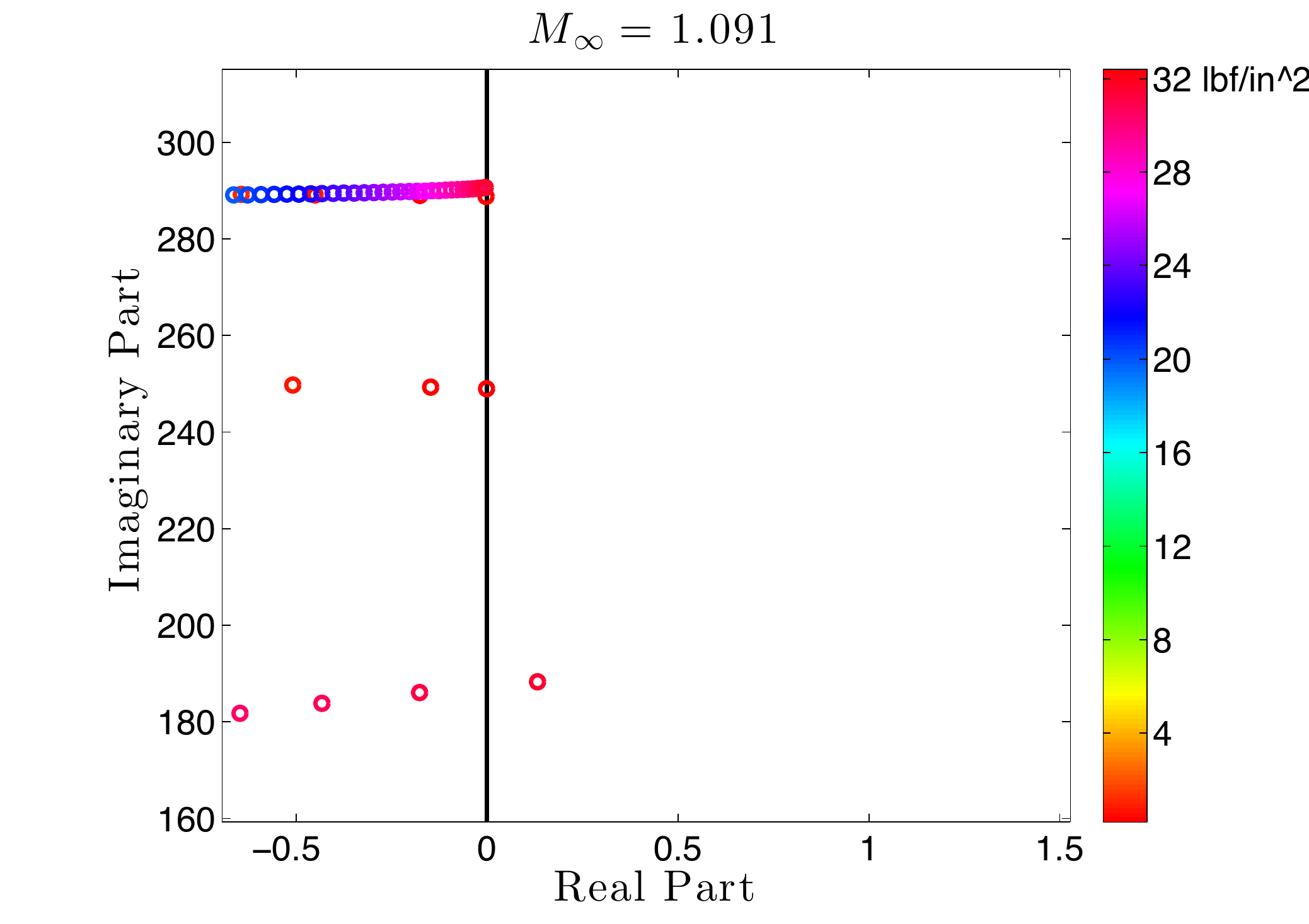}      
 \includegraphics[width=0.4\textwidth]{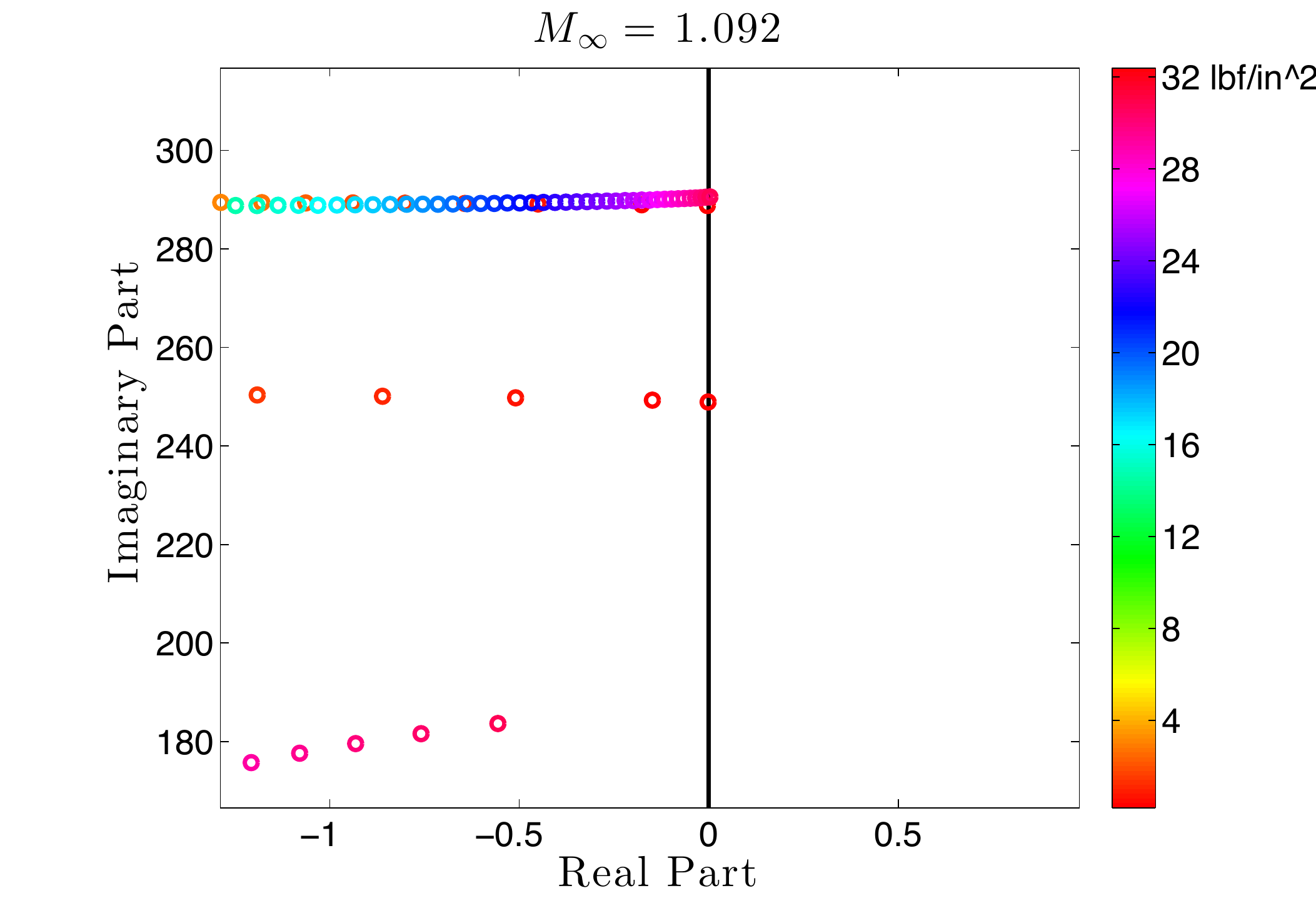}    \includegraphics[width=0.4\textwidth]{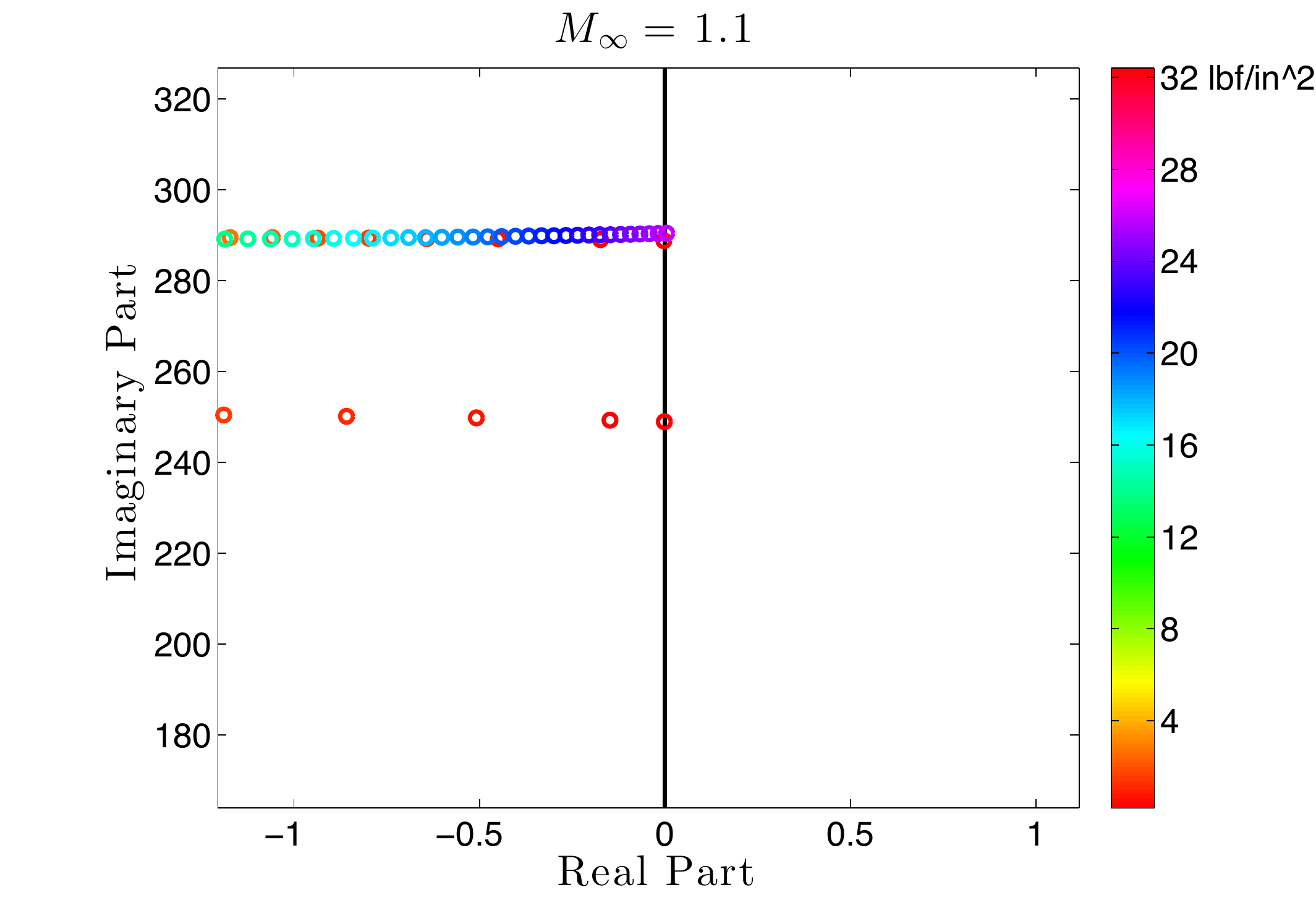}      
\caption{Aeroelastic matrix eigenvalues loci at $0\%$ fill level for various free-stream Mach numbers (zoom). }
\label{fig:bifurcation_zoom} 
\end{figure}

The flutter behavior at supersonic speed for low fuel levels is then studied more in detail by predicting the FSI at an additional fuel level for $15 \%$ tank fill. The corresponding FSI are reported in Figure~\ref{fig:FSI_predictions_zoom}. In order to understand the physical phenomenon at play, the eigenvalues of the interpolated aeroelastic ROM matrices are computed for increasing values of free-stream pressure until flutter is reached. By following the wet structural modes in the complex plane, one can determine which structural mode flutters, that is which one is the first to cross the imaginary axis. Results are reported in Figures~\ref{fig:bifurcation}  and~\ref{fig:bifurcation_zoom} for four different values of the free-stream mach number and an empty tank. One can observe that for $M_\infty=1.075$ and $M_\infty=1.091$, the first mode is the first to cross the imaginary axis while for $M_\infty = 1.092$ and $M_\infty = 1.1$, the third mode is the first one to flutter. These results clearly show that a bifurcation phenomenon is at play. Being able to perform such analyses demonstrates another  clear advantage of the proposed method over RSE. For that same fill level, the HDM predicts the same phenomenon, that is a peak of FSI in function of the free-stream Mach number $M_\infty$ between $M_\infty = 1.09$ and $M_\infty = 1.092$, which is in perfect agreement with the results arising from the interpolated ROM.

%
%
 
  \begin{figure}[htbp] \centering
 \includegraphics[width=0.45\textwidth]{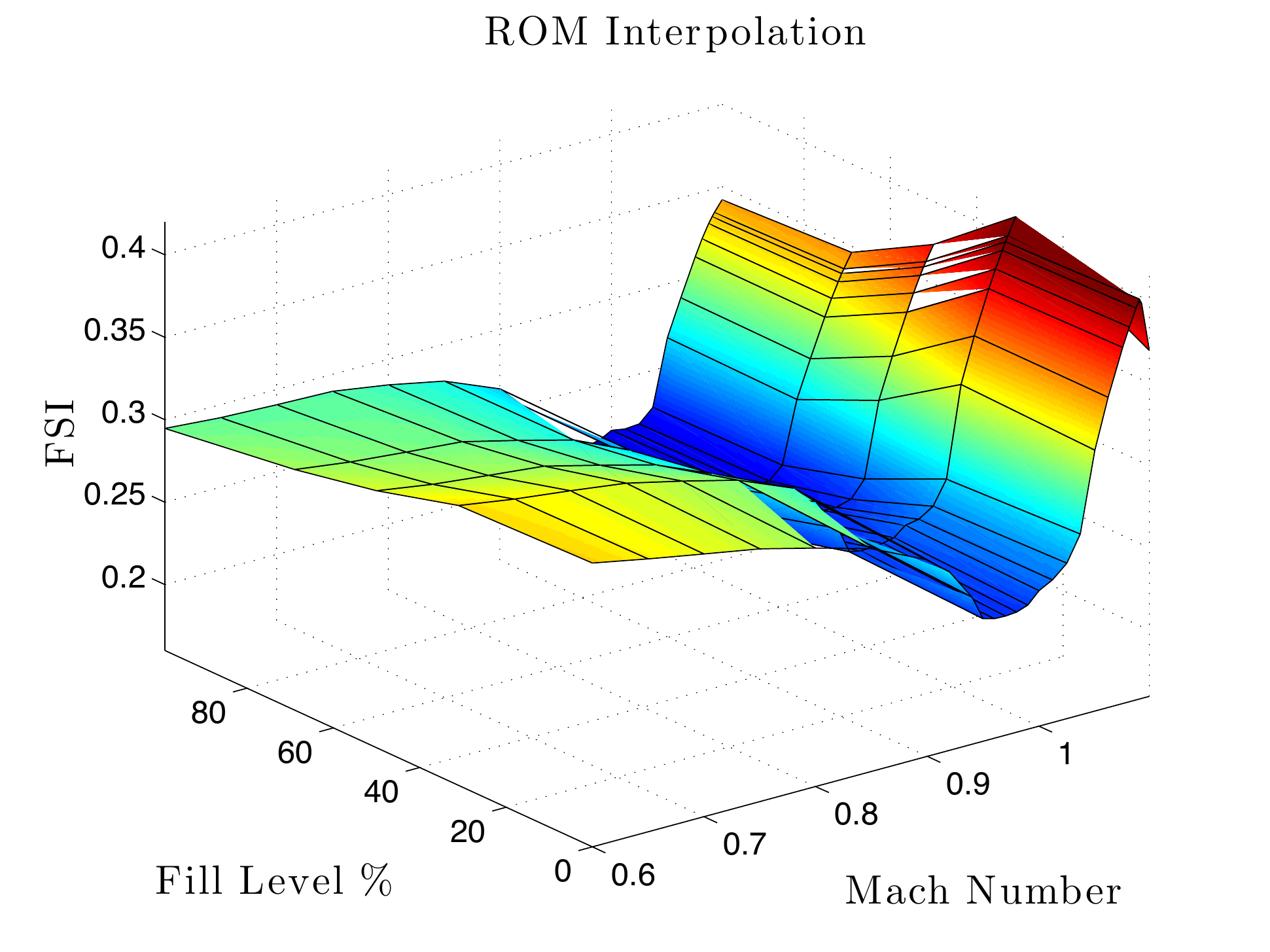}     
  \includegraphics[width=0.45\textwidth]{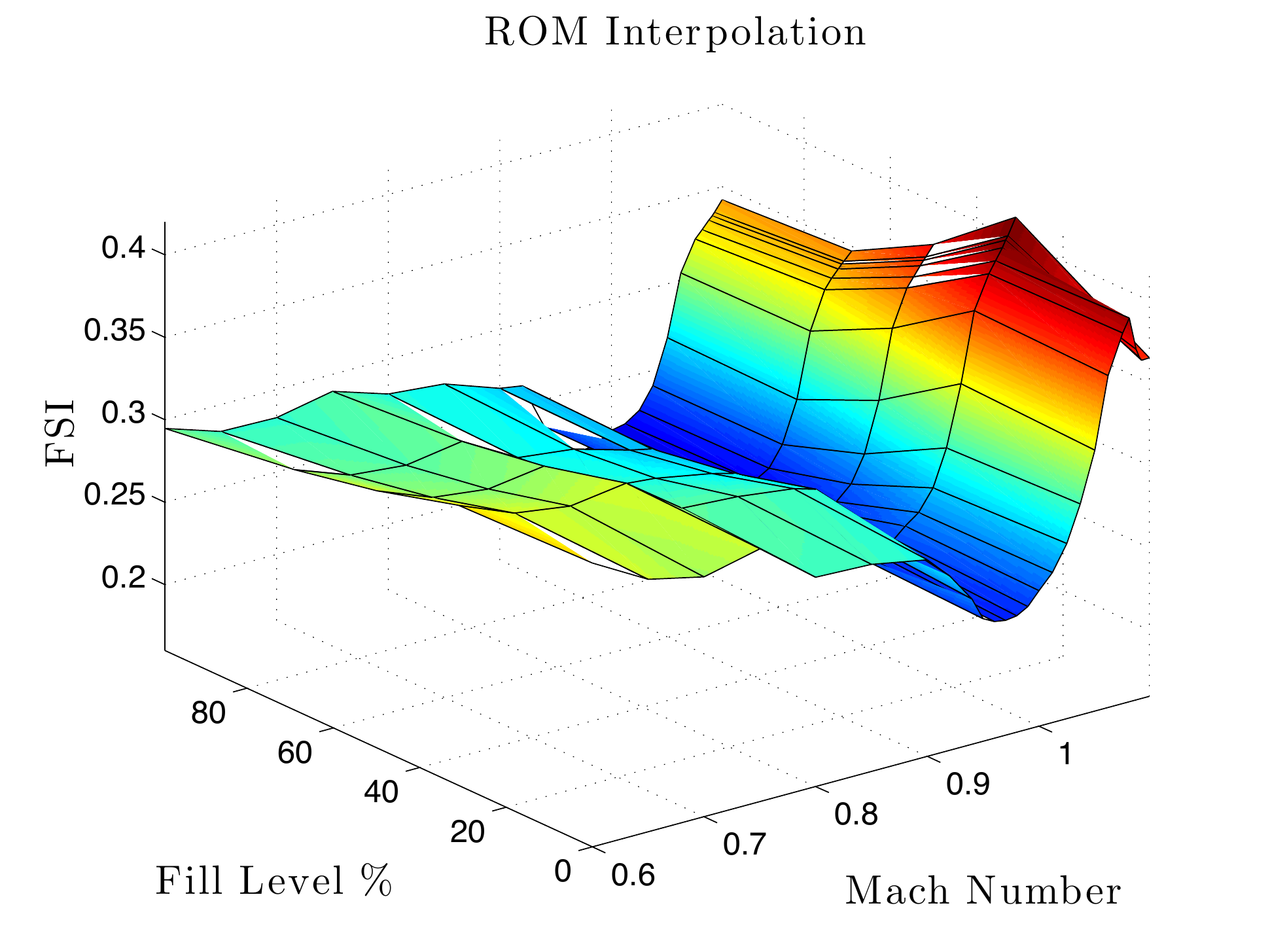}   
\caption{Predicted flutter speed indices using ROM interpolation without the manifold choice heuristic: interpolation of the reduced fluid operator on $\mathbb{R}^{k^{(f)}\times k^{(f)}}$ (left) and on $\text{GL}(k^{(f)})$ (right).}
\label{fig:FSI_predictions_noheuristic} 
\end{figure}

 Next, the effect of the manifold choice heuristic is studied by considering interpolation on $\mathbb{R}^{k^{(f)}\times k^{(f)}}$ and on $\text{GL}(k^{(f)})$, respectively, for the reduced fluid operator at every operating point. The corresponding results are reported in Figure~\ref{fig:FSI_predictions_noheuristic}  for interpolation using Choice 1 and 3, respectively. One can observe that the predicted FSI are more accurate when the heuristic is used, especially for subsonic and transonic flight conditions.

 \begin{figure}[htbp] \centering
 \includegraphics[width=0.45\textwidth]{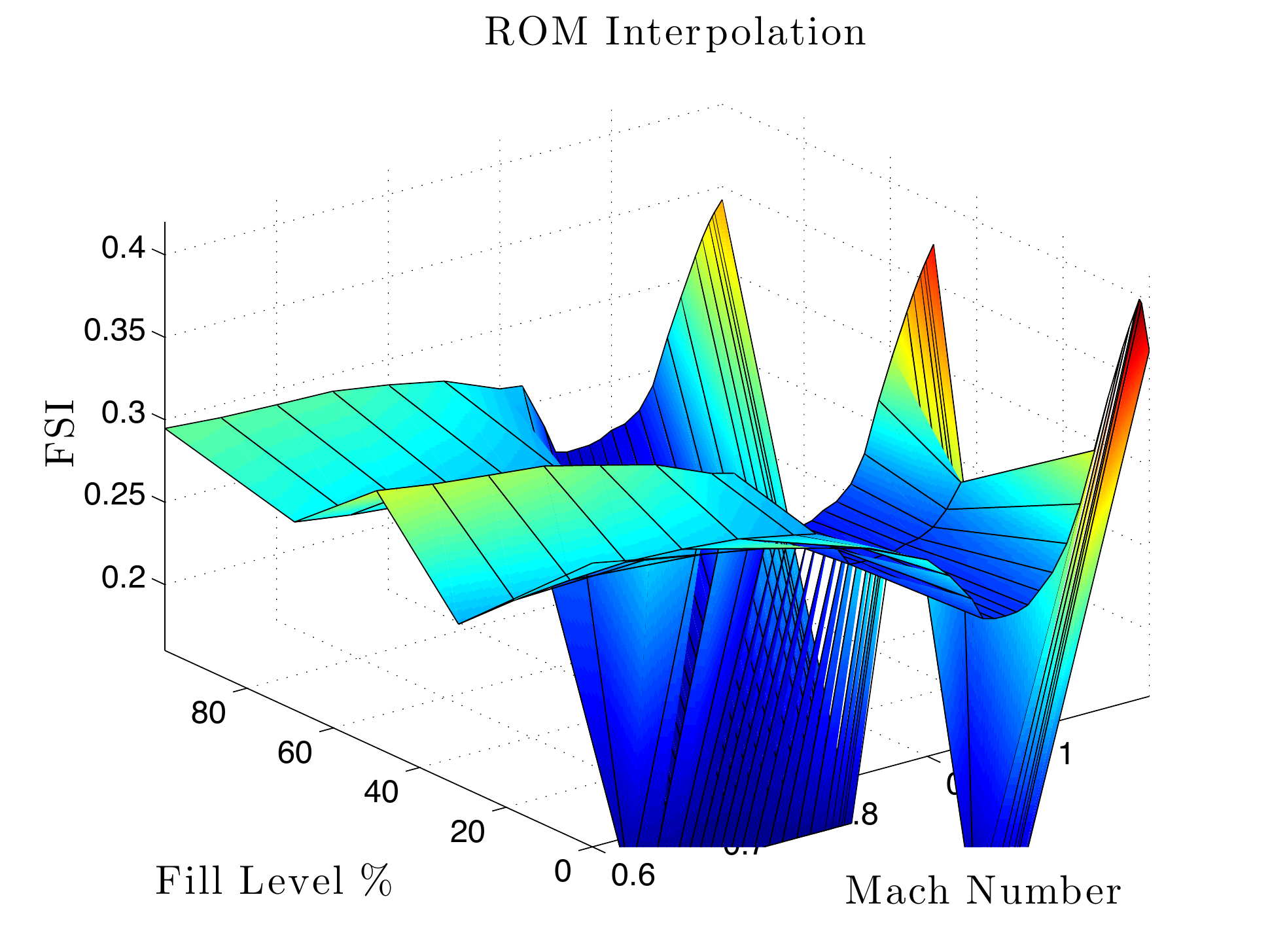}    \includegraphics[width=0.45\textwidth]{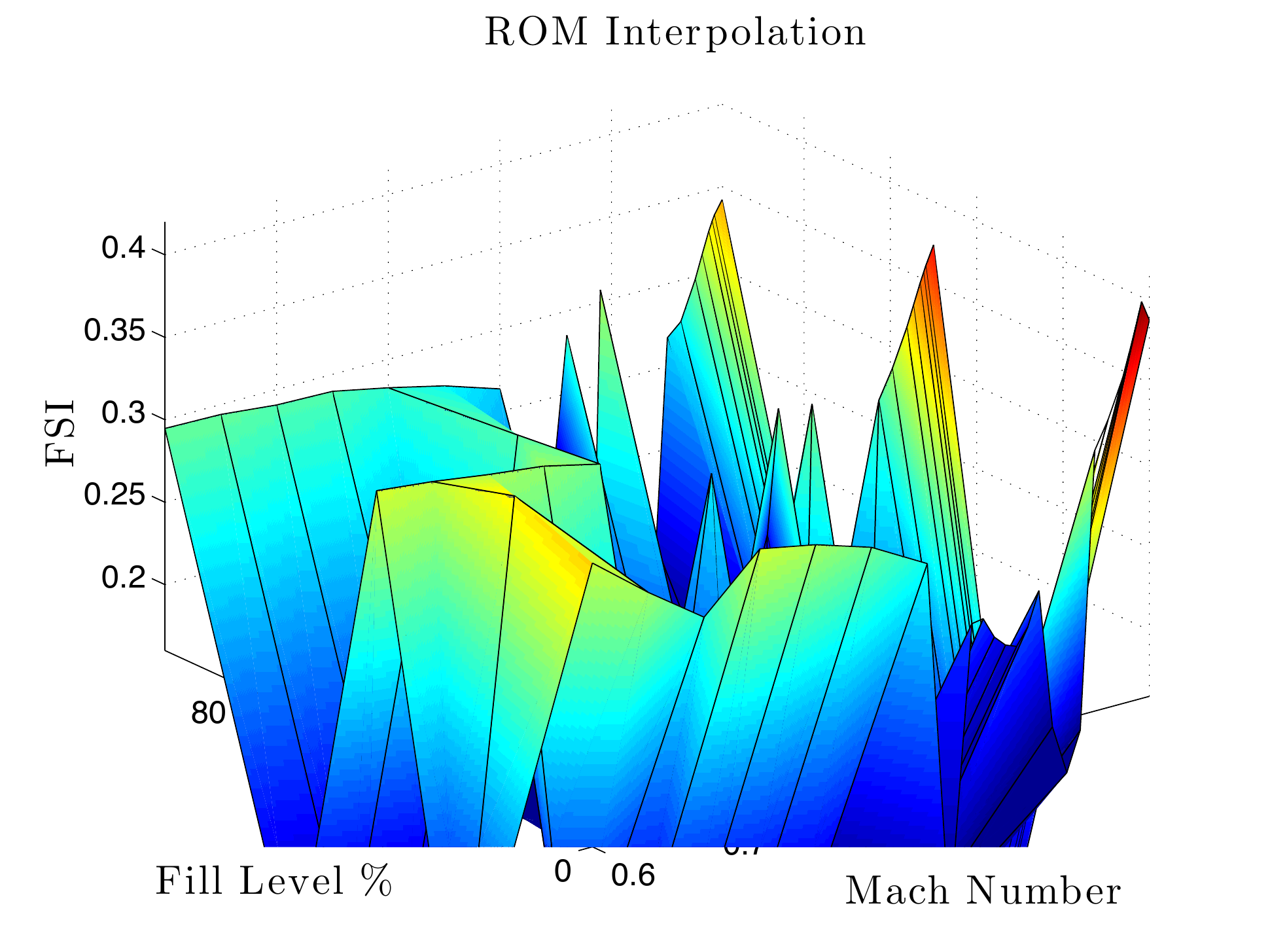}     
    \includegraphics[width=0.45\textwidth]{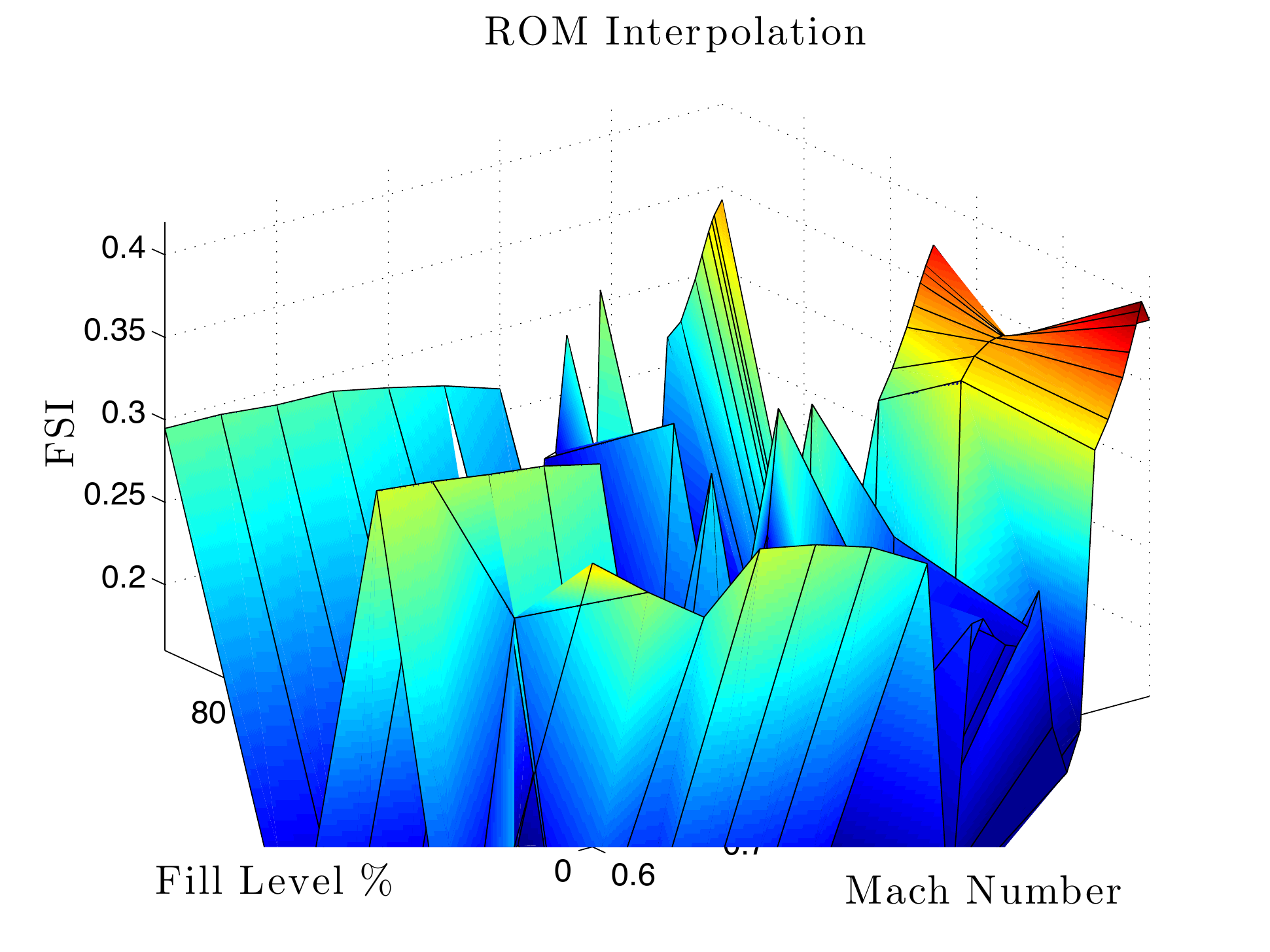}       
\caption{Predicted flutter speed indices using ROM interpolation  without consistency enforcement: for the structural operators only (top left), for the fluid operators only (top right), for all operators (bottom).}
\label{fig:FSI_predictions_norotation} 
\end{figure}

  Finally, the effect of consistency on the results accuracy is investigated. Inconsistent interpolation for both the structural and fluid subsystems and for only one of those two subsystems is performed and the corresponding FSI results reported in Figure~\ref{fig:FSI_predictions_norotation}. The reader can observe the crucial effect of consistency as none of the interpolation  of   inconsistent ROMs leads to accurate predicted FSI.

\begin{figure}[htbp] \centering
 \includegraphics[width=0.45\textwidth]{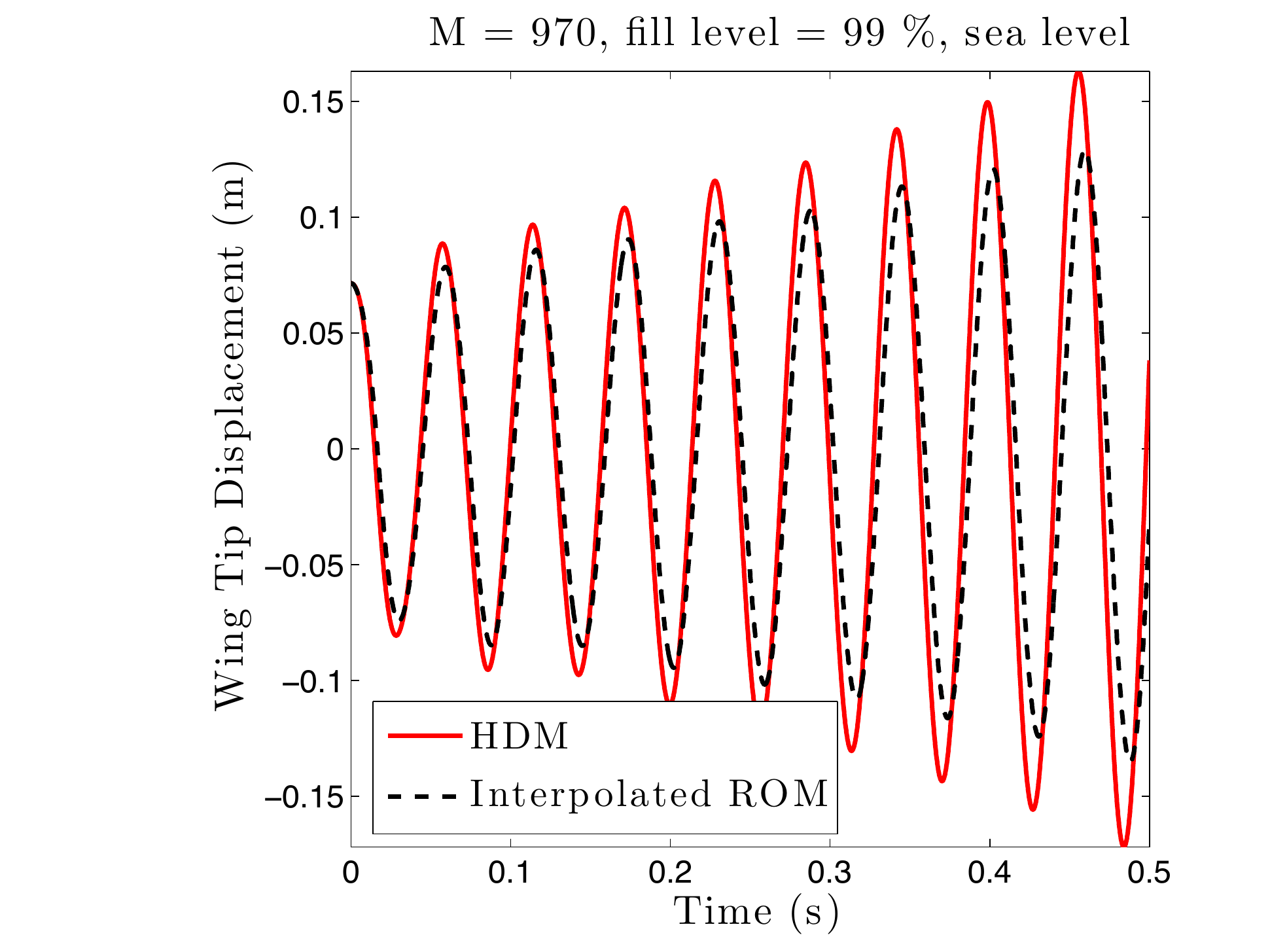}    \includegraphics[width=0.45\textwidth]{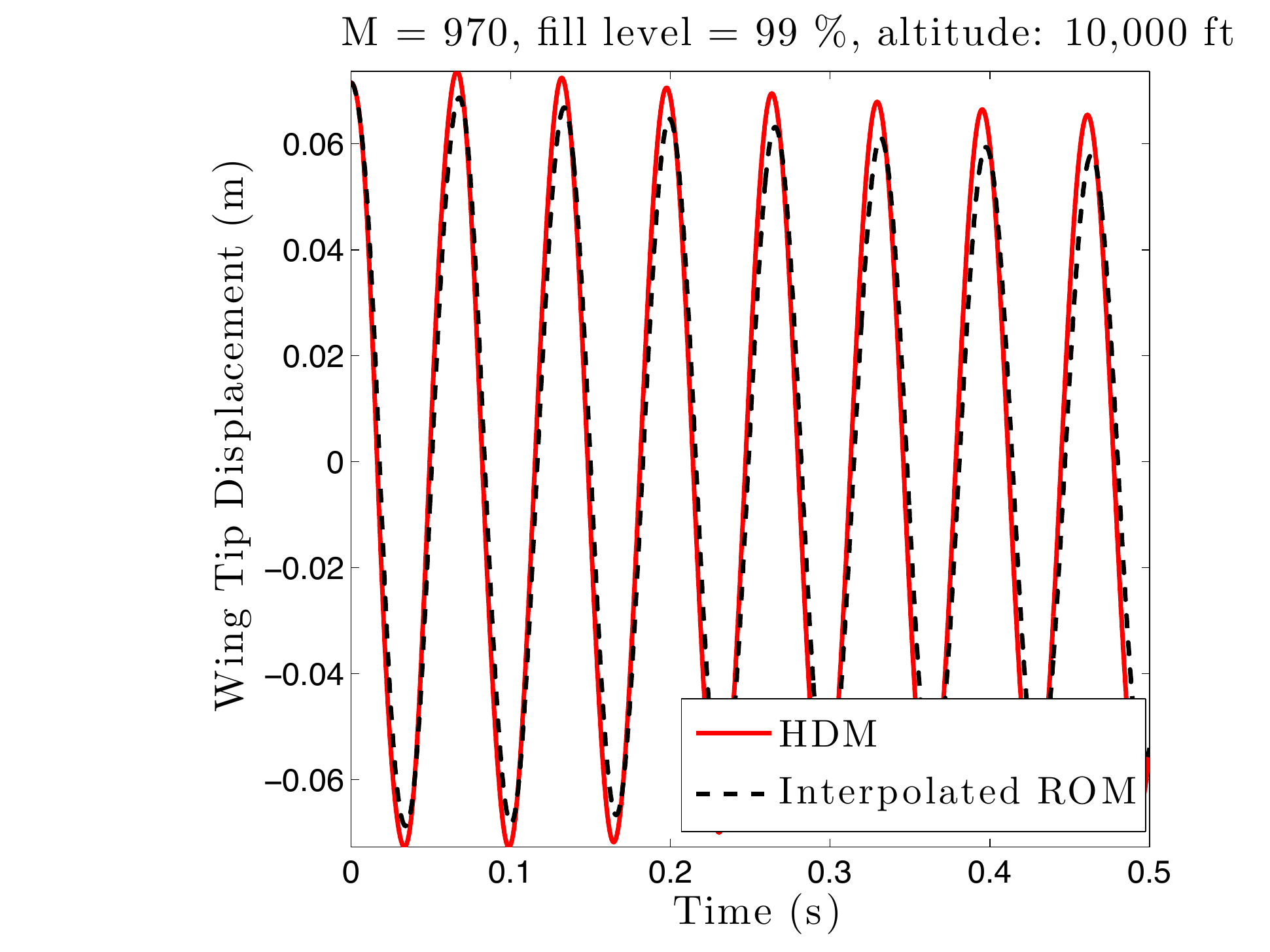}      
 \includegraphics[width=0.45\textwidth]{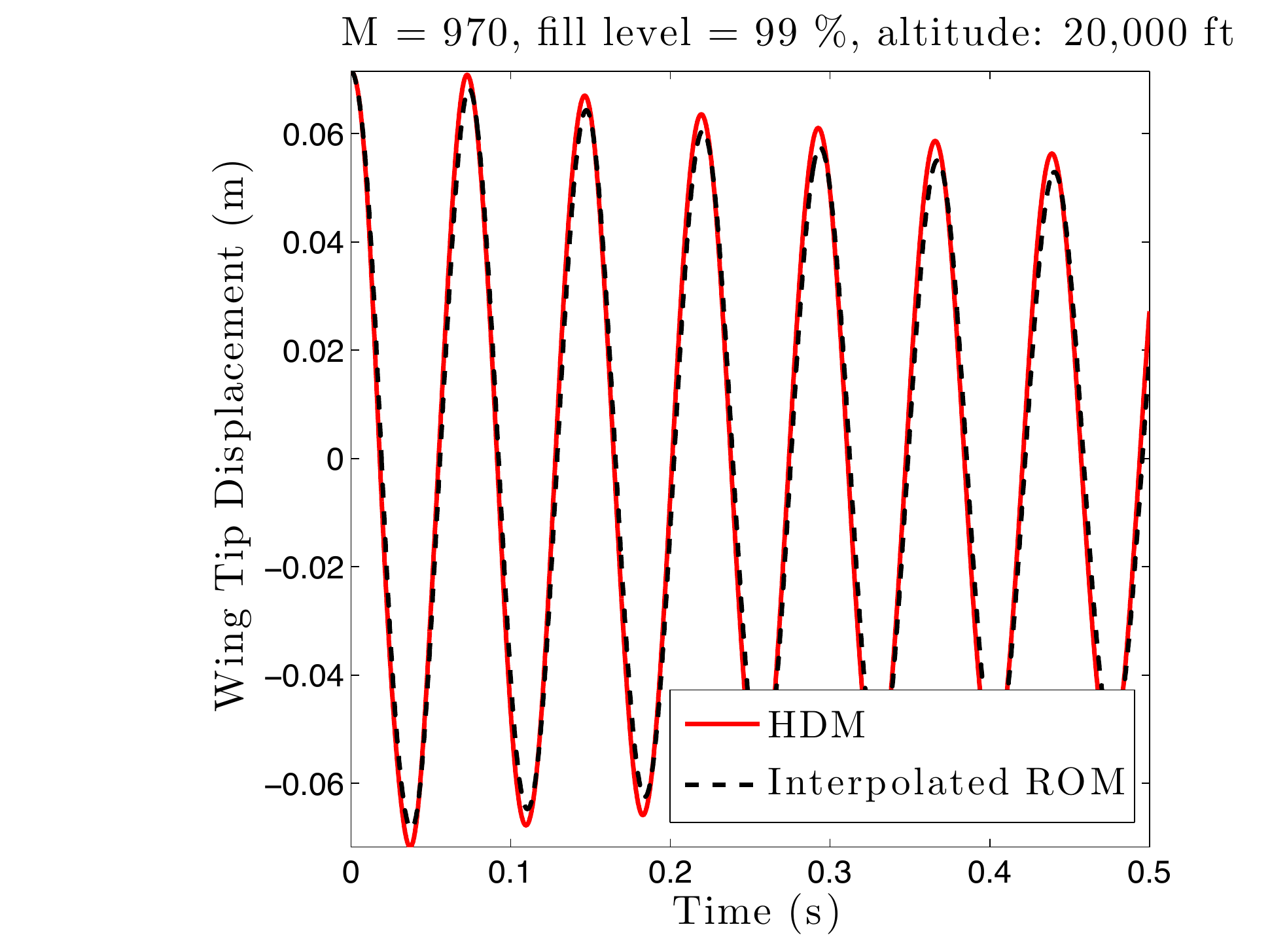}         
\caption{Wing tip displacement at $M_\infty = 0.97$ with full tank predicted using the high-dimensional model and an interpolated ROM at various altitudes. }
\label{fig:aeroelastic_responses} 
\end{figure}

The interpolated ROM can also be used to predict the displacement at a given location of the wing-tank system. Here wing tip displacements time histories are predicted using the interpolated ROM at the transonic tip, that is $M_\infty = 0.97$ for an empty tank. The corresponding results are reported in Figure~\ref{fig:aeroelastic_responses}  for three different cruise altitudes, and compared to predictions using the HDM. Good agreements can be observed.

\begin{figure}[htbp] \centering
\includegraphics[width=0.15\textheight]{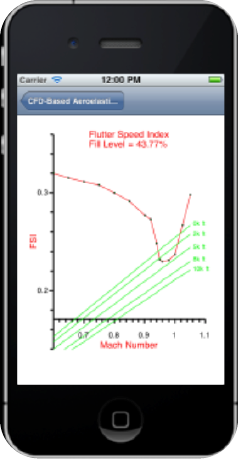}   
\caption{Screenshot of the iPhone application depicting the FSI for fill level $f=43.77\%$.}
\label{fig:iphone_app} \end{figure}

Finally, in order to demonstrate the capability of the proposed method to operate on mobile devices, an iPhone application is implemented for the aeroelastic system of interest. A screenshot of the application is displayed in Figure~\ref{fig:iphone_app}. The application can operate in the following two modes, based on the database of $N_{\text{DB}}=21$ points considered in this section: (1) In the first mode, for a given value of the fill level $f$, the FSI is compiled for $M_\infty\in[0.6,1.1]$. This is the mode depicted in Figure~\ref{fig:iphone_app}. (2) In the second mode, for a given combination $(M_\infty,f)$, the smallest aeroelastic damping ratio is computed for the altitude range $h\in[0,40000]$ ft.

\section{Conclusions}\label{sec:conclu}
This work presents a framework for real-time predictions based on a database of linear reduced-order models. It is based on the offline pre-computation of reduced-order models and their online interpolation at unsampled values of the parameters. A pre-processing step is first established to enforce the consistency of the set of generalized coordinates each reduced operator is defined by. The present paper presents such a step for both systems defined on common and arbitrary underlying meshes. The operators are then interpolated on the tangent space to a matrix manifold to enforce properties associated with each operator. The framework is then applied to two challenging multi-physics applications, demonstrated its capability  to lead to real-time and accurate predictions.
\section{Acknowledgements}
The authors would like to thank Mark Potts for the implementation of the iPhone application. The authors acknowledge partial support by the Army Research Laboratory through the Army High
Performance Computing Research Center under Cooperative Agreement W911NF-07-2-0027, and
partial support by the Office of Naval Research under grants no. N00014-11-1-0707 and N00014-14-1-0233.  This document
does not necessarily reflect the position of these institutions, and no official endorsement should be inferred.
\section{Appendix 1: proof of Theorem 1}
The objective function can be written as
\begin{equation}
\mathcal{J}_G(\Qbold) = \epsilon  \langle \Qbold^T \Ebold_{ri} \Qbold ,\Ebold_{r}^0\rangle + \alpha  \langle \Qbold^T \Abold_{ri} \Qbold ,\Abold_{r}^0\rangle + \langle \Fbold, \Qbold\rangle,
\end{equation}
where $\Fbold = \beta\Bbold_{ri}\left(\Bbold_{r}^0\right)^T +\gamma\Cbold_{ri}^T\Cbold_{r}^0$.

The Lagrangian of the optimization problem is then
\begin{gather}
\begin{split}
\mathcal{L}(\Qbold,\Sbold) &= \mathcal{J}_G(\Qbold) + \left\langle\frac{1}{2} \Sbold, \Ibold_k -\Qbold^T\Qbold\right\rangle\\
&= \epsilon  \langle \Qbold^T \Ebold_{ri} \Qbold ,\Ebold_{r}^0\rangle + \alpha  \langle \Qbold^T \Abold_{ri} \Qbold ,\Abold_{r}^0\rangle + \langle \Fbold, \Qbold\rangle + \left\langle\frac{1}{2} \Sbold, \Ibold_k -\Qbold^T\Qbold\right\rangle,
\end{split}
\end{gather}
where $\frac{1}{2}\Sbold\in\mathbb{R}^{k\times k}$ is a symmetric matrix of Lagrangian multipliers. Using the following identities~\cite{fraikin08},
\begin{eqnarray}
\nabla_\Qbold \langle \Mbold, \Qbold\rangle &=& \Mbold\\
\nabla_\Qbold \langle \Mbold, \Qbold^T\Qbold\rangle&=& \Qbold(\Mbold+\Mbold^T)\\
\nabla_\Qbold \langle \Qbold^T\Mbold\Qbold,\Nbold\rangle &=& \Mbold\Qbold\Nbold^T + \Mbold^T \Qbold \Nbold,
\end{eqnarray}
the gradient of the Lagrangian with respect to $\Qbold$ is obtained as
\begin{equation}
\nabla_\Qbold \mathcal{L}(\Qbold,\Sbold) = \epsilon\left(\Ebold_{ri}\Qbold\left(\Ebold_{r}^0\right)^T + \Ebold_{ri}^T \Qbold \Ebold_{r}^0 \right)+ \alpha\left(\Abold_{ri}\Qbold\left(\Abold_{r}^0\right)^T + \Abold_{ri}^T \Qbold \Abold_{r}^0 \right)+ \Fbold - \Qbold\Sbold
\end{equation}
which leads to the first-order optimality condition 
\begin{equation}
\Qbold\Sbold = \epsilon\left(\Ebold_{ri}\Qbold\left(\Ebold_{r}^0\right)^T + \Ebold_{ri}^T \Qbold \Ebold_{r}^0 \right)+ \alpha\left(\Abold_{ri}\Qbold\left(\Abold_{r}^0\right)^T + \Abold_{ri}^T \Qbold \Abold_{r}^0 \right)+ \Fbold,
\end{equation}
together with the constraint $\Qbold^T\Qbold = \Ibold_k$ and the property that $\Sbold$ is symmetric.

\section{Appendix 2: proof of Theorem 2}

The goal of this section is to prove that the set of fixed points of the proposed recursive algorithm is equal to the set of critical points of the objective function $\mathcal{J}_G$. 

Let $\widehat\Qbold$ denote a  fixed point of the recursive method defined in Algorithm 2. Then $\widehat\Qbold$ satisfies $\widehat\Qbold = \widehat\Ubold\widehat\Vbold^T$ where
\begin{equation}
\widehat\Ubold \widehat\Sigmabold \widehat\Vbold^T =\epsilon\left(\Ebold_{ri}\widehat\Qbold\left(\Ebold_{r}^0\right)^T + \Ebold_{ri}^T \widehat\Qbold \Ebold_{r}^0 \right)+ \alpha\left(\Abold_{ri}\widehat\Qbold\left(\Abold_{r}^0\right)^T + \Abold_{ri}^T \widehat\Qbold \Abold_{r}^0 \right) +s\widehat\Qbold \Fbold
\end{equation}
is a singular value decomposition. Since $\widehat\Vbold$ is an orthogonal matrix,
\begin{equation}
\widehat\Ubold \widehat\Vbold^T \widehat\Vbold\widehat\Sigmabold \widehat\Vbold^T = \epsilon\left(\Ebold_{ri}\widehat\Qbold\left(\Ebold_{r}^0\right)^T + \Ebold_{ri}^T \widehat\Qbold \Ebold_{r}^0 \right)+ \alpha\left(\Abold_{ri}\widehat\Qbold\left(\Abold_{r}^0\right)^T + \Abold_{ri}^T \widehat\Qbold \Abold_{r}^0 \right)  + s\widehat\Qbold + \Fbold,
\end{equation}
that is
\begin{equation}
\widehat\Qbold \Sbold = \epsilon\left(\Ebold_{ri}\widehat\Qbold\left(\Ebold_{r}^0\right)^T + \Ebold_{ri}^T \widehat\Qbold \Ebold_{r}^0 \right)+ \alpha\left(\Abold_{ri}\widehat\Qbold\left(\Abold_{r}^0\right)^T + \Abold_{ri}^T \widehat\Qbold \Abold_{r}^0 \right)  + \Fbold,
\end{equation}
where $\Sbold = \widehat\Vbold\widehat\Sigmabold \widehat\Vbold^T-s\Ibold_k$ is a symmetric matrix. $\widehat\Ubold$ and $\widehat\Vbold$ being orthogonal, $\widehat\Qbold$ is orthogonal as well and therefore meets the requirements of Theorem 2. The set of fixed points of Algorithm 2 is included in the set of critical point of $\mathcal{J}_G$.

Conversely, let $\Qbold^\star$ be a critical point of $\mathcal{J}_G$. Then, there exists a symmetric matrix $\Sbold$ such that Eq.~(\ref{eq:criticalPts}) holds with $\Qbold^\star$ orthogonal. Since, $\Sbold$ is real and symmetric, it is diagonalizable as
\begin{equation}
\Sbold = \Ubold \Lambdabold \Ubold^T,
\end{equation}
the eigenvalues in $\Lambdabold$ being real and ordered decreasingly and $\Ubold$ an orthogonal matrix. Then,
\begin{equation}\label{eq:Qstar}
\Qbold^\star \Ubold \Lambdabold \Ubold^T +s\Qbold^\star  = \epsilon\left(\Ebold_{ri}\Qbold ^\star\left(\Ebold_{r}^0\right)^T + \Ebold_{ri}^T\Qbold ^\star\Ebold_{r}^0\right) + \alpha\left(\Abold_{ri}\Qbold ^\star\left(\Abold_{r}^0\right)^T + \Abold_{ri}^T\Qbold ^\star\Abold_{r}^0\right) + \Fbold +s\Qbold^\star 
\end{equation}
which can also be written as
\begin{equation}
(\Qbold^\star \Ubold)( \Lambdabold +s\Ibold)\Ubold^T  = \epsilon\left(\Ebold_{ri}\Qbold ^\star\left(\Ebold_{r}^0\right)^T + \Ebold_{ri}^T\Qbold ^\star\Ebold_{r}^0\right) + \alpha\left(\Abold_{ri}\Qbold ^\star\left(\Abold_{r}^0\right)^T + \Abold_{ri}^T\Qbold ^\star\Abold_{r}^0\right)  +s\Qbold^\star + \Fbold.
\end{equation}
In order to conclude the proof, it remains to show that $(\Qbold^\star \Ubold)( \Lambdabold +s\Ibold)\Ubold^T$ is a singular value decomposition. $\Qbold^\star \Ubold$ and $\Ubold$  being orthogonal matrices, and $\Lambdabold +s\Ibold$ being a diagonal matrix, it is sufficient to show that $\Lambdabold +s\Ibold$ has all diagonal positive entries.

From Eq.~(\ref{eq:Qstar}),
\begin{equation}
 \|\Lambdabold\|_2 = \left\| \Qbold^\star \Ubold \Lambdabold \Ubold^T\right\|_2 =  \left\|\epsilon\left(\Ebold_{ri}\Qbold ^\star\left(\Ebold_{r}^0\right)^T + \Ebold_{ri}^T\Qbold ^\star\Ebold_{r}^0\right) + \alpha\left(\Abold_{ri}\Qbold ^\star\left(\Abold_{r}^0\right)^T + \Abold_{ri}^T\Qbold ^\star\Abold_{r}^0\right)  + \Fbold\right\|_2,
 \end{equation}
 and
\begin{gather}
\begin{split}
 \|  \Lambdabold \|_2 &\leq \epsilon\left\|\Ebold_{ri}\Qbold ^\star\left(\Ebold_{r}^0\right)^T + \Ebold_{ri}^T\Qbold ^\star\Ebold_{r}^0\right\|_2 + \alpha\left\|\Abold_{ri}\Qbold ^\star\left(\Abold_{r}^0\right)^T + \Abold_{ri}^T\Qbold ^\star\Abold_{r}^0\right\|_2 + \|\Fbold\|_2 \\
 &\leq  \epsilon\left(\|\Ebold_{ri}\|_2\|\Qbold^\star\|_2\left\|\left(\Ebold_{r}^0\right)^T\right\|_2 + \|\Ebold_{ri}^T\|_2\|\Qbold^\star\|_2\|\Ebold_{r}^0\|_2\right) \\&~~~~~~+  \alpha\left(\|\Abold_{ri}\|_2\|\Qbold^\star\|_2\left\|\left(\Abold_{r}^0\right)^T\right\|_2 + \|\Abold_{ri}^T\|_2\|\Qbold^\star\|_2\|\Abold_{r}^0\|_2\right) + \|\Fbold\|_2 \\
 &\leq 2\epsilon\|\Ebold_{ri}\|_2\|\Ebold_{r}^0\|_2  + 2\alpha\|\Abold_{ri}\|_2\|\Abold_{r}^0\|_2  + \|\Fbold\|_2\\
 &\leq s_{\text{min,G}}
 \end{split}
 \end{gather}
 by definition of $s_{\text{min,G}}$ in Eq.~(\ref{eq:smindef}). Denoting by $\lambda_i~i=1,\cdots,k$ the diagonal entries in $\Lambdabold$, this implies that
 \begin{equation}
 -s_{\text{min,G}} \leq \lambda_i \leq s_{min,G},~i=1,\cdots,k,
 \end{equation}
 and, since $s>s_{\text{min,G}}$, 
 \begin{equation}
 \lambda_i + s > 0,~i=1,\cdots,k.
 \end{equation}
The set of critical points of $\mathcal{J}_G$ is included in the set of fixed points of Algorithm 2 and the two sets are equal,  concluding the proof.

\section*{References}
\bibliographystyle{model1-num-names.bst}
\bibliography{paper}

\begin{thebibliography}{37}
\expandafter\ifx\csname natexlab\endcsname\relax\def\natexlab#1{#1}\fi
\providecommand{\bibinfo}[2]{#2}
\ifx\xfnm\relax \def\xfnm[#1]{\unskip,\space#1}\fi
\bibitem[{Moore(1981)}]{moore81}
\bibinfo{author}{B.~Moore},
\newblock \bibinfo{title}{{Principal component analysis in linear systems:
  Controllability, observability, and model reduction}},
\newblock \bibinfo{journal}{IEEE Transactions on Automatic Control}
  \bibinfo{volume}{26} (\bibinfo{year}{1981}) \bibinfo{pages}{17--32}.
\bibitem[{Sirovich(1987)}]{sirovich87}
\bibinfo{author}{L.~Sirovich},
\newblock \bibinfo{title}{{Turbulence and the dynamics of coherent structures.
  Part I: coherent structures}},
\newblock \bibinfo{journal}{Quarterly of applied mathematics}
  \bibinfo{volume}{45} (\bibinfo{year}{1987}) \bibinfo{pages}{561--571}.
\bibitem[{Ryckelynck(2005)}]{ryckelynck05}
\bibinfo{author}{D.~Ryckelynck},
\newblock \bibinfo{title}{{A priori hyperreduction method: an adaptive
  approach}},
\newblock \bibinfo{journal}{Journal of Computational Physics}
  \bibinfo{volume}{202} (\bibinfo{year}{2005}) \bibinfo{pages}{346--366}.
\bibitem[{Chaturantabut and Sorensen(2010)}]{chaturantabut10}
\bibinfo{author}{S.~Chaturantabut}, \bibinfo{author}{D.~Sorensen},
\newblock \bibinfo{title}{{Nonlinear model reduction via discrete empirical
  interpolation}},
\newblock \bibinfo{journal}{SIAM Journal on Scientific Computing}
  \bibinfo{volume}{32} (\bibinfo{year}{2010}) \bibinfo{pages}{2737--2764}.
\bibitem[{Carlberg et~al.(2011)Carlberg, Bou-Mosleh, and Farhat}]{carlberg11}
\bibinfo{author}{K.~Carlberg}, \bibinfo{author}{C.~Bou-Mosleh},
  \bibinfo{author}{C.~Farhat},
\newblock \bibinfo{title}{{Efficient non-linear model reduction via a
  least-squares Petrov--Galerkin projection and compressive tensor
  approximations}},
\newblock \bibinfo{journal}{International Journal for Numerical Methods in
  Engineering} \bibinfo{volume}{86} (\bibinfo{year}{2011})
  \bibinfo{pages}{155--181}.
\bibitem[{Amsallem et~al.(2012)Amsallem, Zahr, and
  Farhat}]{amsallem12:localROB}
\bibinfo{author}{D.~Amsallem}, \bibinfo{author}{M.~J. Zahr},
  \bibinfo{author}{C.~Farhat},
\newblock \bibinfo{title}{{Nonlinear model order reduction based on local
  reduced-order bases}},
\newblock \bibinfo{journal}{International Journal for Numerical Methods in
  Engineering} \bibinfo{volume}{92} (\bibinfo{year}{2012})
  \bibinfo{pages}{891--916}.
\bibitem[{Carlberg et~al.(2013)Carlberg, Farhat, Cortial, and
  Amsallem}]{carlberg13}
\bibinfo{author}{K.~Carlberg}, \bibinfo{author}{C.~Farhat},
  \bibinfo{author}{J.~Cortial}, \bibinfo{author}{D.~Amsallem},
\newblock \bibinfo{title}{{The GNAT method for nonlinear model reduction:
  effective implementation and application to computational fluid dynamics and
  turbulent flows}},
\newblock \bibinfo{journal}{Journal of Computational Physics}
  \bibinfo{volume}{242} (\bibinfo{year}{2013}) \bibinfo{pages}{623--647}.
\bibitem[{Farhat et~al.(2014)Farhat, Avery, Chapman, and Cortial}]{farhat14}
\bibinfo{author}{C.~Farhat}, \bibinfo{author}{P.~Avery},
  \bibinfo{author}{T.~Chapman}, \bibinfo{author}{J.~Cortial},
\newblock \bibinfo{title}{{Dimensional reduction of nonlinear finite element
  dynamic models with finite rotations and energy-based mesh sampling and
  weighting for computational efficiency}},
\newblock \bibinfo{journal}{International Journal for Numerical Methods in
  Engineering} \bibinfo{volume}{98} (\bibinfo{year}{2014})
  \bibinfo{pages}{625--662}.
\bibitem[{Veroy and Patera(2005)}]{veroy05}
\bibinfo{author}{K.~Veroy}, \bibinfo{author}{A.~T. Patera},
\newblock \bibinfo{title}{{Certified real-time solution of the parametrized
  steady incompressible Navier-Stokes equations: rigorous reduced-basis a
  posteriori error bounds}},
\newblock \bibinfo{journal}{International Journal for Numerical Methods in
  Fluids} \bibinfo{volume}{47} (\bibinfo{year}{2005})
  \bibinfo{pages}{773--788}.
\bibitem[{Amsallem and Farhat(2008)}]{amsallem08}
\bibinfo{author}{D.~Amsallem}, \bibinfo{author}{C.~Farhat},
\newblock \bibinfo{title}{{Interpolation method for adapting reduced-order
  models and application to aeroelasticity}},
\newblock \bibinfo{journal}{AIAA Journal} \bibinfo{volume}{46}
  (\bibinfo{year}{2008}) \bibinfo{pages}{1803--1813}.
\bibitem[{Amsallem et~al.(2010)Amsallem, Cortial, and Farhat}]{amsallem10}
\bibinfo{author}{D.~Amsallem}, \bibinfo{author}{J.~Cortial},
  \bibinfo{author}{C.~Farhat},
\newblock \bibinfo{title}{{Toward real-time computational-fluid-dynamics-based
  aeroelastic computations using a database of reduced-order information}},
\newblock \bibinfo{journal}{AIAA Journal} \bibinfo{volume}{48}
  (\bibinfo{year}{2010}) \bibinfo{pages}{2029--2037}.
\bibitem[{Negri et~al.(2015)Negri, Manzoni, and Amsallem}]{negri15}
\bibinfo{author}{F.~Negri}, \bibinfo{author}{A.~Manzoni},
  \bibinfo{author}{D.~Amsallem},
\newblock \bibinfo{title}{{Efficient model reduction of parametrized systems by
  matrix discrete empirical interpolation}},
\newblock \bibinfo{journal}{Submitted for publication}  (\bibinfo{year}{2015}).
\bibitem[{Paul-Dubois-Taine and Amsallem(2015)}]{pdt14}
\bibinfo{author}{A.~Paul-Dubois-Taine}, \bibinfo{author}{D.~Amsallem},
\newblock \bibinfo{title}{{An adaptive and efficient greedy procedure for the
  optimal training of parametric reduced-order models}},
\newblock \bibinfo{journal}{International Journal for Numerical Methods in
  Engineering} \bibinfo{volume}{102} (\bibinfo{year}{2015})
  \bibinfo{pages}{1262--1292}.
\bibitem[{Amsallem et~al.(2015{\natexlab{a}})Amsallem, Zahr, Choi, and
  Farhat}]{amsallem14:smo}
\bibinfo{author}{D.~Amsallem}, \bibinfo{author}{M.~J. Zahr},
  \bibinfo{author}{Y.~Choi}, \bibinfo{author}{C.~Farhat},
\newblock \bibinfo{title}{{Design Optimization Using Hyper-Reduced-Order
  Models}},
\newblock \bibinfo{journal}{Structural and Multidisciplinary Optimization}
  \bibinfo{volume}{51} (\bibinfo{year}{2015}{\natexlab{a}})
  \bibinfo{pages}{919--940}.
\bibitem[{Amsallem et~al.(2015{\natexlab{b}})Amsallem, Zahr, and
  Washabaugh}]{amsallem14:morepas}
\bibinfo{author}{D.~Amsallem}, \bibinfo{author}{M.~J. Zahr},
  \bibinfo{author}{K.~Washabaugh},
\newblock \bibinfo{title}{{Fast Local Reduced Basis Updates for the Efficient
  Reduction of Nonlinear Systems with Hyper-Reduction }},
\newblock \bibinfo{journal}{Special issue on Model Reduction of Parameterized
  Systems (MoRePaS), Advances in Computational Mathematics}
  (\bibinfo{year}{2015}{\natexlab{b}}) \bibinfo{pages}{1--34}.
\bibitem[{Wu and Hetmaniuk(2015)}]{wu15}
\bibinfo{author}{Y.~Wu}, \bibinfo{author}{U.~Hetmaniuk},
\newblock \bibinfo{title}{{Adaptive training of local reduced bases for
  unsteady incompressible Navier-Stokes flows}},
\newblock \bibinfo{journal}{International Journal for Numerical Methods in
  Engineering, published online}  (\bibinfo{year}{2015})
  \bibinfo{pages}{1--22}.
\bibitem[{Amsallem et~al.(2009)Amsallem, Cortial, Carlberg, and
  Farhat}]{amsallem09}
\bibinfo{author}{D.~Amsallem}, \bibinfo{author}{J.~Cortial},
  \bibinfo{author}{K.~Carlberg}, \bibinfo{author}{C.~Farhat},
\newblock \bibinfo{title}{{A method for interpolating on manifolds structural
  dynamics reduced-order models}},
\newblock \bibinfo{journal}{International Journal for Numerical Methods in
  Engineering} \bibinfo{volume}{80} (\bibinfo{year}{2009})
  \bibinfo{pages}{1241--1258}.
\bibitem[{Amsallem(2010)}]{amsallemthesis}
\bibinfo{author}{D.~Amsallem}, \bibinfo{title}{{Interpolation on Manifolds of
  CFD-Based Fluid and Finite Element-Based Structural Reduced-Order Models for
  On-Line Aeroelastic Predictions}}, Ph.D. thesis, Ph.D. Thesis, Stanford
  University, \bibinfo{year}{2010}.
\bibitem[{Degroote et~al.(2010)Degroote, Vierendeels, and Willcox}]{degroote10}
\bibinfo{author}{J.~Degroote}, \bibinfo{author}{J.~Vierendeels},
  \bibinfo{author}{K.~Willcox},
\newblock \bibinfo{title}{{Interpolation among reduced-order matrices to obtain
  parameterized models for design, optimization and probabilistic analysis}},
\newblock \bibinfo{journal}{International Journal for Numerical Methods in
  Fluids} \bibinfo{volume}{63} (\bibinfo{year}{2010})
  \bibinfo{pages}{207--230}.
\bibitem[{Panzer et~al.(2010)Panzer, Mohring, Eid, and Lohmann}]{panzer10}
\bibinfo{author}{H.~Panzer}, \bibinfo{author}{J.~Mohring},
  \bibinfo{author}{R.~Eid}, \bibinfo{author}{B.~Lohmann},
\newblock \bibinfo{title}{{Parametric Model Order Reduction by Matrix
  Interpolation}},
\newblock \bibinfo{journal}{at-Automatisierungstechnik} \bibinfo{volume}{58}
  (\bibinfo{year}{2010}) \bibinfo{pages}{475--484}.
\bibitem[{Amsallem and Farhat(2011)}]{amsallem11}
\bibinfo{author}{D.~Amsallem}, \bibinfo{author}{C.~Farhat},
\newblock \bibinfo{title}{{An online method for interpolating linear parametric
  reduced-order models}},
\newblock \bibinfo{journal}{SIAM Journal on Scientific Computing}
  \bibinfo{volume}{33} (\bibinfo{year}{2011}) \bibinfo{pages}{2169--2198}.
\bibitem[{Berkooz et~al.(1993)Berkooz, Holmes, and Lumley}]{berkooz93}
\bibinfo{author}{G.~Berkooz}, \bibinfo{author}{P.~Holmes},
  \bibinfo{author}{J.~L. Lumley},
\newblock \bibinfo{title}{{The proper orthogonal decomposition in the analysis
  of turbulent flows}},
\newblock \bibinfo{journal}{Annual Review of Fluid Mechanics}
  \bibinfo{volume}{25} (\bibinfo{year}{1993}) \bibinfo{pages}{539--575}.
\bibitem[{Grimme(1997)}]{grimme97}
\bibinfo{author}{E.~J. Grimme}, \bibinfo{title}{{Krylov projection methods for
  model reduction}}, Ph.D. thesis, Ph.D. Thesis, University of Illinois at
  Urbana Champaign, \bibinfo{year}{1997}.
\bibitem[{Hetmaniuk et~al.(2012)Hetmaniuk, Tezaur, and Farhat}]{hetmaniuk12}
\bibinfo{author}{U.~Hetmaniuk}, \bibinfo{author}{R.~Tezaur},
  \bibinfo{author}{C.~Farhat},
\newblock \bibinfo{title}{{Review and assessment of interpolatory model order
  reduction methods for frequency response structural dynamics and acoustics
  problems}},
\newblock \bibinfo{journal}{International Journal for Numerical Methods in
  Engineering} \bibinfo{volume}{90} (\bibinfo{year}{2012})
  \bibinfo{pages}{1636--1662}.
\bibitem[{Amsallem and Farhat(2014)}]{amsallem14:book}
\bibinfo{author}{D.~Amsallem}, \bibinfo{author}{C.~Farhat}, \bibinfo{title}{{On
  the Stability of Reduced-Order Linearized Computational Fluid Dynamics Models
  Based on POD and Galerkin Projection: Descriptor vs Non-Descriptor Forms}},
  volume~\bibinfo{volume}{9}, \bibinfo{publisher}{Reduced Order Methods for
  Modeling and Computational Reduction, MS{\&}A - Modeling, Simulation and
  Applications, Springer}, \bibinfo{year}{2014}.
\bibitem[{Fraikin et~al.(2008)Fraikin, Nesterov, and Dooren}]{fraikin08}
\bibinfo{author}{C.~Fraikin}, \bibinfo{author}{Y.~Nesterov},
  \bibinfo{author}{P.~V. Dooren},
\newblock \bibinfo{title}{{Optimizing the Coupling Between Two Isometric
  Projections of Matrices}},
\newblock \bibinfo{journal}{SIAM Journal on Matrix Analysis and Applications}
  \bibinfo{volume}{30} (\bibinfo{year}{2008}) \bibinfo{pages}{324--345}.
\bibitem[{Helmke and Barratt~Moore(1994)}]{helmkebook}
\bibinfo{author}{U.~Helmke}, \bibinfo{author}{J.~Barratt~Moore},
  \bibinfo{title}{{Optimization and Dynamical Systems}},
  \bibinfo{publisher}{Springer}, \bibinfo{year}{1994}.
\bibitem[{Zimmermann(2014)}]{zimmermann14}
\bibinfo{author}{R.~Zimmermann},
\newblock \bibinfo{title}{{A Locally Parametrized Reduced-Order Model for the
  Linear Frequency Domain Approach to Time-Accurate Computational Fluid
  Dynamics}},
\newblock \bibinfo{journal}{SIAM Journal on Scientific Computing}
  \bibinfo{volume}{36} (\bibinfo{year}{2014}) \bibinfo{pages}{B508--B537}.
\bibitem[{Choi et~al.(2015)Choi, Amsallem, and Farhat}]{choi15}
\bibinfo{author}{Y.~Choi}, \bibinfo{author}{D.~Amsallem},
  \bibinfo{author}{C.~Farhat},
\newblock \bibinfo{title}{{Gradient-Based Constrained Optimization Using a
  Database of Linear Reduced-Order Models}},
\newblock \bibinfo{journal}{submitted to Arxiv}  (\bibinfo{year}{2015})
  \bibinfo{pages}{1--21}.
\bibitem[{Grepl and Patera(2005)}]{grepl05}
\bibinfo{author}{M.~A. Grepl}, \bibinfo{author}{A.~T. Patera},
\newblock \bibinfo{title}{{A posteriori error bounds for reduced-basis
  approximations of parametrized parabolic partial differential equations}},
\newblock \bibinfo{journal}{ESAIM: Mathematical Modelling and Numerical
  Analysis} \bibinfo{volume}{39} (\bibinfo{year}{2005})
  \bibinfo{pages}{157--181}.
\bibitem[{Bui-Thanh et~al.(2008)Bui-Thanh, Willcox, and Ghattas}]{buithanh08}
\bibinfo{author}{T.~Bui-Thanh}, \bibinfo{author}{K.~Willcox},
  \bibinfo{author}{O.~Ghattas},
\newblock \bibinfo{title}{{Parametric reduced-order models for probabilistic
  analysis of unsteady aerodynamic applications}},
\newblock \bibinfo{journal}{AIAA Journal} \bibinfo{volume}{46}
  (\bibinfo{year}{2008}) \bibinfo{pages}{2520--2529}.
\bibitem[{Amsallem and Hetmaniuk(2015)}]{amsallem14:expintegrators}
\bibinfo{author}{D.~Amsallem}, \bibinfo{author}{U.~Hetmaniuk},
\newblock \bibinfo{title}{{A posteriori error estimators for linear reduced
  order models using Krylov-based integrators}},
\newblock \bibinfo{journal}{International Journal for Numerical Methods in
  Engineering} \bibinfo{volume}{102} (\bibinfo{year}{2015})
  \bibinfo{pages}{1238--1261}.
\bibitem[{Colton and Kress(2013)}]{coltonbook}
\bibinfo{author}{D.~Colton}, \bibinfo{author}{R.~Kress},
  \bibinfo{title}{{Inverse acoustic and electromagnetic scattering theory}},
  \bibinfo{publisher}{Springer}, \bibinfo{year}{2013}.
\bibitem[{Berenger(1994)}]{berenger94}
\bibinfo{author}{J.~P. Berenger},
\newblock \bibinfo{title}{{A perfectly matched layer for the absorption of
  electromagnetic waves}},
\newblock \bibinfo{journal}{Journal of Computational Physics}
  \bibinfo{volume}{114} (\bibinfo{year}{1994}) \bibinfo{pages}{185--200}.
\bibitem[{Chiu and Farhat(2009)}]{chiu09}
\bibinfo{author}{E.~K.-y. Chiu}, \bibinfo{author}{C.~Farhat},
\newblock \bibinfo{title}{{Effects of fuel slosh on flutter prediction}},
\newblock \bibinfo{journal}{AIAA 2009-2682, 50th AIAA/ASME/ASCE/AHS/ASC
  Structures, Structural Dynamics, and Materials Conference}
  (\bibinfo{year}{2009}).
\bibitem[{Farhat et~al.(2013)Farhat, Chiu, Amsallem, Schott{\'e}, and
  Ohayon}]{farhat13}
\bibinfo{author}{C.~Farhat}, \bibinfo{author}{E.~K.-y. Chiu},
  \bibinfo{author}{D.~Amsallem}, \bibinfo{author}{J.-S. Schott{\'e}},
  \bibinfo{author}{R.~Ohayon},
\newblock \bibinfo{title}{{Modeling of Fuel Sloshing and its Physical Effects
  on Flutter}},
\newblock \bibinfo{journal}{AIAA Journal} \bibinfo{volume}{51}
  (\bibinfo{year}{2013}) \bibinfo{pages}{2252--2265}.
\bibitem[{Yates(1987)}]{yates87}
\bibinfo{author}{E.~C. Yates}, \bibinfo{title}{{AGARD Standard Aeroelastic
  Configurations For Dynamic Response - 1 - Wing 445.6}},
  \bibinfo{howpublished}{NASA}, \bibinfo{year}{1987}.

\end{thebibliography}

\end{document}